\documentclass[11pt, dina4,  amssymb]{article}
\usepackage{amssymb, amsmath,amsthm}
\usepackage{graphicx}
\usepackage{diagbox}
\usepackage{multirow}
\usepackage{epstopdf}
\usepackage{color}
\usepackage[margin=1in]{geometry}
\usepackage[normalem]{ulem}
\usepackage{xcolor}

\newcommand\redsout{\bgroup\markoverwith{\textcolor{red}{\rule[0.5ex]{2pt}{0.4pt}}}\ULon}

\newcommand{\bxi}{\textit{{\boldmath$\xi$}}}

\newcommand{\p}{\partial}
\newcommand{\og}{\omega}
\newcommand{\Og}{\Omega}
\newcommand{\fl}[2]{\frac{#1}{#2}}
\newcommand{\dt}{\delta}
\newcommand{\gm}{\gamma}
\newcommand{\nn}{\nonumber}
\newcommand{\ap}{\alpha}

\newcommand{\veps}{\varepsilon}

\newcommand{\Dt}{\Delta}

\newcommand{\be}{\begin{equation}}
\newcommand{\ee}{\end{equation}}
\newcommand{\ba}{\begin{array}}
\newcommand{\ea}{\end{array}}
\newcommand{\bea}{\begin{eqnarray}}
\newcommand{\eea}{\end{eqnarray}}
\newcommand{\beas}{\begin{eqnarray*}}
\newcommand{\eeas}{\end{eqnarray*}}
\newtheorem{remark}{Remark}[section]
\newtheorem{lemma}{Lemma}[section]
\newtheorem{theorem}{Theorem}[section]

\newtheorem{definition}{Definition}[section]

\newtheorem{property}{Property}[section]
\newcommand{\bx}{{\bf x} }

\definecolor{ForestGreen}{rgb}{0.0, 0.5, 0.0}

\title{Accurate numerical methods for two and three dimensional integral fractional Laplacian with applications}
\author{
Siwei Duo\thanks{ Department of Mathematics,  University of South Carolina, Columbia, SC 29208 (Email: duo@mailbox.sc.edu)} \  \  and 
Yanzhi Zhang\thanks{Department of Mathematics and Statistics, Missouri University of Science and Technology, Rolla, MO 65409-0020 (Email: zhangyanz@mst.edu) }
}
\date{}

\begin{document}
\maketitle

% =======================================================================
% Abstract;
% =======================================================================
\begin{abstract}
In this paper, we propose accurate and efficient finite difference methods to discretize the two- and three-dimensional fractional Laplacian $(-\Dt)^{\fl{\ap}{2}}$ ($0 < \ap < 2$) in hypersingular integral form. 
The proposed finite difference methods provide a fractional analogue of the central difference schemes to the fractional Laplacian, and as $\ap \to 2^-$, they collapse to the central difference schemes of the classical Laplace operator $-\Dt$.  
We prove that our methods are consistent if $u \in C^{\lfloor\ap\rfloor,\,\ap-\lfloor\ap\rfloor+\veps}({\mathbb R}^d)$, and the local truncation error is ${\mathcal O}(h^\veps)$,  with $\veps > 0$ a small constant and $\lfloor \cdot \rfloor$ denoting the floor function. 
If $u \in C^{2+\lfloor\ap\rfloor,\,\ap-\lfloor\ap\rfloor+\veps}({\mathbb R}^d)$, they can achieve the second order of accuracy for any $\ap \in (0, 2)$. 
These results hold for any dimension $d \ge 1$ and thus improve the existing error estimates for the finite difference method of the one-dimensional fractional Laplacian. 
Extensive numerical experiments are provided and confirm our analytical results.  
We then apply our method to solve the fractional Poisson problems and the fractional Allen--Cahn equations. 
Numerical simulations suggest that to achieve the second order of accuracy, the solution of the fractional Poisson problem should {\it at most} satisfy  $u \in C^{1,1}({\mathbb R}^d)$. 
One merit of our methods is that they yield a multilevel Toeplitz stiffness matrix, an appealing property for the development of fast algorithms via the fast Fourier transform (FFT). 
Our studies of the two- and three-dimensional fractional Allen--Cahn equations demonstrate the efficiency of our methods in solving the high-dimensional fractional problems. 
\end{abstract}


% =======================================================================
% Introduction;
% =======================================================================
\section{Introduction}
\label{section1}

%%% %%%% %%%% %%%%
In classical partial differential equation (PDE) models, diffusion is described by the classical Laplace operator $\Dt = (\p_{xx} + \p_{yy} + \p_{zz})$, characterizing the transport mechanics due to Brownian motion. 
Recently, it has been suggested that many complex (e.g., biological and chemical) systems are instead characterized by L{\' e}vy motion, rather than Brownian motion; see \cite{del-Castillo-Negrete2003, Hanert2011, Metzler2000, Cusimano2015} and references therein. 
Hence, the classical PDE models fail to properly describe the phenomena in these  systems.  
To circumvent such issues, the fractional models were proposed, where the classical Laplace operator is replaced by the fractional Laplacian $(-\Dt)^{\fl{\ap}{2}}$ \cite{Metzler2000, del-Castillo-Negrete2003}. 
In contrast to the classical diffusion models, the fractional models possess significant advantages for describing problems with long-range interactions, enabling one to describe the power law invasion profiles that have been observed in many applications \cite{Clark2001,Paradis2002,Tackenberg2003}. 
However, the nonlocality of the fractional Laplacian introduces considerable challenges in both mathematical analysis and numerical simulations.

%%%%%%%% Problem setup %%%%%%%%
In this study, we focus on the finite difference discretization of the integral fractional Laplacian $(-\Delta)^{\fl{\alpha}{2}}$ (for $0 < \ap < 2$). 
Probabilistically, the fractional Laplacian represents the infinitesimal generator of a symmetric $\alpha$-stable L\'{e}vy process. 
It is defined in terms of the hypersingular integral \cite{Landkof1972, Samko}:
\bea\label{fL-nD}
(-\Delta)^{\fl{\ap}{2}}u({\bf x}) 
= c_{d,\ap} \ {\rm P.V.}\int _{\mathbb{R}^d} \frac{u({\bf x})-u({\bf x}')}{|{\bf x}-{\bf x}^\prime|^{d+\alpha}}\,d{\bf x'}, \qquad \mbox{for \ $\alpha \in (0,2)$},
\eea
where P.V. stands for the Cauchy principal value, and $|\bx - \bx'|$ denotes the Euclidean distance between points $\bx$ and $\bx'$. 
The normalization constant $c_{d,\ap}$ is defined as
\beas\label{DefC1ap}
 c_{d,\ap} = \fl{2^{\ap-1}\ap\,\Gamma\big(({d+\ap})/{2}\big)}{\sqrt{\pi^d}\,\Gamma\big(\displaystyle 1-{\ap}/{2}\big)}, \qquad  \mbox{for} \  \  \ d = 1, 2, 3
\eeas
with $\Gamma(\cdot)$ denoting the Gamma function.
In the literature, the fractional Laplacian is also defined via a pseudo-differential operator with symbol $|{\bf k}|^\alpha$  \cite{Landkof1972, Samko}, i.e., 
\begin{equation}
\label{pseudo}
(-\Delta)^{\fl{\alpha}{2}}u({\bx}) = \mathcal{F}^{-1}\big[|{\bf k}|^\alpha \mathcal{F}[u]\big], \qquad \mbox{for} \ \ \ap > 0,
\end{equation}
where $\mathcal{F}$ represents the Fourier transform, and $\mathcal{F}^{-1}$ denotes its inverse.  
Over the entire space $\mathbb{R}^d$, the fractional Laplacian (\ref{fL-nD}) is equivalent to the pseudo-differential operator (\ref{pseudo}) \cite{DiNezza2012, Samko, Kwasnicki2017, DuoWangZhang}. 
As pointed out in \cite{Izsak2017}, the $d$-dimensional (for any $d\ge1$) fractional Laplacian is rotational invariant, which is a crucial property in modeling isotropic anomalous diffusion in many applications \cite{Kullberg2013}. 
The rotational invariance distinguishes the fractional Laplacian from many other fractional derivatives when $d\geq 2$ \cite{Izsak2017}. 

%%%% %%%% %%%% %%%%
Due to its nonlocality and strong singularity, numerical methods for the integral fractional Laplacian (\ref{fL-nD}) still remain scant. 
In \cite{Acosta2017, Acosta2017-2D, Bonito0017, Bonito2018, Ainsworth0017},  various finite element  methods (FEM) are proposed to solve problems with the integral fractional Laplacian (\ref{fL-nD}), and numerical results are presented for one-, two-, and three-dimensional cases.
In \cite{Zhang0018}, a spectral Galerkin method is presented for the one-dimensional fractional reaction-diffusion equations. 
In contrast to the progress of FEM, the development of  finite difference methods to the integral fractional Laplacian still falls behind. 
So far, although several finite difference methods exist in the literature \cite{Duo2018, Huang2014, Gao2014}, they are all limited to one-dimensional (i.e., $d = 1$) cases. 
To the best of our knowledge, finite difference methods for the high-dimensional (i.e., $d > 1$) fractional Laplacian (\ref{fL-nD}) are still missing in the literature. 
%Moreover, no numerical reports can be found for the problems with the three-dimensional fractional Laplacian, due to the  discretization difficulty and huge computational costs and storage memory requirements. 

%%%% %%%% %%%% %%%% 
In this paper, we propose accurate and efficient finite difference methods to discretize the two- and three-dimensional integral fractional Laplacian (\ref{fL-nD}). 
Our methods  fill the gap in the literature on the finite difference methods for the fractional Laplacian $(-\Dt)^{\fl{\ap}{2}}$.
More importantly,  they yield a symmetric {multilevel Toeplitz} stiffness matrix and thus enable one to develop fast algorithms via the fast Fourier transform (FFT) for their efficient implementation. 
The proposed finite difference methods provide a fractional analogue of the central difference schemes to the fractional Laplacian, and as $\ap \to 2^-$, they collapse to the five-point (for $d = 2$) or {seven}-point (for $d = 3$) finite difference schemes to the classical Laplace operator $-\Dt$.  
We carry out both numerical analysis and simulations to understand the performance of our   methods. 
It shows that our methods are consistent if $u \in C^{\lfloor\ap\rfloor,\,\ap-\lfloor\ap\rfloor+\veps}({\mathbb R}^d)$, and the local truncation error is ${\mathcal O}(h^\veps)$,  with $\veps > 0$ a small constant and $\lfloor \cdot \rfloor$ denoting the floor function. 
If $u \in C^{2+\lfloor\ap\rfloor,\,\ap-\lfloor\ap\rfloor+\veps}({\mathbb R}^d)$, our methods can achieve the second order of accuracy for any $\ap \in (0, 2)$.
Extensive numerical experiments confirm our analytical results. 
Note that our results improve the error estimates in \cite{Duo2018} and provide  much sharper consistent conditions for the finite difference methods of the fractional Laplacian. 
We then apply our methods to solve the fractional Poisson problems and the fractional Allen--Cahn equations. 
Our numerical studies suggest that our finite difference methods have a second order of accuracy in solving the fractional Poisson equation,  if the solution $u \in C^{1, 1}({\mathbb R}^d)$. 
To demonstrate their effectiveness in solving high-dimensional problems, we implement our methods via FFT-based fast algorithms to solve both two- and three-dimensional fractional Allen--Cahn equations.
It shows that at each time step the computational complexity is ${\mathcal O}(M\log M)$, and memory storage requirement is ${\mathcal O}(M)$, with $M$ the total number of spatial unknowns. 

%%%%%%%% Paper organization %%%%%%%%
The paper is organized as follows. 
In Section \ref{section2}, we propose a finite difference method for the two-dimensional fractional Laplacian (\ref{fL-nD}), and the detailed error estimates are provided in Section \ref{section3}.  
In Section \ref{section6},  we generalize our method and analysis in Sections \ref{section2}--\ref{section3} to the three dimensions.  
Numerical examples are presented in Section \ref{section7} to test the performance of our methods and study the fractional problems, including the fractional Poisson equation and the fractional Allen--Cahn equation. 
Finally, we draw conclusions in Section \ref{section8}. 

% ======================================================================
\section{Finite difference method in two dimensions}
\label{section2}
\setcounter{equation}{0}

%%%% %%%% %%%% %%%% 
So far, numerical methods for the integral fractional Laplacian (\ref{fL-nD}) still remain very limited, especially in high dimensions (i.e., $d > 1$), and the main challenges come from its strong singularity. 
Recently, several finite difference methods have been proposed to discretize the one-dimensional (1D) fractional Laplacian; see \cite{Duo2018, Huang2014, Gao2014} and references therein.  
Among them, the method in \cite{Duo2018} is the current state-of-the-art finite difference method for the 1D fractional Laplacian -- it can achieve the second order of accuracy uniformly for any $\ap\in(0, 2)$. 
However, the finite  difference method for the high-dimensional fractional Laplacian (\ref{fL-nD}) is still missing in the literature. 
In this section, we present a finite difference method for the two-dimensional (2D) fractional Laplacian, and its generalization to the three-dimensional  (3D) cases can be found in Section \ref{section6}. 

%%%%%%%%%%%
We consider the fractional Laplacian (\ref{fL-nD}) on the bounded domain $\Og = (a_x, b_x)\times (a_y,b_y)$ with the extended homogeneous Dirichlet boundary conditions on $\Og^c = {\mathbb R}^2\backslash\Og$. 
Introducing two new variables $\xi = |x-x'|$ and  $\eta = |y-y'|$ and denoting the vector $\bxi = (\xi,\eta)$, we can reformulate the 2D fractional Laplacian (\ref{fL-nD}) as:  
\begin{eqnarray}\label{fL-0}
(-\Delta)^{\fl{\ap}{2}}u({\bf x})=-c_{2,\ap}\int _{0}^{\infty}\int _{0}^{\infty} \fl{1}{|\bxi|^{2+\ap}}\bigg(\sum_{m,n=0, 1}
u\big(x+(-1)^{m}\xi,\, y+(-1)^{n}\eta\big)-4 u({\bf x})\bigg)d\xi d\eta.
\end{eqnarray}
Here, we introduce a splitting parameter $\gm \in (\ap, 2]$, and define the function 
\begin{eqnarray}\label{psi-2D}
\psi_{\gamma}({\bf x},\textit{{\boldmath$\xi$}}) = \fl{1}{|\bxi|^\gm}
\Big(\sum_{m,  n=0, 1} u\big(x+(-1)^{m}\xi, \, y+(-1)^{n}\eta\big)-4 u(\bx) \Big).
\end{eqnarray}
Then, the fractional Laplacian in (\ref{fL-0}) can be further written as the weighted integral of the central difference quotient  $\psi_\gm$, i.e., 
\bea\label{fL}
(-\Delta)^{\fl{\ap}{2}}u({\bf x}) = -c_{2,\ap}\int _{0}^{\infty} \int _{0}^{\infty}
\psi _{\gamma}({\bf x},\textit{{\boldmath$\xi$}}) \og_{\gm}(\bxi)\,d\xi d\eta, \qquad\mbox{with} \quad
\og_\gm(\textit{{\boldmath$\xi$}}) = |\textit{{\boldmath$\xi$}}|^{\gm-(2+\ap)},
\eea
where $\og_\gm$ can be viewed as the weight function. 
%%%% %%%% %%%% %%%%
The key novelty of our approach is to reformulate the fractional Laplacian (\ref{fL-nD}) as a weighted integral of the central difference quotient in (\ref{fL}), that is, splitting the strong-singular kernel function $|\bxi|^{-(2+\ap)}$ into two weak-singular parts, i.e., $|\bxi|^\gm$ and $|\bxi|^{\gm-(2+\ap)}$. 
This idea of splitting the kernel function was introduced in \cite{Duo2018, Duo2015} for the 1D fractional Laplacian, and it has been applied to solve the fractional Schr\"{o}dinger equation in an infinite potential well \cite{Duo2015}. 
Note that the splitting parameter $\gm$ plays a crucial role in determining the accuracy of our method, which will be further  discussed in Sections \ref{section3} and \ref{section7}. 

%%%% %%%% %%%% %%%%
Let  $L = \max \{b_x-a_x, \, b_y - a_y\}$. 
Denote $D_1 = (0, L)^2$ and $D_2 = {\mathbb R}^{2}_+\backslash  D_1 = \{(\xi, \eta)\mid \xi \ge L \ \mbox{or}\  \eta \ge L\}$. 
We then divide the integral in (\ref{fL}) into two parts:  
\begin{eqnarray}\label{fL1-2D}
(-\Delta)^{\fl{\alpha}{2}}u({\bf x})=-c_{2,\alpha} \bigg(
\int_{D_1}\psi_{\gamma}({\bf x},\, \textit{{\boldmath$\xi$}}) \og_{\gm}(\bxi)\,d\xi d\eta 
+ \int_{D_2}\psi_{\gamma}({\bf x},\, \textit{{\boldmath$\xi$}}) \og_{\gm}(\bxi)\,d\xi d\eta\bigg).
\end{eqnarray} 
Due to the extended homogeneous Dirichlet boundary conditions over $\Og^c$, we have that for ${\bf x}\in \Omega$, and $\xi \ge L$ or $\eta \ge L$,   the point $(x + (-1)^m\xi,\, y + (-1)^n\eta)\in \Omega^c$, for $m, n = 0, 1$.
Thus,  we obtain  for any $\bx\in \Og$ and $\bxi \in D_2$, the function $u(x + (-1)^m\xi, y + (-1)^n\eta) \equiv 0$ and immediately $\psi_\gm(\bx, \textit{{\boldmath$\xi$}})  = -4u(\bx)|\textit{{\boldmath$\xi$}}|^{-\gm}$.  
Then, the integral over $D_2$ is simplified as: 
\begin{eqnarray}\label{ILinf-2D}
\int_{D_2}\psi_{\gamma}({\bf x},\,\textit{{\boldmath$\xi$}})\,\og_{\gm}(\bxi)\,d\xi d\eta 
= -4 u({\bf x}) \int_{D_2} |\textit{{\boldmath$\xi$}}|^{-(2+\alpha)}\,d\xi d\eta.
\eea
It is easy to see that the only possible error of (\ref{ILinf-2D}) is from evaluating the integral of $|\textit{{\boldmath$\xi$}}|^{-(2+\alpha)}$, which can be  easily controlled.

%%%% %%%% %%%% %%%%
We now move to approximate the first integral of (\ref{fL1-2D}).  
As discussed previously, the main difficulty is from the strong singular kernel, and the traditional quadrature rules fail to provide a satisfactory approximation. 
To achieve high accuracy,  we propose a weighted trapezoidal method, i.e., retaining part of the singularity in the integral. 
%%%% %%%% %%%% %%%%
Choose an integer $N > 0$, and define the mesh size $h = L/N$.
Denote grid points $\xi_i = ih$ and $\eta_{j} = jh$, for $0\leq i, j \leq N$. 
For notational simplicity, we denote $\bxi_{ij} = (\xi_i,\,\eta_j)$ and then $|\bxi_{ij}|=\sqrt{\xi_{i}^2+\eta_{j}^2}$.   
Additionally,  we define the element 
\beas
{I}_{ij}:=[ih, \, (i+1)h]\times [jh, (j+1)h], \qquad  \mbox{for} \ \  0\leq i,\,j \leq N-1.
\eeas
It is easy to see that $D_1 = \cup_{i, j = 0}^{N-1}\, I_{ij}$.  
Hence,  the  integral over $D_1$ can be formulated as: 
\bea\label{I0L-2D-WT}
\int _{D_1} \psi_{\gamma}({\bf x},\bxi) \og_{\gm}(\bxi)\,d\xi d\eta =\sum _{i, j = 0}^{N-1} \int _{{I}_{ij}} \psi_{\gamma}({\bf x},\bxi) \og_{\gm}(\bxi)\,d\xi d\eta.
\eea
Next, we focus on the approximation to the integral over each element $I_{ij}$. 

%%%% %%%% %%%% %%%%
For $i \neq 0$ or $j \neq 0$,  the weighted trapezoidal rule is applied,  and we obtain the approximation: 
\begin{eqnarray}\label{I-2D}
\int _{I_{ij}}\psi_{\gamma}({\bf x},\bxi) \og_{\gm}(\bxi)\,d\xi d\eta
\approx \frac{1}{4}\bigg(\sum_{m, n = 0, 1} \psi_\gm\big(\bx, \bxi_{(i+m)(j+n)}\big)\bigg)\int _{I_{ij}} \og_{\gm}(\bxi)\,d\xi d\eta. %\quad \ \ \mbox{if} \ \ i+j\neq0.
\end{eqnarray}
While $i = j = 0$, the approximation of the integral over $I_{00}$ is not as straightforward as that in (\ref{I-2D}). 
Using the weighted trapezoidal rule, we get
\bea\label{I00}
\int _{I_{00}}
\psi_{\gamma}({\bf x},\bxi) \og_{\gm}(\bxi)\,d\xi d\eta 
\approx\frac{1}{4}\bigg(\lim_{\bxi \rightarrow {\bf 0}}\psi_{\gamma}\big({\bf x},\,{\bxi}\big) + \sum_{\substack{m, n = 0, 1\\ m+n\neq 0}} \psi_\gm\big(\bx, \bxi_{(i+m)(j+n)}\big)\bigg) 
\int _{I_{00}} \og_{\gm}(\bxi)\,d\xi d\eta.
\eea
Assuming the limit in (\ref{I00}) exists, then it depends on the splitting parameter $\gm$.  
We will divide our discussion into two cases: $\gm \in (\ap, 2)$ and $\gm = 2$. 
If $\gm = 2$,  the limit can be approximated by: 
\bea
\label{lim2}
\lim_{\bxi \rightarrow {\bf 0}}\psi_{2}({\bf x},\,{\bxi}) \approx \psi_{2}({\bf x},\,{\bxi}_{10}) + \psi_{2}({\bf x},\,{\bxi}_{01})- \psi_{2}({\bf x},\,{\bxi}_{11}).
\eea
If  $\gm \in (\ap, 2)$, we have
\bea\label{lim1}
\lim_{\bxi \rightarrow {\bf 0}}\psi_{\gamma}({\bf x},\,{\bxi}) = \lim_{\bxi \rightarrow {\bf 0}}\big(\psi_{2}({\bf x},\, {\bxi}) |\bxi|^{2-\gamma}\big) = 0. 
\eea
Substituting (\ref{lim2})--(\ref{lim1}) into (\ref{I00}), we obtain the approximation of the integral over $I_{00}$ as:
\bea\label{I00-final}
\int _{I_{00}}
\psi_{\gamma}({\bf x},\bxi) \og_\gm(\bxi)\,d\xi d\eta \approx 
\fl{1}{4}\bigg(\sum_{\substack{m, n = 0,1\\ m+n\neq 0}} c_{mn}^\gm\, \psi_\gm\big(\bx, \bxi_{mn}\big)\bigg) \int _{I_{00}} \og_\gm(\bxi)\,d\xi d\eta,\ \  
\eea
where the coefficient $c_{mn}^\gm = 1$ for $\gm \in (\ap, 2)$, while $c_{10}^\gm = c_{01}^\gm = 2$ and $c_{11}^\gm = 0$ for $\gm = 2$. 
%%%% %%%% %%%%% %%%%
Let $T_{ij}$ denote the collection of all elements that have  $\bxi_{ij}$ as a vertex,  i.e.,
\beas
{T}_{ij} = \big({I}_{(i-1)(j-1)} \cup {I}_{(i-1)j}  \cup {I}_{i(j-1)} \cup {I}_{ij} \big) \cap D_1,  \quad \ \ \mbox{for} \   \ 0\leq i,\,j \leq N-1.
\eeas
%Especially,  $T_{00} = I_{00}$.
Then, combining (\ref{fL1-2D})--(\ref{I-2D}) and (\ref{I00-final}) and reorganizing the terms, we obtain the approximation to the 2D fractional Laplacian (\ref{fL-nD}) as follows:
\bea\label{I0L-2D}  
&&(-\Delta)_{h,\gamma}^{\fl{\alpha}{2}} u(\bx) = 
-\frac{c_{2,\ap}}{4}\bigg(\left\lfloor\fl{\gm}{2}\right\rfloor\Big(\psi_{\gamma}\big({\bf x}, {\bxi}_{10}\big)+\psi_{\gamma}\big({\bf x}, {\bxi}_{01}\big) -\psi_{\gamma}\big({\bf x},{\bxi}_{11}\big)\Big)\int _{I_{00}}\og_\gm(\bxi)\,d\xi d\eta \qquad\qquad\nn\\
&&\hspace{1.6cm}+ {\sum_{\substack{i, j = 0 \\ i + j \neq 0}}^{N-1}}
\psi_{\gamma}\big({\bf x}, \bxi_{ij}\big)\int _{T_{ij}}\og_\gm(\bxi)\,d\xi d\eta - 16u(\bx)\int_{D_2} |\bxi|^{-(2+\ap)} d\xi d\eta\bigg),  \quad\ \ \bx\in \Og.
\eea

%%%% %%%% %%%% %%%%
Without loss of generality, we assume that $N_x = N$, and choose $N_y$  as the smaller integer such that  $a_y + N_yh \geq b_y$. 
Define the grid points $x_i = a_x + ih$ for $0 \le i \le N_x$,  and $y_j = a_y + jh$ for $0 \le j \le N_y$. 
Let $u_{ij}$ be the numerical approximation of $u(x_i,y_j)$. 
Noticing the definition of  $\psi_\gm$ in (\ref{psi-2D}), we obtain the finite difference discretization to the 2D fractional Laplacian $(-\Dt)^{\fl{\ap}{2}}$ as: 
 \bea\label{fLh-2D}
&&(-\Delta)_{h,\gamma}^{\fl{\alpha}{2}}u_{ij}  = -c_{2,\alpha}
\bigg[a_{00} u_{ij}+\sum_{m=0}^{i-1}\bigg(\sum_{\substack{n=0 \\ m+n\neq 0}}^{j-1}a_{mn}u_{(i-m)(j-n)} + \sum_{n=1}^{N_y-1-j} a_{mn} u_{(i-m)(j+n)} \bigg)\qquad\qquad  \nn\\
&&\hspace{2.5cm}  +\sum_{m=0}^{N_x-1-i} \bigg(\sum_{\substack{n = 0 \\ m+n\neq 0}}^{j-1}a_{mn}u_{(i+m)(j-n)}+\sum_{n=1}^{N_y-1-j}a_{mn}u_{(i+m)(j+n)}\bigg)\bigg],
\eea
for $1 \le i \le N_x-1$ and $1 \le j \le N_y-1$. 
The scheme (\ref{fLh-2D}) shows that the discretized fractional Laplacian at point $(x_i, y_j)$ depends on all points in the domain $\Og$,  consistent with the nonlocal characteristic of the fractional Laplacian. 
The coefficient $a_{mn}$ depends on the choice of the splitting parameter $\gm$. 
For $m, n \ge 0$ but $m + n > 0$,  there is 
\bea\label{amn}
a_{mn} =  
\frac{2^{\sigma(m, n)}}{4|\bxi_{mn}|^{\gamma}}
\bigg(\int_{T_{mn}} |\bxi|^{\gamma-(2+\alpha)}\,d\xi d\eta
+ \bar{c}_{mn}\left\lfloor\fl{\gm}{2}\right\rfloor\int_{I_{00}} |\bxi|^{\gamma-(2+\alpha)}\,d\xi d\eta\bigg), \qquad \eea
where $\sigma(m, n)$ denotes the number of zeros of $m$ and $n$, and the constant $\bar{c}_{01} = \bar{c}_{10} = -\bar{c}_{11} = 1$, and  $\bar{c}_{mn} \equiv 0$ for other $m, n$. 
For $m = n = 0$, the coefficient
\bea\label{a00}
a_{00} = -2\sum_{m=1}^{N} \big(a_{m0} + a_{0m}\big)
-4\sum_{m, n =1}^{N} a_{mn}
-4 \int_{D_2} |\bxi|^{-(2+\ap)}\,d\xi d\eta.
\eea

%%%% %%%% %%%% %%%% Remark 
\begin{remark}
In the limit of $\ap \to 2^-$, our finite difference scheme in (\ref{fLh-2D}) with $\gm = 2$ reduces to  the five-point finite difference scheme of the Laplace operator $-\Dt = -(\p_{xx} + \p_{yy})$, consistent with the limit behavior of the fractional Laplacian $(-\Dt)^{\fl{\ap}{2}}$.
In fact,  if $\gm = 2$ we obtain
\beas
\lim_{\ap \to2^-} \big(c_{2,\ap}  a_{10}\big) = \lim_{\ap \to2^-} \big(c_{2,\ap}  a_{01}\big)= \fl{1}{h^2}; \quad \ 
\lim_{\ap \to2^-} \big(c_{2,\ap}  a_{mn}\big) = 0, \ \ \mbox{for} \ \, m+n >1, 
\eeas
due to the following properties:
\beas
\lim_{\ap \to2^-} \bigg(c_{2,\ap}\int_{I_{00} }|\bxi|^{-\ap}d\bxi\bigg) = 1,  \qquad
{\lim_{\ap \to2^-} c_{2,\ap} = 0}.
\eeas
Immediately, the relation in (\ref{a00}) implies
\beas
\lim_{\ap \to2^-} \big(c_{2,\ap}  a_{00}\big) = -\fl{4}{h^2}. 
\eeas
%Additionally, even though the splitting parameter $\gm \in (\ap, 2]$,  our extensive studies show that the optimal choice is $\gm = 2$, which usually leads to the smallest numerical errors.
\end{remark}

%%%% %%%% %%%% %%%%
We can write the scheme (\ref{fLh-2D}) into matrix-vector form. 
 For $1\le j \le N_y-1$, denote the vector ${\bf u}_{x, j} = (u_{1, j}, u_{2, j},\ldots,u_{N_x-1, j})$, and let the block vector $
{\bf u} = ({\bf u}_{x,1}, {\bf u}_{x,2},\,\ldots,{\bf u}_{x,\,N_y-1})^T$. 
Then the matrix-vector form of the  scheme (\ref{fLh-2D}) is given by 
\bea\label{product}
(-\Delta)^{\fl{\alpha}{2}}_{h,\gamma}{\bf u} = A_2{\bf u}, 
\eea
where the matrix $A_2$ is a symmetric block Toeplitz matrix, defined as
\begin{eqnarray}\label{A-2D}
{{A_2}}= 
\left(
\begin{array}{cccccc}
A_{x,0} & A_{x,1} & \ldots &  A_{x,N_y-3} & A_{x,N_y-2}  \\
A_{x,1}& A_{x,0} & A_{x,1}  &  \cdots & A_{x,N_y-3}  \\
\vdots & \ddots  & \ddots & \ddots & \vdots \\
A_{x,N_y-3} &  \ldots  & A_{x,1} & A_{x,0} & A_{x,1} \\
A_{x,N_y-2} & A_{x_1,N_y-3}  & \ldots & A_{x,1} & A_{x,0}
\end{array}
\right) _{M\times M}
\end{eqnarray}
with $M = (N_x-1)(N_y-1)$ being the total number of unknowns, and each block $A_{x, j}$ (for $0 \le j \le N_y-2$) is a symmetric Toeplitz matrix, defined as 
\begin{eqnarray}\label{Aj}
{A}_{x,j} = 
-c_{2, \ap}\left(
\begin{array}{cccccc}
a_{0j} & a_{1j} & \ldots &  a_{(N_x-3)j} & a_{(N_x-2)j}  \\
a_{1j} & a_{0j} & a_{1j}  &  \cdots & a_{(N_x-3)j}  \\
\vdots & \ddots  & \ddots & \ddots & \vdots \\
a_{(N_x-3)j} &  \ldots  & a_{1j} & a_{0j} & a_{1j} \\
a_{(N_x-2)j} & a_{(N_x-3)j}  & \ldots & a_{1j} & a_{0j}
\end{array}
\right) _{(N_x-1)\times (N_x-1)}.\nn
\end{eqnarray}
It is easy to verify that the matrix $A_2$ is positive definite.  
In contrast to the differentiation matrix of the classical Laplacian,  the matrix $A_2$ in (\ref{A-2D}) is a large dense matrix, which usually causes considerable challenges not only for storing the matrix but also for computing matrix-vector products. 
However, thanks to the block-Toeplitz-Toeplitz-block structure of $A_2$, one can develop a FFT-based fast algorithm to efficiently compute matrix-vector multiplication in (\ref{product}), which greatly save not only the computational costs but also  storage requirements \cite{Du2015}.

% ======================================================================
\section{Error analysis}
\label{section3}
\setcounter{equation}{0}
 
%%%% %%%% %%%% %%%%
In this section, we provide the error estimates for our finite difference method in discretizing the 2D fractional Laplacian. 
The main technique used in our proof is an extension of the weighted Montgomery identity (see Lemma \ref{lemma2-2D}). 
The Montgomery identity is the framework of developing many classical inequalities, such as the Ostrowski, Chebyshev,  and Gr{\"u}ss type inequalities. 
%\cite{Dragomir1997, Dragomir1998, Cerone2003}. 
As an extension, the weighted Montgomery identity, first introduced in \cite{Pecaric1980, Khan2013}, plays an important role in the study of weighted integrals. 
Here, we will begin with introducing an extension function of the generalized Peano kernel. 

%%%% %%%% %%%% %%%%% Definition -- K
\begin{definition}\label{K-2D} 
Let $w: [a, b]\times[c,\,d] \to {\mathbb R}$ be an integrable function. For $m, n \in {\mathbb N}^0$, define 
\begin{eqnarray*}
\Theta_{[a,\,b]\times[c,\,d]}^{(m, n)}(x,y) =
\sum_{(s,  t) \in S} \int_{t}^{y} \int_{s}^{x}
w(\xi,\eta)\,\frac{(x-{\xi})^{m}(y-{\eta})^{n}}{m!\ n!}\,d{\xi}d{\eta}, \quad \ \, (x, y) \in [a,\,b]\times[c,\,d],
\end{eqnarray*}
where the set $S = \big\{(a,c), (a,d), (b,c), (b,d)\big\}$.
\end{definition}   

The function $\Theta$ has the following properties: 

%%%% %%%% %%%% %%%% Property -- K
\begin{property}\label{K-Property} Let  $m, n \in {\mathbb N}^0$, and $(x, y) \in [a,\,b]\times[c,\,d]$.
\vskip 6pt 

\noindent (i) \ If $w(x, y) = w(y, x)$, then 
\begin{eqnarray}
\Theta_{[a,\,b]\times[c,\,d]}^{(m, n)}(x, y) = \Theta_{[c,\,d]\times[a,\,b]}^{(n, m)}(y, x). \nn
\end{eqnarray}
\noindent (ii) \  There exists a positive constant $C$, such that
\begin{eqnarray}
	\big|\Theta_{[a,\,b]\times[c,\,d]}^{(m, n)}(x, y)\big| \leq C(b-a)^m(d-c)^n \int_c^d \int_a^b {|w(\xi,\eta)|}\,d\xi d\eta \nn
\end{eqnarray}
\noindent (iii) \ For $0 \le k \le m$ and $0 \le l \le n$, we have
\beas
\p_{k, l} \Theta^{(m,n)}_{[a, b] \times [c, d]}(x, y) = \Theta^{(m-k,\, n-l)}_{[a, b] \times [c, d]}(x, y).
\eeas
Here, we  denote $\p_{m, n} f(x, y) = \p_x^m\p_y^n f(x, y)$
%\beas
%\p_{m, n} f(x, y) = \fl{\p^{m+n}f(x, y)}{\p x^m \p y^n}, \qquad\mbox{for} \ \ m, n \in {\mathbb N}^0.
%\eeas
as a partial derivative of $f$.
\end{property}
The properties (i) and (ii) are implied from its definition, and the property (iii) can be obtained by using the Leibniz integral rule. Here, we will omit their proofs for brevity. 
Next, we introduce the following lemma from the weighted Montgomery identity of two variables. 
%%%% %%%% %%%% %%%% Lemma 2
\begin{lemma}[Extension of the weighted Montgomery identity]\label{lemma2-2D}
Let $w, f: [a,\,b]\times[c,\,d]\rightarrow \mathbb{R}$ be integrable functions, and $m, n \in {\mathbb N^0}$. 
Define 
\beas
Q = \int_{c}^{d}\int_{a}^{b} \Big(4f(x, y)- \big(f(a,c) + f(a,d) +f(b,c) + f(b,d)\big)
\Big)w(x, y)\,dx dy. 
\eeas
Assume that the derivative $\p_{m,n}f$  exists and is integrable for the following $m$ and $n$. 
Then, 
\vskip 5pt
\noindent (i) \ If  $m, n = 0, 1$ and $m + n \le1$,  we have
\begin{eqnarray}
Q &=& -\fl{1}{2}\int_{c}^{d}\int_{a}^{b} \left(\p_{0,1}\Theta_{[a,\,b]\times[c,\,d]}^{(0,0)}(x, y)\, \partial_{1,0}f(x, y) + \p_{1, 0}\Theta_{[a,\,b]\times[c,\,d]}^{(0,0)}(x, y)\, \partial_{0,1}f(x, y)\right)dx dy \nn\\
&&-\fl{1}{2}
\int_{a}^{b} \left( \Theta_{[a,\,b]\times[c,\,d]}^{(0,0)}(x,d)  
\, \partial_{1,0}f(x,d) - \Theta_{[a,\,b]\times[c,\,d]}^{(0,0)}(x,c)  
\, \partial_{1,0}f(x,c) \right)dx \nn\\
&&-\fl{1}{2}
 \int_{c}^{d}
\left(\Theta_{[a,\,b]\times[c,\,d]}^{(0,0)}(b, y)
\, \partial_{0,1}f(b, y)-\Theta_{[a,\,b]\times[c,\,d]}^{(0,0)}(a, y)
\, \partial_{0,1}f(a, y)\right) dy. \qquad\qquad\qquad\qquad  \nn
\end{eqnarray}
%%%%
\noindent(ii) \ If  $m, n = 0, 1, 2$ and $m + n \le 2$,  we have
\bea
Q &=& \int_{c}^{d}\int_{a}^{b} \Theta_{[a,\,b]\times[c,\,d]}^{(0,0)}(x,y)\,\partial_{1,1}f(x, y) dx dy\nn\\
&&+\int_{a}^{b}\left(\Theta_{[a,\,b]\times[c,\,d]}^{(1,0)}(x,d)\partial_{2, 0}f(x,d)
-\Theta_{[a,\,b]\times[c,\,d]}^{(1, 0)}(x,c)\partial_{2, 0}f(x,c)\right)dx\qquad\qquad\qquad \nn\\
&&+ \int_{c}^{d}  
\left(\Theta_{[a,\,b]\times[c,\,d]}^{(0,1)}(b, y)\partial_{0, 2}f(b, y)
- \Theta_{[a,\,b]\times[c,\,d]}^{(0,1)}(a, y)\partial_{0,2}f(a, y)\right) dy \nn\\
&&- 
\sum_{\substack{m, n = 0, 1 \\ m+n\neq 0,2}}\left(\Theta_{[a,\,b]\times[c,\,d]}^{(m, n)}(b,d)\,\partial_{m, n}f(b,d) - 
\Theta_{[a,\,b]\times[c,\,d]}^{(m, n)}(b,c)\,\partial_{m,n}f(b,c) \right. \nn\\
&&\hspace{1.9cm}
- \Theta_{[a,\,b]\times[c,\,d]}^{(m, n)}(a,d)\partial_{m, n}f(a,d)
+\Theta_{[a,\,b]\times[c,\,d]}^{(m, n)}(a,c)\partial_{m, n}f(a,c)\Big).\qquad\qquad\qquad\qquad  \nn
\eea
%%%%
\noindent (iii) \ If  $m, n = 0, 1, 2$ and $m + n \le 3$,  we have 
\bea
Q &=& -\fl{1}{2}\int_{c}^{d}\int_{a}^{b} \left(\Theta_{[a,\,b]\times[c,\,d]}^{(0,1)}(x,y)\,\partial_{1,2}f(x, y) + \Theta_{[a,\,b]\times[c,\,d]}^{(1,0)}(x,y)\,\partial_{2,1}f(x, y)\right)dx dy\nn\\
&&+\sum_{n=0, 1}  \left(-\fl{1}{2}\right)^{n}
\bigg(\int_{a}^{b}\left(\Theta_{[a,\,b]\times[c,\,d]}^{(1,n)}(x,d)\partial_{2, n}f(x,d)
-\Theta_{[a,\,b]\times[c,\,d]}^{(1, n)}(x,c)\partial_{2, n}f(x,c)\right)dx\qquad\qquad\qquad \nn\\
&&\hspace{2.5cm}+ \int_{c}^{d}  \left(\Theta_{[a,\,b]\times[c,\,d]}^{(n,1)}(b, y)\partial_{n, 2}f(b, y)
- \Theta_{[a,\,b]\times[c,\,d]}^{(n,1)}(a, y)\partial_{n,2}f(a, y)\right) dy\bigg) \nn\\
&&+ \sum_{\substack{m, n = 0, 1 \\ m+n\neq 0}} (-1)^{m+n}
\left(\Theta_{[a,\,b]\times[c,\,d]}^{(m, n)}(b,d)\,\partial_{m, n}f(b,d) - 
\Theta_{[a,\,b]\times[c,\,d]}^{(m, n)}(b,c)\,\partial_{m,n}f(b,c) \right. \nn\\
&&  \hspace{3.0cm} - \Theta_{[a,\,b]\times[c,\,d]}^{(m, n)}(a,d)\partial_{m, n}f(a,d)
+\Theta_{[a,\,b]\times[c,\,d]}^{(m, n)}(a,c)\partial_{m, n}f(a,c)\Big).\nn\eea
\end{lemma}
\begin{proof}
The proof of Lemma \ref{lemma2-2D} can be  done by  first averaging the weighted Montgomery identity \cite[Theorem 2.2]{Khan2013} at points $(a, c)$, $(a, d)$, $(b, c)$ and $(b, d)$, and then using integration by parts.  
\end{proof}
For the sake of completeness, we will review the Chebyshev integral inequality for two-variable functions as follows, which will be frequently used in the proof of our theorems. 
The Chebyshev integral inequality for  multiple variable functions can be found in \cite[Theorem A]{Beesack1985}. 
%%%#######
\begin{lemma}[Chebyshev integral inequality]\label{lemma4}
Let $f, g : [a, b]\times[c, d] \to {\mathbb R}$ be continuous, nonnegative, and similarly ordered, i.e., $\big(f(x_1, y_1)-f(x_2, y_2)\big)\big(g(x_1, y_1)-g(x_2, y_2)\big) \ge 0$, for any points $(x_1, y_1)$ and $(x_2, y_2)$. Then, it follows that 
\begin{eqnarray*}\label{Chebyshev}
\bigg(\int_c^d \int_a^b f(x, y)\,dx dy \bigg)\bigg(\int_c^d \int_a^b g(x, y)\,dx dy\bigg)
\leq (b-a)(d-c) \int_c^d \int_a^b f(x, y)g(x,y) \,dx dy. 
\end{eqnarray*}
\end{lemma}
%\begin{proof}
%Since the functions $f$ and $g$ have the same monotonicity, we have
%\begin{eqnarray}
%0 &\leq& \int_c^d \int_a^b \int_c^d \int_a^b \big(f(x_1 y_1)-f(x_2, y_2)\big)\big(g(x_1, y_1)-g(x_2, y_2)\big)\,dx_2\,dy_2\,dx_1\,dy_1 \nn\\
%&=& \int_c^d \int_a^b \int_c^d \int_a^b \big(f(x_1, y_1)g(x_1, y_1)-f(x_1, y_1)g(x_2, y_2)\nn\\
%&&\hspace{2.2cm}-f(x_2, y_2)g(x_1, y_1)+f(x_2, y_2)g(x_2, y_2)\big) \,dx_2\,dy_2 \,dx_1\,dy_1 \nn\\
%&=& 2\bigg((b-a)(d-c)\int_c^d \int_a^b f  g\,dx\,dy - \Big(\int_c^d \int_a^b f(x, y)\,dx\,dy\Big)\Big(\int_c^d \int_a^b g(x, y) \,ds\,ds\Big)\bigg),\nn
%\end{eqnarray}
%which implies (\ref{Chebyshev}) immediately.
%\end{proof}

To prepare our main theorems, we will first discuss the properties of function $\psi_\gm(\bx, \bxi)$. 
For notational simplicity, we will omit ${\bf x}$, and let $\psi_{\gamma}({\bxi}) := \psi_{\gamma}({\bf x}, \bxi)$ in the rest of this section.
 
%%%% %%%% %%%% %%%% Definition
\begin{definition}
For  $k \in {\mathbb N}^0$ and $s \in (0, 1]$, let $C^{k, s}({\mathbb R}^d)$ denote the space that consists of all functions $u: {\mathbb R}^d \to {\mathbb R}$ with continuous partial derivatives of order less than or equal to $k$, whose $k$-th partial derivatives are uniformly H\"older continuous with exponent $s$.
\end{definition}

%%%% %%%% %%%% %%%% Lemma \psi
\begin{lemma}\label{lemma1-2D}
Let $0 < s \le 1$, $\bxi \in {\mathbb R}_+^2\backslash\{\bf 0\}$,  and $m, n \in {\mathbb N}^0$.
\begin{enumerate}\itemsep -1pt
%%%% %%%%
\item[] (i)  If $u \in C^{0, s}(\mathbb{R}^2)$, there exists a constant $C > 0$, such that $\big|\psi_0(\bxi)\big| \leq C|\bxi|^{s}$. 
\item[] (ii) If $u \in C^{1, s}(\mathbb{R}^2)$, there exists a constant $C > 0$, such that for any $\gm \in (\ap, 2]$,
\beas
\big| \partial_{m,n}\psi_{\gamma}(\bxi)\big|
\leq C|\bxi|^{s+1-\gamma-(m+n)},\qquad \mbox{for \ $m + n \le 1$}.
\eeas
\item[] (iii) Let $u \in C^{2+\lfloor\ap\rfloor, \, s}({\mathbb R}^2)$ for $\ap\in(0, 2)$.
Then we have
\begin{eqnarray}\label{psi-0}
\big|\partial_{m,n}\psi_2(\xi,\eta) + \partial_{n,m}\psi_2(\eta,\xi)\big|
\leq C|\bxi|^{s + \lfloor\ap\rfloor- (m+n)}, \qquad \mbox{for} \ \ 0 < m + n \le 2+\lfloor\ap\rfloor,
\end{eqnarray}
with the constant $C > 0$. 
%Moreover, if one of $\xi$ and $\eta$ equals to zero, we have
%\bea\label{psi-1}
%\big|\p_{1,0}\psi_2(\xi, 0)\big| \le C\xi^{s -1}, \qquad 
%\big|\p_{0,1}\psi_2(0, \eta)\big| \le C\eta^{s -1}.
%\eea
\end{enumerate} 

\end{lemma}
The proof of the above properties can be done by directly applying Taylor's theorem.

%%%% %%%% %%%% %%%% Theorem 1 -- Minimum smoothness requirement
\begin{theorem}[{\bf Minimum consistent conditions}]\label{thm1}
Suppose that $u$ has finite support on the domain $\Og \subset {\mathbb R}^2$.
Let $(-\Dt)_{h,\gamma}^{\fl{\ap}{2}}$ be the finite difference approximation of the fractional Laplacian $(-\Dt)^{\fl{\ap}{2}}$, with $h$ a small mesh size.  
If $u \in C^{\lfloor\ap\rfloor,\,\ap-\lfloor\ap\rfloor+\veps}({\mathbb R}^2)$ with $0 < \veps \le  1 + \lfloor\ap\rfloor -\ap$, then for any splitting parameter $\gm \in (\ap, 2]$,  the local truncation error satisfies
\bea \label{thm1-error}
 \big\|(-\Dt)^{\fl{\ap}{2}} u  - (-\Dt)^{\fl{\ap}{2}}_{h,\gamma} u\big\|_\infty \le Ch^{\veps},  \qquad  \mbox{for} \ \ \ap\in(0, 2)
\eea
with $C$ a positive constant independent of $h$. 
\end{theorem}
%%%% %%%% %%%% %%%% Proof
\begin{proof}
Introduce the error function at point $\bx \in\Og$ as: 
\bea\label{error}
e_{\alpha,\gamma}^{h}({\bf x}) &=& (-\Dt)^{\fl{\ap}{2}} u({\bf x}) - (-\Dt)^{\fl{\ap}{2}}_{h,\gamma} u({\bf x})\nn\\
&=&-\fl{c_{2,\alpha}}{4}
\bigg[\int_{I_{00}}\Big(4\psi_{\gamma}(\bxi) 
- \sum_{(m, n)\in\widetilde{\varkappa}_1}c_{mn}^\gm\psi_{\gamma}\big(\bxi_{mn}\big)\Big)
w_{\gm}(\bxi)d\xi d\eta \nn\\
&&+ \sum_{(i,j) \in \widetilde{\varkappa}_{N-1}} \int_{I_{ij}} \Big( 4\psi_{\gamma}(\bxi) - \sum_{(m, n)\in \varkappa_1}\psi_\gm\big(\bxi_{(i+m)(j+n)}\big)\Big) w_{\gm}(\bxi) d\xi d\eta \bigg]\qquad  \nn\\
&=& -\fl{c_{2,\alpha}}{4} (I+II), \qquad \bx \in \Og,
\eea
which is obtained from {(\ref{fL1-2D})--(\ref{I-2D}) and (\ref{I00-final})}. 
For simplicity,  we denote the index set
\beas
\varkappa_{N} = \big\{(i, j)\,|\, i, j = 0, 1, \ldots, N\}, \qquad \widetilde{\varkappa}_{N} = \varkappa_{N}\backslash(0, 0).
\eeas
We will then prove the cases of $\ap \in (0, 1)$ and $\ap \in [1, 2)$ separately. 
\bigskip 

\noindent{\bf Case (i) (For $u\in C^{0, \ap + \veps}({\mathbb R}^2)$ with $\alpha \in (0,1)$): }  
Here, we rewrite  $\psi_\gm(\bxi) = |\bxi|^{-\gm}\psi_0(\bxi)$. 
%%%% %%%% Term I
Using the triangle inequality and Lemma \ref{lemma1-2D} (i) with $s = \ap+\veps$ to the term $I$, we  obtain
\bea\label{thm11-termI}
|\,I \, | &=& \bigg|\int_0^h \int_0^h
\Big(4\psi_0(\bxi)|\bxi|^{-(2+\alpha)} - \sum_{(m, n) \in \widetilde{\varkappa}_1}c_{mn}^\gm\,\psi_0(\bxi_{mn})\,|\bxi_{mn}|^{-\gm} |\bxi|^{\gamma-(2+\alpha)}\,d\xi d\eta\bigg|\nn\\
&\le& C\bigg(\int_0^h \int_0^h\big|\psi_0(\bxi)\big| |\bxi|^{-(2+\alpha)}\,d\xi d\eta + \sum_{(m, n) \in \widetilde{\varkappa}_1}\int_0^h \int_0^h\big|\psi_0(\bxi_{mn})\big|\,|\bxi_{mn}|^{-\gm} |\bxi|^{\gamma-(2+\alpha)}d\xi d\eta\bigg)\quad \nn\\
&\le& C\bigg(\int_0^h \int_0^h |\bxi|^{-(2-\veps)}d\xi d\eta  + h^{\ap+\veps-\gm}\int_0^h \int_0^h |\bxi|^{\gm-(2+\ap)}d\xi d\eta \bigg)  \le \, Ch^{\varepsilon},
\eea
where the last inequality is obtained by using the following property:  for $p < 2$,  we have
\bea\label{xi0}
\int_0^h\int_0^h |\bxi|^{-p} d\xi d\eta \,= \, 2\int_0^{\frac{\pi}{4}}\int_0^{h\sec \theta} r^{-p+1} \,dr d\theta = \fl{2}{2-p}h^{2-p}\int_0^{\fl{\pi}{4}}(\sec\theta)^{2-p}d\theta\, \leq \, C{h}^{2-p}\quad
\eea 
with $C$ a positive constant independent of $h$. 

%%%% %%%% Term II
For term $II$, we first rewrite it by adding and subtracting terms and then use the triangle inequality, {Taylor's} theorem, and Lemma \ref{lemma1-2D} (i) with $s = a+\veps$, to obtain
\bea \label{thm11-termII}
|\,II\,| &=& \bigg|\sum_{(i,j) \in \widetilde{\varkappa}_{N-1}} \int_{I_{ij}} |\bxi|^{-\gm}\Big( 4\psi_0(\bxi) - \sum_{(m, n)\in\varkappa_1}\psi_0\big(\bxi_{(i+m)(j+n)}\big)\Big) |\bxi|^{\gamma-(2+\alpha)} d\bxi \nn\\
&&+\sum_{(i,j) \in \widetilde{\varkappa}_{N-1}} \int_{I_{ij}} \Big(\sum_{(m, n)\in\varkappa_1}\big(|\bxi|^{-\gm}-|\bxi_{(i+m)(j+n)}|^{-\gm}\big) \psi_0\big(\bxi_{(i+m)(j+n)}\big)\Big)|\bxi|^{\gamma-(2+\alpha)}\, d\bxi \bigg|\nn\\
&\leq& \bigg|\sum_{(i,j) \in \widetilde{\varkappa}_{N-1}} \int_{I_{ij}} |\bxi|^{-\gm}\Big( 4\psi_0(\bxi) - \sum_{(m, n)\in\varkappa_1}\psi_0\big(\bxi_{(i+m)(j+n)}\big)\Big) |\bxi|^{\gamma-(2+\alpha)} d\bxi\bigg| \nn\\
&&+\bigg|\sum_{(i,j) \in \widetilde{\varkappa}_{N-1}} \int_{I_{ij}} \Big(\sum_{(m, n)\in\varkappa_1}\big(|\bxi|^{-\gm}-|\bxi_{(i+m)(j+n)}|^{-\gm}\big) \psi_0\big(\bxi_{(i+m)(j+n)}\big)\Big)|\bxi|^{\gamma-(2+\alpha)} d\bxi \bigg|\nn\\
&\leq& C\bigg(\sum_{(i,j) \in \widetilde{\varkappa}_{N-1}} \int_{I_{ij}} |\bxi|^{-\gm} h^{\alpha+\varepsilon} |\bxi|^{\gamma-(2+\alpha)}\,d\bxi 
+ \sum_{(i,j) \in \widetilde{\varkappa}_{N-1}} \int_{I_{ij}} h|\bxi|^{-\gm-1}|\bxi|^{\alpha+\varepsilon}|\bxi|^{\gamma-(2+\alpha)}\, d\bxi \bigg)\nn\\
&\leq& Ch^{\varepsilon}, \qquad 
\eea
where the last inequality is obtained from the following property:  for $p > 0$, we have
\bea\label{prop2}
&&\sum_{(i, j) \in \widetilde{\varkappa}_{N-1}} \int_{\eta_{j}}^{\eta_{j+1}}\int_{\xi_{i}}^{\xi_{i+1}}   
|\bxi|^{-p}\,d\xi d\eta=\int_h^L\int_h^L|\bxi|^{-p}\,d\xi d\eta + 2\int_h^L\int_0^h |\bxi|^{-p}\,d\xi d\eta \qquad\qquad\quad \nn\\
&& \hspace{5.1cm} \le \bigg(\int_h^L\xi^{-\fl{p}{2}}\,d\xi\bigg)^2 + 2\int_h^L\int_0^h |\bxi|^{-p}\,d\xi d\eta \, \le \,
Ch^{\min\{0, 2-p\}}.  \qquad
\eea
Combining (\ref{thm11-termI}) with (\ref{thm11-termII}) leads to (\ref{thm1-error}) immediately. \\

%%%% %%%% %%%% %%%% Case 2
\noindent {\bf Case {(ii)} (For $u\in C^{1, \ap -1 + \veps}({\mathbb R}^2)$ with $\alpha \in [1, 2)$): }  
Using the extension of the weighted Montgomery identity in Lemma \ref{lemma2-2D} (i) to the error function (\ref{error}) with $w(\bxi) = |\bxi|^{\gm-(2+\ap)}$, we  get 
\bea\label{error1}
e_{\alpha,\gamma}^{h}({\bf x})&=&-\fl{c_{2,\alpha}}{4}
\bigg[\int_{I_{00}}
\Big(4\psi_{\gamma}(\bxi) 
-\sum_{(m, n)\in\widetilde{\varkappa}_1}c_{mn}^\gm\,\psi_\gm(\bxi_{mn})\Big)
|\bxi|^{\gamma-(2+\alpha)} d\xi d\eta \nn \qquad \\
&&\hspace{0.5cm}-\fl{1}{2}\sum _{(i,\, j) \in \widetilde{\varkappa}_{N-1}}
\int_{I_{ij}} \Big(\p_{0,1}\Theta_{I_{ij}}^{(0,0)} (\xi, \eta)\, \partial_{1,0}\psi_{\gamma}(\xi, \eta) + 
\p_{1,0}\Theta_{I_{ij}}^{(0, 0)} (\xi, \eta)\, \partial_{0,1}\psi_{\gamma}(\xi, \eta)\Big)\,d\xi d\eta \nn\\
&&\hspace{0.5cm}-\fl{1}{2}\sum _{(i,\, j) \in \widetilde{\varkappa}_{N-1}}
\int_{\xi_{i}}^{\xi_{i+1}}  \Big(\Theta_{I_{ij}}^{(0,0)}(\xi,\eta_{j+1})\,\partial_{1,0}\psi_{\gamma}(\xi,\eta_{j+1})-\Theta_{I_{ij}}^{(0,0)}(\xi,\eta_{j})\,\partial_{1,0}\psi_{\gamma}(\xi,\eta_{j})\Big)d\xi\nn\\ 
&&\hspace{0.5cm}-\fl{1}{2}\sum _{(i,\, j) \in \widetilde{\varkappa}_{N-1}}
\int_{\eta_{j}}^{\eta_{j+1}}\Big(\Theta_{I_{ij}}^{(0,0)}(\xi_{i+1},\eta)\,\partial_{0,1}\psi_{\gamma}(\xi_{i+1},\eta)-\Theta_{I_{ij}}^{(0,0)}(\xi_{i},\eta)\,\partial_{0,1}\psi_{\gamma}(\xi_{i},\eta)\Big)d\eta\nn\\
&=& -\fl{c_{2,\ap}}{4} \big(I + II + III + IV\big).
\eea

%%%% %%%% Term I
For term $I$, we use the triangle inequality and then Lemma \ref{lemma1-2D} (ii) with $m = n = 0$ to obtain
\bea\label{termI}
|\,I \, | &=& \bigg|\int_0^h \int_0^h
\Big(4\psi_\gm(\bxi)|\bxi|^{\gm-(2+\alpha)} - \sum_{(m, n)\in\widetilde{\varkappa}_1}c_{mn}^\gm\,\psi_\gm(\bxi_{mn})\Big)\, |\bxi|^{\gamma-(2+\alpha)}\,d\xi d\eta\bigg|\nn\\
&\le& C\bigg(\int_0^h \int_0^h\big|\psi_\gm(\bxi)\big| |\bxi|^{\gm-(2+\alpha)}\,d\xi d\eta + \sum_{(m, n)\in\widetilde{\varkappa}_1}\int_0^h \int_0^h\big|\psi_\gm(\bxi_{mn})\big| |\bxi|^{\gamma-(2+\alpha)}d\xi d\eta\bigg)\quad \nn\\
&\le& C\bigg(\int_0^h \int_0^h |\bxi|^{-(2-\veps)}d\xi d\eta  + h^{\ap+\veps-\gm}\int_0^h \int_0^h |\bxi|^{\gm-(2+\ap)}d\xi d\eta \bigg)  \le \, Ch^{\varepsilon},
\eea
where the last inequality is obtained by using the property (\ref{xi0}). 

%%%% %%%% Term II
For term $II$, by the triangle inequality, Property \ref{K-Property}, Lemma \ref{lemma1-2D} (ii) with $s = \ap - 1 +\veps$, and the Chebyshev integral inequality, we obtain 
\bea\label{termII}
|\,II\,| &=& \fl{1}{2}\bigg|\sum _{(i,\, j) \in \widetilde{\varkappa}_{N-1}}
\int_{I_{ij}} \Big(\p_{0,1}\Theta_{I_{ij}}^{(0,0)} (\xi, \eta)\, \partial_{1,0}\psi_{\gamma}(\xi, \eta) + 
\p_{1,0}\Theta_{I_{ij}}^{(0, 0)} (\xi, \eta)\, \partial_{0,1}\psi_{\gamma}(\xi, \eta)\Big)\,d\xi d\eta\bigg|\nn\\
&\le&\fl{1}{2}\sum _{(i,\, j) \in \widetilde{\varkappa}_{N-1}}
\int_{I_{ij}} \Big(\big|\p_{0,1}\Theta_{I_{ij}}^{(0, 0)} (\xi, \eta)\big| \big|\partial_{1,0}\psi_{\gamma}(\xi, \eta)\big| + 
\big|\p_{1,0}\Theta_{I_{ij}}^{(0, 0)} (\xi, \eta)\big| \big|\partial_{0,1}\psi_{\gamma}(\xi, \eta)\big|\Big) d\xi d\eta\nn\\
&\leq& \fl{C}{h}\sum_{(i, j) \in \widetilde{\varkappa}_{N-1}}
\bigg(\int_{\eta_{j}}^{\eta_{j+1}}\int_{\xi_{i}}^{\xi_{i+1}} |\bxi| ^{\gamma-(2+\alpha)}\,d\xi d\eta
\bigg)
\bigg(\int_{\eta_{j}}^{\eta_{j+1}}\int_{\xi_{i}}^{\xi_{i+1}} |\bxi|^{\ap+\veps-(1+\gamma)}\,d\xi d\eta\bigg)\nn\\
&\leq& Ch\sum_{(i, j) \in \widetilde{\varkappa}_{N-1}}
\int_{\eta_{j}}^{\eta_{j+1}}\int_{\xi_{i}}^{\xi_{i+1}} |\bxi| ^{-(3-\veps)} \,d\xi d\eta \ \le \ Ch^{\veps}
\eea
where the last inequality is obtained by the property (\ref{prop2}) with $p = 3-\veps$.

%%%% %%%% Term III
Following similar lines as in obtaining (\ref{termII}), i.e., using the triangle inequality,  Property \ref{K-Property} (ii), Lemma \ref{lemma1-2D} (ii), and the Chebyshev integral inequality, we obtain the estimate of term $III$ as:
\bea\label{termIII}
|\,III\,| &=& \fl{1}{2}\bigg|\sum_{(i, j) \in \widetilde{\varkappa}_{N-1}}
\int_{\xi_{i}}^{\xi_{i+1}}  
\Big(\Theta_{I_{ij}}^{(0,0)}(\xi,\eta_{j+1})\,\partial_{1,0}\psi_{\gamma}(\xi,\eta_{j+1})-\Theta_{I_{ij}}^{(0,0)}(\xi,\eta_{j})\,\partial_{1,0}\psi_{\gamma}(\xi,\eta_{j})
\Big) d\xi \,\bigg|\nn\\
&\le&  C\sum_{(i, j) \in \widetilde{\varkappa}_{N-1}}\int_{\xi_{i}}^{\xi_{i+1}} 
\bigg(\int_{\eta_{j}}^{\eta_{j+1}}\int_{\xi_{i}}^{\xi_{i+1}}   
|\bxi|^{\gamma-(2+\alpha)} d\xi d\eta\bigg)
\bigg(\big|\partial_{1,0}\psi_{\gamma}(\xi,\eta_{j})\big| + \big|\partial_{1,0}\psi_{\gamma}(\xi,\eta_{j+1})\big|\bigg)d\xi\nn\\
&\le& C\sum_{(i, j) \in \widetilde{\varkappa}_{N-1}}\bigg(\int_{\eta_{j}}^{\eta_{j+1}}\int_{\xi_{i}}^{\xi_{i+1}}   
|\bxi|^{\gamma-(2+\alpha)}\,d\xi d\eta\bigg)
\bigg(h^{-1}\int_{\eta_{j}}^{\eta_j+1}\int_{\xi_{i}}^{\xi_{i+1}} |\bxi|^{\ap+\veps-(1+\gm)} d\xi d\eta\bigg)\nn\\
&\le& C h\sum_{(i, j) \in \widetilde{\varkappa}_{N-1}} \int_{\eta_{j}}^{\eta_{j+1}}\int_{\xi_{i}}^{\xi_{i+1}}   
|\bxi|^{-(3-\veps)}\,d\xi d\eta \ \le \ Ch^{\veps},
\eea
by the property (\ref{prop2}). Following the same lines, we can obtain the estimate of term $IV$ as:
\bea\label{termIV} 
|\,IV\,| &\le& Ch^{\veps}.
\eea
%%%% %%%% %%%% %%%% 
Combining (\ref{error1}) with (\ref{termI}), (\ref{termII})--(\ref{termIV}) yields the error estimate in (\ref{thm1-error}).\\
\end{proof}
Theorem \ref{thm1} shows that our finite difference method is consistent if the function $u \in C^{\lfloor\ap\rfloor,\,\ap-\lfloor\ap\rfloor+\veps}({\mathbb R}^2)$, for $0 < \veps \le 1+\lfloor\ap\rfloor-\ap$. 
Furthermore, it implies that the consistent condition that is required when $\ap < 1$ is much weaker than that for $\ap \ge 1$.

%%%% %%%% %%%% %%%% Theorem 2 -- Maximum smoothness requirement
\begin{theorem}[{\bf Second order of accuracy}]\label{thm2}
Suppose that $u$ has finite support on the domain $\Og \subset {\mathbb R}^2$.
Let $(-\Dt)_{h,\gamma}^{\fl{\ap}{2}}$ be the finite difference approximation of the fractional Laplacian $(-\Dt)^{\fl{\ap}{2}}$, with $h$ a small mesh size.  
If the splitting parameter $\gamma = 2$ and $u \in C^{2+\lfloor\ap\rfloor,\,\ap-\lfloor\ap\rfloor+\veps}({\mathbb R}^2)$ with $0 < \veps \le 1+\lfloor\ap\rfloor-\ap$,  then the local truncation error satisfies
\bea\label{thm2-error}
 \big\|(-\Dt)^{\fl{\ap}{2}} u  - (-\Dt)^{\fl{\ap}{2}}_{h,\gamma} u\big\|_\infty \le Ch^2,  \qquad  \mbox{for} \ \ \ap\in(0, 2)
\eea
with $C$ a positive constant independent of $h$. 
\end{theorem}

\begin{proof} 
In the following, we will use $\gm = 2$ and prove the cases of $\ap < 1$ and $\ap \ge 1$ separately. \\

\noindent{\bf Case (i) (For $u\in C^{2, \ap+\veps}({\mathbb R}^2)$ with $\alpha \in (0,1)$): } 
Taking $\gm = 2$ and using  Lemma \ref{lemma2-2D} (ii) with $\og(\bxi) = |\bxi|^{-\ap}$ to the error function (\ref{error}), we get: 
\begin{eqnarray} \label{Error2}
&&e_{\alpha, 2}^{h}({\bf x}) = -\fl{c_{2,\ap}}{4}
\bigg[\int_{I_{00}} 
\Big(4\psi_{2}(\bxi)-2 \big(\psi_{2}(\xi_1,0) + \psi_{2}(0,\eta_1)\big)\Big)|\bxi|^{-\alpha}\,d\xi d\eta \nn \\
&&\hspace{1cm}+\sum_{(i, j) \in \widetilde{\varkappa}_{N-1}} \int_{I_{ij}} \Theta_{I_{ij}}^{(0,0)} (\xi,\eta)\partial_{1,1}{\psi_2}(\xi,\eta)\,d\xi d\eta\nn\\
&&\hspace{1cm}+\sum_{(i, j) \in \widetilde{\varkappa}_{N-1}}\bigg(
\int_{\xi_{i}}^{\xi_{i+1}}  \Big(\Theta_{I_{ij}}^{(1,0)}(\xi,\eta_{j+1})\partial_{2,0}{\psi_2}(\xi,\eta_{j+1})
-\Theta_{I_{ij}}^{(1,0)}(\xi,\eta_{j})\partial_{2, 0}{\psi_2}(\xi,\eta_{j})\Big)d\xi \nn\\
&&\hspace{2.9cm}
+\int_{\eta_{j}}^{\eta_{j+1}} \Big(
\Theta_{I_{ij}}^{(0,1)}(\xi_{i+1},\eta)\partial_{0, 2}{\psi_2}(\xi_{i+1},\eta)- 
\Theta_{I_{ij}}^{(0,1)}(\xi_{i},\eta)\partial_{0, 2}{\psi_2}(\xi_{i},\eta)\Big)d\eta\bigg)\nn\\
&&\hspace{1cm}-\sum_{(i, j) \in \widetilde{\varkappa}_{N-1}} \sum_{\substack{m+n=1}}\Big(\Theta_{I_{ij}}^{(m,n)}(\bxi_{(i+1)(j+1)})\partial_{m,n}\psi_2(\bxi_{(i+1)(j+1)})+\Theta_{I_{ij}}^{(m,n)}(\bxi_{ij})\partial_{m,n}\psi_2(\bxi_{ij})\nn\\
&&\hspace{3.68cm}-\Theta_{I_{ij}}^{(m,n)}(\bxi_{i(j+1)})\partial_{m,n}\psi_2(\bxi_{i(j+1)})-\Theta_{I_{ij}}^{(m,n)}(\bxi_{(i+1)j})\partial_{m,n}\psi_2(\bxi_{(i+1)j})\Big) \qquad\quad \nn\\ 
&&\hspace{1.5cm}=-\fl{c_{2,\alpha}}{4}\big(I + II +III + IV\big).
\end{eqnarray}
%%%% %%%% Term 1
Denote $\Psi_2(\xi,\eta) := \psi_2(\xi,\eta) + \psi_2(\eta,\xi)$, and rewrite the term $I$ as
\begin{eqnarray}\label{term1}
|\,I\,| &=& 2\bigg|\int_0^h\int_0^h \Big(\psi_2(\xi,\eta)+\psi_2(\eta,\xi) - \big(\psi_2(h,0)+\psi_2(0,h)\big)\Big) |\bxi|^{-\alpha}d\xi d\eta\,\bigg| \nn\\
&=& 2\bigg|\int_0^h\int_0^h \Big(\Psi_2(\xi,\eta)-\Psi_2(\xi,0)+\Psi_2(\xi,0)-\Psi_2(h,0)\Big)|\bxi|^{-\alpha}d\xi d\eta\bigg| \nn\\
&=& 2\bigg|\int_0^h\int_0^h \bigg(\int_0^{\eta}\partial_{0,1}\Psi_2(\xi,\widetilde{\eta})\,d\widetilde{\eta}+\int^h_{\xi}\partial_{1,0}\Psi_2(\widetilde{\xi},0)\,d\widetilde{\xi}\bigg) |\bxi|^{-\alpha}d\xi d\eta\bigg| \nn
\eea
by Taylor's theorem. 
Then noticing the definition of $\Psi_2$, using Lemma \ref{lemma1-2D} (iii) and Chebyshev integral inequality {for one variable}, we obtain
\begin{eqnarray}\label{term1}
|\,I\,| %& \le & C\int_0^h\int_0^h \bigg(\int_0^{\eta}|\partial_{1,0}\psi_2(\widetilde{\eta},\xi)+\partial_{0,1}\psi_2(\xi,\widetilde{\eta})|\,d\widetilde{\eta}+\int_{\xi}^h|\partial_{1,0}\psi_2(\widetilde{\xi},0)+\partial_{0,1}\psi_2(0,\widetilde{\xi})|\,d\widetilde{\xi}\bigg)|\bxi|^{-\alpha}d\xi d\eta \nn\\
&\leq & C\int_0^h\int_0^h \bigg(\int_0^{h}\big(\xi^2 + \widetilde{\eta}^2\big)^{\fl{\alpha+\varepsilon-1}{2}}\,d\widetilde{\eta}+\int_0^{h}\widetilde{\xi}^{\alpha+\varepsilon-1}\,d\widetilde{\xi}\bigg)|\bxi|^{-\alpha}d\xi d\eta \nn\\
&\leq & Ch\int_0^h\int_0^h |\bxi|^{\veps-1} d\xi d\eta+Ch^{\alpha+\varepsilon}\int_0^h\int_0^h |\bxi|^{-\alpha}\,d\xi d\eta 
\leq  Ch^{2+\varepsilon},
\eea
where the last inequality is obtained by the property (\ref{xi0}).

%%%% %%%% Term 2
For term $II$, we first rewrite it by using Property \ref{K-Property} (i) as
\bea
|\,II \, | %&=& \fl{1}{2}\bigg|\sum_{(i, j) \in \widetilde{\varkappa}_{N-1}} \int_{I_{ij}} \Theta_{I_{ij}}^{(0,0)} (\xi,\eta)\partial_{1,1}{\psi_2}(\xi,\eta) + \Theta_{I_{ij}}^{(0,0)} (\xi,\eta)\partial_{1,1}{\psi_2}(\xi,\eta)\Big)\,d\xi d\eta\bigg|\nn\\
&=& \frac{1}{2}\bigg|\sum_{(i, j) \in \widetilde{\varkappa}_{N-1}} \int_{I_{ji}} \Theta_{I_{ij}}^{(0,0)} (\eta,\xi)\partial_{1,1}{\psi_2}(\eta,\xi)\,d\xi d\eta
  +\sum_{(i, j) \in \widetilde{\varkappa}_{N-1}}\int_{I_{ij}} \Theta_{I_{ij}}^{(0,0)} (\xi,\eta)\partial_{1,1}{\psi_2}(\xi,\eta)\,d\xi d\eta \bigg|\nn\\
&\le& \frac{1}{2}\sum_{(i, j) \in \widetilde{\varkappa}_{N-1}} \int_{I_{ij}} \big|\Theta_{I_{ij}}^{(0,0)} (\xi,\eta)\big|\big|\partial_{1,1}{\psi_2}(\eta,\xi) +\partial_{1,1}{\psi_2}(\xi,\eta)\big|\,d\xi d\eta, \nn
\eea
Then, using the triangle inequality, Property \ref{K-Property} (i),  Lemma \ref{lemma1-2D} {(iii)}, and the Chebyshev integral inequality, we obtain 
\begin{eqnarray}\label{term2}
|\,II\,| &\leq & C\sum_{(i,j)\in \widetilde{\varkappa}_{N-1}}\bigg(
\int_{\eta_{j}}^{\eta_{j+1}}\int_{\xi_{i}}^{\xi_{i+1}} |\bxi| ^{-\alpha} d\xi d\eta \bigg)
\bigg(\int_{\eta_{j}}^{\eta_{j+1}}\int_{\xi_{i}}^{\xi_{i+1}}  |\bxi|^{\ap-2+\veps} \,d\xi d\eta\bigg) \nn\\
&\leq & Ch^2\sum_{(i,j)\in \widetilde{\varkappa}_{N-1}}\int_{\eta_{j}}^{\eta_{j+1}}\int_{\xi_{i}}^{\xi_{i+1}}  |\bxi|^{-(2-\veps)} d\xi d\eta \,
\le \, C h^{2}, 
\end{eqnarray}
where the last inequality is obtained by (\ref{prop2}) with $p = 2-\veps$.

%%%% %%%% Term III
For term $III$, noticing that $\xi_i = \eta_i$ and applying Property \ref{K-Property} (i) and following the same lines as in obtaining (\ref{term2}), we get 
\bea\label{term3}
|\, III \,| 
&=&\bigg|\sum_{(i, j) \in \widetilde{\varkappa}_{N-1}}\bigg(
\int_{\xi_{i}}^{\xi_{i+1}}  \Big(\Theta_{I_{ij}}^{(1,0)}(\xi,\eta_{j+1})\partial_{2,0}{\psi_2}(\xi,\eta_{j+1})
-\Theta_{I_{ij}}^{(1,0)}(\xi,\eta_{j})\partial_{2, 0}{\psi_2}(\xi,\eta_{j})\Big)d\xi \nn\\
&&\hspace{1.60cm}
+\int_{\eta_{j}}^{\eta_{j+1}} \Big(
\Theta_{I_{ij}}^{(0,1)}(\xi_{i+1},\eta)\partial_{0, 2}{\psi_2}(\xi_{i+1},\eta)- 
\Theta_{I_{ij}}^{(0,1)}(\xi_{i},\eta)\partial_{0, 2}{\psi_2}(\xi_{i},\eta)\Big)d\eta\bigg)\bigg|\nn\\
&=& \bigg|\sum_{(i, j) \in \widetilde{\varkappa}_{N-1}}
\bigg(\int_{\xi_{i}}^{\xi_{i+1}} \Theta_{I_{ij}}^{(1,0)}(\xi,\eta_{j+1})
\big(\partial_{2,0}{\psi_2}(\xi,\eta_{j+1})+\partial_{0,2}{\psi_2}(\eta_{j+1},\xi)\big)\,d\xi\nn\\
&&\hspace{1.7cm}-\int_{\xi_{i}}^{\xi_{i+1}}  \Theta_{I_{ij}}^{(1,0)}(\xi,\eta_{j})
\big(\partial_{2,0}{\psi_2}(\xi,\eta_{j})+\partial_{0,2}{\psi_2}(\eta_{j},\xi)\big)\,d\xi\bigg)\bigg|\nn\\
&\le& C\sum_{(i, j) \in \widetilde{\varkappa}_{N-1}}
\bigg(h\int_{\eta_{j}}^{\eta_{j+1}}\int_{\xi_{i}}^{\xi_{i+1}} |\bxi| ^{-\alpha} d\xi d\eta \bigg)
\bigg(h^{-1}\int_{\eta_{j}}^{\eta_{j+1}}\int_{\xi_{i}}^{\xi_{i+1}}  |\bxi|^{\ap+\veps-2} \,d\xi d\eta\bigg) \nn\\
&\le& C\sum_{(i, j) \in \widetilde{\varkappa}_{N-1}}
\bigg(h^{2}\int_{\eta_{j}}^{\eta_{j+1}}\int_{\xi_{i}}^{\xi_{i+1}} |\bxi| ^{\veps-2} d\xi d\eta \bigg) \, \le \, Ch^{2},
\eea
by the property (\ref{prop2}).  
%%%% %%%% %%%% %%%% Term IV
Rewrite term $IV$ as
\beas
IV &=& \sum_{(i, j) \in \widetilde{\varkappa}_{N-1}} \sum_{\substack{m+n=1}}\Big(\Theta_{I_{ij}}^{(m,n)}(\bxi_{i(j+1)})\partial_{m,n}\psi_2(\bxi_{i(j+1)})+\Theta_{I_{ij}}^{(m,n)}(\bxi_{(i+1)j})\partial_{m,n}\psi_2(\bxi_{(i+1)j})\nn\\
&&\hspace{2.6cm}-\Theta_{I_{ij}}^{(m,n)}(\bxi_{(i+1)(j+1)})\partial_{m,n}\psi_2(\bxi_{(i+1)(j+1)})-\Theta_{I_{ij}}^{(m,n)}(\bxi_{ij})\partial_{m,n}\psi_2(\bxi_{ij})\Big) \qquad\quad\nn\\
&=& IV_1 + IV_2 + IV_3 + IV_4 + IV_5 + IV_6 + IV_7 + IV_8,
\eeas
where the terms
\beas
IV_1 &=&\sum_{m+n =1}\sum_{\substack{i, j = 1\\ i+j \neq 2}}^{N-1} \Big(\Theta^{(m, n)}_{I_{i(j-1)}}(\bxi_{ij}) -\Theta^{(m, n)}_{I_{ij}}(\bxi_{ij})\Big)\big(\p_{m,n}\psi_2(\xi_i, \eta_j) + \p_{n,m}\psi_2(\eta_j, \xi_i)\big),\nn\\
IV_2 &=&\sum_{m+n =1}\sum_{\substack{i, j = 1\\ i+j \neq 2}}^{N-1} \Big(\Theta^{(m, n)}_{I_{(i-1)j}}(\bxi_{ij}) -\Theta^{(m, n)}_{I_{(i-1)(j-1)}}(\bxi_{ij})\Big)\big(\p_{m,n}\psi_2(\xi_i, \eta_j) + \p_{n,m}\psi_2(\eta_j, \xi_i)\big),\nn\\
IV_3 &=&\sum_{m+n =1}\sum_{i=2}^{N-1} 
\Big(\Theta_{I_{(i-1)0}}^{(m, n)}(\bxi_{i0}) - \Theta_{I_{i0}}^{(m, n)}(\bxi_{i0})\Big)
\big(\p_{m,n}\psi_2(\xi_N, \eta_i) + \p_{n,m}\psi_2(\eta_i, \xi_N)\big),\nn\\
IV_4 &=&\sum_{m+n =1}\sum_{i=1}^{N-1} 
\Big(\Theta_{I_{(N-1)i}}^{(m, n)}(\bxi_{Ni}) - \Theta_{I_{(N-1)(i-1)}}^{(m, n)}(\bxi_{Ni})\Big)
\big(\p_{m,n}\psi_2(\xi_i, \eta_0) + \p_{n,m}\psi_2(\eta_0, \xi_i)\big),\nn\\
IV_5 &=&\sum_{m+n =1}\Theta_{I_{(N-1)0}}^{(m,n)}(\bxi_{N0}) \big(\p_{m,n}\psi_2(\xi_N, \eta_0)+\p_{n,m}\psi_2(\eta_0, \xi_N)\big),\nn\\
IV_6 &=&\sum_{m+n =1}-\Theta_{I_{10}}^{(m,n)}(\bxi_{10}) \big(\p_{m,n}\psi_2(\xi_1, \eta_0)+\p_{n,m}\psi_2(\eta_0, \xi_1)\big),\nn \\
IV_7 &=&\Theta_{I_{11}}^{(1,0)}(\bxi_{11}) \big(\p_{1, 0} \psi_2(\xi_1, \eta_1) + \p_{0, 1} \psi_2(\eta_1, \xi_1) \big),\nn\\
IV_8 &=&\Theta^{(1, 0)}_{I_{(N-1)(N-1)}}(\bxi_{NN})\big(\p_{1, 0} \psi_2(\eta_N, \xi_N) + 
\p_{0,1} \psi_2(\xi_N, \eta_N)\big).
\eeas
Here, the terms $IV_1$, $IV_2$, {$IV_3$, and $IV_4$}  can be estimated in the same manner. 
To show it, we will introduce the following property of $\Theta$.
Define an auxiliary function
\begin{eqnarray*}
 G^{(1, n )}(x) =\int_{\eta_j}^{\eta_{j+1}}\int_{x}^{{\xi_{i}}} |\bxi|^{-\alpha}(\xi_i - {\xi})(\eta_j-{\eta})^n\,d\xi d\eta, \qquad \mbox{for} \ \ n = 0,  1. 
\end{eqnarray*}
We then can write 
\begin{eqnarray}\label{eq100}
\Big|\Theta_{I_{(i-1)j}}^{(1, n)}(\xi_i, \eta_j) -\Theta_{I_{ij}}^{(1, n)}(\xi_{i},\eta_{j})\Big| &=& { \big|G^{(1,n)}(\xi_{i+1}) - G^{(1,n)}(\xi_{i-1})\big|/n!}\nn\\
&\leq& Ch^{4+n}\Big(\sqrt{\xi_i^2 + \eta_{j}^2}\Big)^{-(1+\alpha)}, \qquad \mbox{for} \ \ n = 0,  1.\qquad 
\end{eqnarray}
Using the triangle inequality, Lemma \ref{lemma1-2D} {(iii)} and (\ref{eq100}), we obtain      
\bea
|\, IV_1 \,| &\le&\sum_{m+n=1}\sum_{\substack{i, j=1 \\ i+j\neq 2}}^{N-1} \Big|\Theta_{I_{(i-1)j}}^{(m, n)}(\xi_i, \eta_j) -\Theta_{I_{ij}}^{(m, n)}(\xi_{i},\eta_{j})\Big|\big|\p_{m, n} \psi_2(\xi_i, \eta_j) + \p_{n, m} \psi_2(\eta_j, \xi_i) \big|\nn\\
&\le& C h^{4} \sum_{\substack{i, j=1 \\ i+j\neq 2}}^{N-1} \Big(\sqrt{\xi_i^2 + \eta_{j}^2}\Big)^{-(2-\veps)} \le  Ch^2 \int_h^L \int_h^L |\bxi|^{-(2-\veps)}{\,d\xi d\eta}  \le Ch^{2}.
\eea
Following similar lines, we can get 
\beas
|\, IV_2\,|, {|\,IV_3\,|, |\,IV_4\,|} \le Ch^2. 
\eeas
Using Property \ref{K-Property} (ii) and Lemma \ref{lemma1-2D} (iii),  we obtain 
\beas
 {|\, IV_5\,|, |\,IV_6\,|, |\,IV_7\,|, |\, IV_8\,|} \le Ch^{2+\veps}.
\eeas
To avoid redundancy, we will omit their {proofs}. Hence, we get 
\bea\label{term4}
|\,IV\,| \le Ch^2.
\eea
Combining (\ref{term1})--(\ref{term3}) and (\ref{term4}) yields (\ref{thm2-error}).

\bigskip
\noindent{\bf Case (ii)  (For $u\in C^{3, \ap-1+\veps}({\mathbb R}^2)$ with $\alpha \in [1,2)$): } 
By the extension of the weighted Montgomery identity in Lemma \ref{lemma2-2D} (iii) with  $\og(\bxi) = |\bxi|^{-\ap}$, we obtain 
\begin{eqnarray} \label{Error2}
&&e_{\alpha,2}^{h}({\bf x}) = -\fl{c_{2,\alpha}}{4}
\bigg[\int_{I_{00}} 
\Big(4\psi_{2}(\xi,\eta)-2 \big(\psi_{2}(\xi_1,0) + \psi_{2}(0,\eta_1)\big)\Big)|\bxi|^{-\alpha}\,d\xi d\eta\nn \\
&&\hspace{.6cm}-\fl{1}{2}\sum_{(i, j) \in \widetilde{\varkappa}_{N-1}} \int_{I_{ij}} \Big(\Theta_{I_{ij}}^{(1,0)} (\xi,\eta)\partial_{2,1}{\psi_2}(\xi,\eta) + \Theta_{I_{ij}}^{(0,1)} (\xi,\eta)\partial_{1,2}{\psi_2}(\xi,\eta)\Big)\,d\xi\,d\eta \nn\\
&&\hspace{0.6cm}+\sum_{(i, j) \in \widetilde{\varkappa}_{N-1}}\sum_{n=0,1}\left(-\fl{1}{2}\right)^{n}\bigg(
\int_{\xi_{i}}^{\xi_{i+1}}  \Big(\Theta_{I_{ij}}^{(1,n)}(\xi,\eta_{j+1})\partial_{2,n}{\psi_2}(\xi,\eta_{j+1})
-\Theta_{I_{ij}}^{(1,n)}(\xi,\eta_{j})\partial_{2, n}{\psi_2}(\xi,\eta_{j})\Big)d\xi \nn\\
&&\hspace{2.6cm}
+\int_{\eta_{j}}^{\eta_{j+1}} \Big(
\Theta_{I_{ij}}^{(n,1)}(\xi_{i+1},\eta)\partial_{n, 2}{\psi_2}(\xi_{i+1},\eta)- 
\Theta_{I_{ij}}^{(n,1)}(\xi_{i},\eta)\partial_{n, 2}{\psi_2}(\xi_{i},\eta)\Big)d\eta\bigg)\nn\\
&&\hspace{.6cm}-\sum_{(i, j) \in \widetilde{\varkappa}_{N-1}}\sum_{(m, n)\in\widetilde{\varkappa}_1}(-1)^{m+n}\Big(\Theta_{I_{ij}}^{(m,n)}(\bxi_{(i+1)j})\partial_{m,n}\psi_2(\bxi_{(i+1)j})- \Theta_{I_{ij}}^{(m,n)}(\bxi_{ij})\partial_{m,n}\psi_2(\bxi_{ij})\nn\\
&&\hspace{2.6cm}+\Theta_{I_{ij}}^{(m,n)}(\bxi_{i(j+1)})\partial_{m,n}\psi_2(\bxi_{i(j+1)})-\Theta_{I_{ij}}^{(m,n)}(\bxi_{(i+1)(j+1)})\partial_{m,n}\psi_2(\bxi_{(i+1)(j+1)})\Big) \nn\\ 
&&\hspace{1.4cm} = -\fl{c_{2,\alpha}}{4}\big(I + II + III +IV\big).
\end{eqnarray}
Here, the estimates of term $I$--$IV$ can be done by following similar lines as those in the proof of Case (i). 
For the purpose of brevity, we will omit their {proofs}.
\end{proof}

Theorem \ref{thm2} shows that our finite difference method can achieve the second order of accuracy for any $\ap \in (0, 2)$. 
For $\ap < 1$, less smoothness condition of $u$ is required to achieve this accuracy. 
As $\ap \to 2^-$, the behavior of our finite difference method that is stated in Theorem \ref{thm2} is consistent with the central difference method for the classical Laplace operator. 

% ====================================================================
\section{Generalization to three dimensions}
\label{section6}
\setcounter{equation}{0}

%%%% %%%% %%% %%%%
In this section, we will generalize our numerical methods in Section \ref{section2} and present a finite difference method for the 3D integral  fractional Laplacian (\ref{fL-nD}). 
We consider the fractional Laplacian $(-\Dt)^{\fl{\ap}{2}}$ on the bounded domain $\Omega = (a_x,b_x)\times(a_y,b_y)\times (a_z,b_z)$ with extended homogeneous boundary conditions. 

%%%% %%%% %%%% %%%%
Following the same lines as in Section \ref{section2}, we first rewrite the 3D  fractional Laplacian (\ref{fL-nD}) as a weighted integral, i.e.,
\begin{eqnarray}\label{fL2-3D}
(-\Delta)^{\fl{\ap}{2}}u({\bf x}) 
= -c_{3,\ap}\int _0^\infty\int_0^\infty\int_0^\infty
\psi_{\gamma}({\bf x},\bxi) \og_\gm(\bxi)\,d\xi  d\eta  d\zeta, \qquad \mbox{with} \quad \og_\gm(\bxi) = |\bxi|^{\gm-(3+\ap)}, 
\end{eqnarray}
where the vector {\boldmath$\xi$}$= (\xi,\,\eta,\,\zeta)$ with  $\xi = |x-x'|$, $\eta = |y-y'|$ and $\zeta = |z-z'|$, and the function
\begin{eqnarray}\label{psi-3D}
\psi_{\gamma}({\bf x},\bxi) =
\frac{1}{|\bxi|^{\gamma}}
\bigg(\sum_{m,n,s = 0, 1}u\big(x+(-1)^m\xi, \, y+(-1)^n\eta, \, z+(-1)^s\zeta \big)- 8 u({\bf x})\bigg). 
\end{eqnarray}
%%%% %%%% %%%% %%%%
Choose a constant  $\displaystyle L = \max \big\{b_x-a_x,\, b_y-a_y, \, b_z-a_z\big\}$. 
Define grid points $\xi_{i} = ih$, $\eta_j = jh$, $\zeta_k = kh$, for $0\leq i,j,k \leq N$, with the mesh size $h = L/N$. 
For notational convenience, we let  $\bxi_{ijk} = (\xi_{i},\eta_{j},\zeta_{k})$, 
and $|\bxi_{ijk}| = \sqrt{\xi_i^2 + \eta_j^2 + \zeta_k^2}$. 
Denote  $D_1 = (0, L)^3 = \cup_{i,j,k=0}^{N-1} I_{ijk}$ with  the element
\beas
I_{ijk}= \big[ih,\,(i+1)h\big]\times\big[jh,\,(j+1)h\big]\times \big[kh,\,(k+1)h\big],\qquad 0 \le i, j, k \le N-1,
\eeas
and $D_2 = \big({\mathbb R}^3_+\big)\backslash D_1$. 
The fractional Laplacian in (\ref{fL2-3D}) can be formulated as: 
\bea
&&(-\Dt)^{\fl{\ap}{2}} u(\bx)=-c_{3,\ap}\bigg(\int_{D_1}\psi_{\gamma}({\bf x},\bxi) \og_\gm(\bxi)\,d\xi d\eta d\zeta + \int_{D_2} \psi_{\gm}(\bx, \bxi) \og_\gm(\bxi)\,d\xi d\eta d\zeta \bigg)\nn\\
\label{eq3D}
&&\hspace{2cm}=-c_{3,\ap}\bigg(\sum_{i, j, k = 0}^{N-1}\int_{I_{ijk}}\psi_{\gamma}({\bf x},\bxi) \og_\gm(\bxi)\,d\xi d\eta d\zeta -8u(\bx) \int_{D_2} |\bxi|^{-(3+\ap)}d\xi d\eta d\zeta \bigg),\qquad\qquad\quad 
\eea
since the extended homogeneous Dirichlet boundary {condition implies} that $\psi_\gm(\bx, \bxi) = -8 u(\bx) {|\bxi|^{-\gm}}$ for any $\bx \in \Og$ and $\bxi \in D_2$.

%%%% %%%% %%%% %%%%
We now focus on the approximation to the integral over each element $I_{ijk}$.
If  $i+j+k \neq 0$, we use the weighted trapezoidal rules and obtain the approximation:
\begin{eqnarray}\label{Iijk-3D}
\int _{I_{ijk}}
\psi_{\gamma}({\bf x},\bxi) \og_\gm(\bxi)\,d\xi d\eta d\zeta\approx\,
\frac{1}{8} \bigg(\sum_{m, n, s = 0, 1} \psi_\gm\big(\bx, \bxi_{(i+m)(j+n)(k+s)}\big)\bigg) \int _{I_{ijk}} \og_\gm(\bxi)\,d\xi d\eta d\zeta. 
\eea
%%%% %%%% %%%% %%%%
If $i = j = k = 0$, using the weighted trapezoidal rules leads to 
\begin{eqnarray}\label{I000-3D}
\int _{I_{000}}
\psi_{\gamma}({\bf x},\bxi) \og_\gm(\bxi)\,d\xi d\eta d\zeta\approx 
\frac{1}{8}
\bigg(\lim_{\bxi \rightarrow {\bf 0}}\psi_{\gamma}({\bf x},\bxi)+\sum_{\substack{m, n, s = 0, 1 \\ 
m+n+s\neq0}} \psi_\gm\big(\bx, \bxi_{mns}\big)\bigg) \int _{I_{000}}\og_\gm(\bxi)\,d\xi d\eta d\zeta. 
\end{eqnarray}
Assuming the above limit exists,  we then divide our discussion into two parts based on the value of the splitting parameter $\gm$.  
For $\gamma = 2$, we can approximate the limit as:
\begin{eqnarray}\label{limit-gam2-3D}
\lim_{\bxi\rightarrow {\bf 0}}\psi_2({\bf x},\bxi) \approx
\frac{5}{3}\sum_{\substack{m,n,s = 0, 1\\ m+n+s = 1}}\psi_{2}\big({\bf x},\bxi_{mns}\big)-\sum_{\substack{m, n, s = 0, 1 \\ 
m+n+s > 1}}\psi_{2}\big({\bf x},\bxi_{mns}\big).
\end{eqnarray}
While  $\gamma \in (\alpha,2)$,  we obtain 
\begin{eqnarray}\label{limit-gam-3D}
	\lim_{\bxi \rightarrow {\bf 0}}\psi_{\gamma}({\bf x},\bxi)
	 = \lim_{\bxi \rightarrow {\bf 0}}\Big(\psi_2({\bf x},\bxi)\,|\bxi|^{2-\gamma}\Big) = 0.
\end{eqnarray}
Substituting the limits (\ref{limit-gam2-3D})--(\ref{limit-gam-3D}) into (\ref{I000-3D}) gives the approximation of the integral over $I_{000}$ as: 
\begin{eqnarray}\label{I000-final}
\int_{I_{000}}\psi_{\gamma}({\bf x},\bxi) \og(\bxi)\,d\xi d\eta d\zeta
\approx \frac{1}{8}\bigg(\sum_{\substack{m, n, s = 0, 1\\ m+n+s\neq 0}} c_{mns}^\gm 
\psi_\gm\big(\bx, \bxi_{mns}\big)\bigg)\int _{I_{000}}\og_\gm(\bxi)\,d\xi d\eta d\zeta,
\end{eqnarray}
where the coefficient 
\beas
c_{mns}^\gm = \left\{\begin{array}{ll}
1, \quad & \mbox{if  $\gm \in (\ap, 2)$},\\
8/3, & \mbox{if  $\gm = 2$, and $m + n + s = 1$,} \\
0, & \mbox{if  $\gm = 2$,  and $m + n + s > 1$.}
\end{array}\right.
\eeas
Combining (\ref{eq3D}) with (\ref{Iijk-3D}) and (\ref{I000-final}), we obtain 
\bea\label{eq-3D1}
&&(-\Dt)^{\fl{\ap}{2}}_{h, \gm} u(\bx) = -\fl{c_{3,\ap}}{8} \bigg[\sum_{\substack{i, j, k = 0\\ i+j+k\neq 0}}^{N-1}\bigg(\sum_{m, n, s = 0, 1} \psi_\gm\big(\bx, \bxi_{(i+m)(j+n)(k+s)}\big)\bigg)\int _{I_{ijk}} \og_\gm(\bxi)\,d\xi d\eta d\zeta\qquad\qquad \nn\\
&&\hspace{.8cm}+ \bigg(\sum_{\substack{m, n, s = 0, 1\\ m+n+s\neq 0}} c_{mns}^\gm \psi_\gm\big(\bx, \bxi_{mns}\big)\bigg)\int _{I_{000}}\og_\gm(\bxi) d\xi d\eta d\zeta - 64\,u(\bx)\int_{D_2}|\bxi|^{-(3+\ap)} d\xi d\eta d\zeta \bigg).
\eea

%%%% %%%% %%%% %%%%
Without loss of generality, we assume that $N_x = N$, and choose $N_y,\,N_z$  as the smaller integers such that  $a_y + N_yh \geq b_y$ and $a_z + N_zh \geq b_z$. 
Define the grid points $x_i = a_x + ih$ for $0 \le i \le N_x$,  $y_j = a_y + jh$ for $0 \le j \le N_y$, and $z_k = a_z + kh$ for $0 \le k \le N_z$. 
Let $u_{ijk}$ represent the numerical solution {of} $u(x_i,y_j,z_k)$. 
Combining (\ref{eq-3D1}) with (\ref{psi-3D}) and simplifying the calculations, we then obtain the finite difference schemes for the 3D fractional Laplacian $(-\Dt)^{\fl{\ap}{2}}$ as: 
\begin{eqnarray}\label{fLh-3D}
&&(-\Delta)_{h,\gamma}^{\frac{\alpha}{2}}u_{ijk}
= -c_{3,\alpha}\bigg[a_{000}\,u_{ijk} + 
\sum_{p=0, 1}\bigg(\sum_{m\in S_i^p}a_{m00}\,u_{[i+(-1)^pm]jk} 
+ \sum_{n\in S_j^p}a_{0n0}u_{i[j+(-1)^pn]k} \nn\\
&&\hspace{1.6cm} + \sum_{s\in S_k^p}a_{00s}u_{ij[k+(-1)^ps]} \bigg)
+\sum_{p, q=0, 1}\bigg(\sum_{s\in S_k^q}\sum_{n\in S_j^p}a_{0ns}u_{i[j+(-1)^pn][k+(-1)^qs]} \nn\\
&&\hspace{1.6cm}+ \sum_{s\in S_k^q}\sum_{m\in S_i^p}a_{m0s}u_{[i+(-1)^pm]j[k+(-1)^qs]}
+ \sum_{n\in S_j^q}\sum_{m\in S_i^p}a_{mn0}u_{[i+(-1)^pm][j+(-1)^qn]\,k} \bigg) \nn\\
&&\hspace{1.6cm}+ \sum_{p,q,r=0,1}\sum_{s\in S_k^r}\sum_{n\in S_j^q}\sum_{m\in S_i^p}a_{mns}u_{[i+(-1)^pm][j+(-1)^qn][k+(-1)^rs]} \bigg], 
\end{eqnarray}
for $1 \le i \le N_x -1$, $1 \le j \le N_y -1$, and $1 \le k \le N_z-1$, where the index sets 
\beas
&&S_i^p = \big\{l\,|\, l\in{\mathbb N}, \  1 \le i + (-1)^pl \le N_x-1\big\},\\
&&S_j^p = \big\{l\,|\, l\in{\mathbb N}, \  1 \le j + (-1)^pl \le N_y-1\big\}, \\
&&S_k^p = \big\{l\,|\, l\in{\mathbb N}, \  1 \le k + (-1)^pl \le N_z-1\big\},\qquad p = 0, \,\mbox{or} \ 1.
\eeas
For $m + n + s > 0$, the coefficient
\beas
a_{mns} = \frac{2^{\sigma(m, n, s)}}{8|\bxi_{mns}|^{\gamma}}\bigg(\int_{T_{mns}} |\bxi|^{\gamma-(3+\alpha)}\,d\bxi - \bar{c}_{mns}\left\lfloor \fl{\gm}{2}\right\rfloor\int_{I_{000}} |\bxi|^{\gamma-(3+\alpha)}\,d\bxi\bigg), 
\eeas
where $\sigma(m, n, s)$ denotes the number of zeros of $m, n$ and $s$. 
For $m, n, s \le 1$,  if $\sigma(m, n, s) = 2$, then $\bar{c}_{mns} = -\fl{5}{3}$;  if $\sigma(m, n, s) = 1$,  then $\bar{c}_{mns} = 1$;  otherwise,  $\bar{c}_{mns} = 0$ if $\sigma(m, n, s) = 0$.  
Here, we denote $T_{mns}$ (for $0 \le m, n, s \le N-1$) as the collection of all  elements  that have the point $\bxi_{mns}$ as a vertex, i.e.,
\beas
T_{mns} = \big([(m-1)h, \,(m+1)h] \times [(n-1)h,\,(n+1)h] \times [(s-1)h,\,(s+1)h]\big)\cap D_1. 
\eeas  
The coefficient $a_{000}$ is computed by:  
\begin{eqnarray*}
&&a_{000} = -2\sum_{m=1}^{N}\big(a_{m00} + a_{0m0} + a_{00m}\big) 
-4\sum_{m,n=1}^{N}\big(a_{0mn} + a_{m0n} + a_{mn0}\big) \nn\\
&&\hspace{1.2cm} - 8\sum_{m,n,s=1}^{N}a_{mns} -8\int_{D_2} |\bxi|^{-(3+\ap)}\,d\xi d\eta d\zeta. 
\end{eqnarray*}

%%%% %%%% %%%% %%%% Remark 
\begin{remark}
As in the 2D case,  the optimal choice of the splitting parameter is $\gm = 2$. 
Moreover,  as $\ap \to 2^-$, the finite difference scheme in (\ref{fLh-3D}) with $\gm = 2$ collapses to  the central difference scheme of the classical 3D Laplace operator $-\Dt = -(\p_{xx} + \p_{yy} + \p_{zz})$. \end{remark}

%%%% %%%% %%%% %%%% 
Following similar arguments in proving Theorems \ref{thm1} and \ref{thm2}, we can obtain the following estimates on the local truncation errors of the finite difference scheme  to the $d$-dimensional  ($d\ge 1$) fractional Laplacian $(-\Dt)^{\fl{\ap}{2}}$.  
Note that our estimates in the following theorems hold for any dimension $d \ge 1$ and greatly improve the 1D error estimates in \cite{Duo2018}. 
We will omit their proofs for brevity. 
\begin{theorem}[{\bf Minimum consistent conditions}]\label{thm3}
Suppose that $u$ has finite support on the domain $\Og \subset {\mathbb R}^d$ (for $d \ge 1$).
Let $(-\Dt)_{h,\gamma}^{\fl{\ap}{2}}$ be the finite difference approximation of the fractional Laplacian $(-\Dt)^{\fl{\ap}{2}}$, with $h$ a small mesh size.  
If $u \in C^{\lfloor\ap\rfloor,\,\ap-\lfloor\ap\rfloor+\veps}({\mathbb R}^d)$ with $0 < \veps \le  1 + \lfloor\ap\rfloor -\ap$, then for any splitting parameter $\gm \in (\ap, 2]$,  the local truncation error satisfies
\bea \label{thm3-error}
 \big\|(-\Dt)^{\fl{\ap}{2}} u  - (-\Dt)^{\fl{\ap}{2}}_{h,\gamma} u\big\|_\infty \le Ch^{\veps},  \qquad  \mbox{for} \ \ \ap\in(0, 2)
\eea
with $C$ a positive constant independent of $h$. 
\end{theorem}

\begin{theorem}[{\bf Second order of accuracy}]\label{thm4}
Suppose that $u$ has finite support on the domain $\Og \subset {\mathbb R}^d$ (for $d \ge 1$).
Let $(-\Dt)_{h,\gamma}^{\fl{\ap}{2}}$ be the finite difference approximation of the fractional Laplacian $(-\Dt)^{\fl{\ap}{2}}$, with $h$ a small mesh size.  
If the splitting parameter $\gamma = 2$ and $u \in C^{2+\lfloor\ap\rfloor,\,\ap-\lfloor\ap\rfloor+\veps}({{\mathbb R}^d})$ with $0 < \veps \le  1 + \lfloor\ap\rfloor -\ap$,  then the local truncation error satisfies
\bea\label{thm4-error}
 \big\|(-\Dt)^{\fl{\ap}{2}} u  - (-\Dt)^{\fl{\ap}{2}}_{h,\gamma} u\big\|_\infty \le Ch^2,  \qquad  \mbox{for} \ \ \ap\in(0, 2)
\eea
with $C$ a positive constant independent of $h$. 
\end{theorem}

%%%% %%%% %%%% %%%%
Denote the vector ${\bf u} = \big({\bf u}_{x,y,1},  {\bf u}_{x, y, 2} \, \ldots,\, {\bf u}_{x,y,N_z-1}\big)^T$. 
Here, the block vector \\${\bf u}_{x,y,k} = \big({\bf u}_{x,1,k}, {\bf u}_{x,2,k} \,\ldots,\,{\bf u}_{x, N_y-1, k}\big)$, with each block ${\bf u}_{x, j, k} = \big({u}_{1jk}, {u}_{2jk}, \,\ldots,\, u_{(N_x-1)jk}\big)$.
Then, the matrix-vector form of (\ref{fLh-3D}) is given by 
\bea\label{A3D}
(-\Dt)^{\fl{\ap}{2}}{\bf u} = A_3 {\bf u}.
\eea
Here,  $A_3$ is the matrix representation of the 3D fractional Laplacian, defined as:
\begin{eqnarray*}\label{A-3D}
{{{A}_3}}= 
\left(
\begin{array}{cccccc}
A_{x,y,0} & A_{x,y,1} & \ldots &  A_{x,y,N_z-3} & A_{x,y,N_z-2}  \\
A_{x,y,1}& A_{x,y,0} & A_{x,y,1}  &  \cdots & A_{x,y,N_z-3}  \\
\vdots & \ddots  & \ddots & \ddots & \vdots \\
A_{x,y,N_z-3} &  \ldots  & A_{x,y,1} & A_{x,y,0} & A_{x,y,1} \\
A_{x,y,{N_z-2}} & A_{x,y,N_z-3}  & \ldots & A_{x,y,1} & A_{x,y,0}
\end{array}
\right),
\end{eqnarray*}
where for $k = 0,1,\dots,N_z-2$, the block matrix
\begin{eqnarray*}\label{Ak}
{{A}}_{x,y,k} = 
\left(
\begin{array}{cccccc}
A_{x,0,k} & A_{x,1,k} & \ldots &  A_{x,N_y-3,k} & A_{x,N_y-2,k}  \\
A_{x,1,k} & A_{x,0,k} & A_{x,1,k}  &  \cdots & A_{x,N_y-3,k}  \\
\vdots & \ddots  & \ddots & \ddots & \vdots \\
A_{x,N_y-3,k} &  \ldots  & A_{x,1,k} & A_{x,0,k} & A_{x,1,k} \\
A_{x,N_y-2,k} & A_{x,N_2-3,k}  & \ldots & A_{x,1,k} & A_{x,0,k}
\end{array}
\right),\nn
\end{eqnarray*}
with 
\begin{eqnarray}\label{Aj-3D}
{{A}}_{x,j,k} = 
\left(
\begin{array}{cccccc}
a_{0jk} & a_{1jk} & \ldots &  a_{(N_x-3)jk} & a_{(N_x-2)jk}  \\
a_{1jk} & a_{0jk} & a_{1jk}  &  \cdots & a_{(N_x-3)jk}  \\
\vdots & \ddots  & \ddots & \ddots & \vdots \\
a_{(N_x-3)jk} &  \ldots  & a_{1jk} & a_{0jk} & a_{1jk} \\
a_{(N_x-2)jk} & a_{(N_x-3)jk}  & \ldots & a_{1jk} & a_{0jk}
\end{array}
\right),  \nn
\end{eqnarray}
for $j = 0,1,\dots,N_y-2$, and $k = 0,1,\dots,N_z-2$.    
Similar to the 2D case, the {matrix-vector} product in (\ref{A3D}) can be efficiently computed via the 3D FFT, and the computational cost is of ${\mathcal O}(M\log M)$, and the memory requirement is ${\mathcal O}(M)$ with $M = (N_x-1)(N_y-1)(N_z-1)$.

% ======================================================================
%%%% %%%% %%%% %%%%
\section{Numerical experiments}
\setcounter{equation}{0}
\label{section7}

%%%% %%%% %%%% %%%% 
In this section,  we on one hand test the numerical accuracy of the proposed finite difference methods to  verify our analytical results in Theorems \ref{thm1}--\ref{thm2} and \ref{thm3}--\ref{thm4}, and on the other hand apply the methods  to study the fractional Poisson problems as well as the  fractional Allen--Cahn equation.  
Unless otherwise stated, we will choose the splitting parameter $\gm = 2$ in this section.
% --------------------------------------------------------------------------------------------------------------------------
\subsection{Numerical accuracy}
\label{section7-1}

To test the accuracy of our method, we consider the function
\bea\label{fun-u}
u(\bx) = \left\{\begin{array}{ll}
\displaystyle\bigg(\prod_{i=1}^d \big(1-(x^{(i)})^2\big)\bigg)^{s},  \qquad & \mbox{if} \ \ \bx \in (-1, 1)^d, \\
0, & \mbox{otherwise}.
\end{array}\right.
\eea
for $s \ge 1$, where $\bx = (x^{(1)}, x^{(2)}, \cdots, x^{(d)})$. 
Define the error $\| e \|_p= \|(-\Dt)^{\fl{\ap}{2}}u - (-\Dt )_{h, \gm}^{\fl{\ap}{2}}u\|_p$. 
Since the exact solution of $(-\Dt)^{\fl{\ap}{2}}u$ is unknown, we will use the numerical solutions with very fine mesh size, i.e., $h = 2^{-12}$ for 2D and {$h = 2^{-8}$} for 3D,  as the ``exact" solution to compute  numerical errors. 
%In the following examples, we will not only verify our analytical results in Theorems \ref{thm1}--\ref{thm2} and \ref{thm3}--\ref{thm4} but  also study the effects of the splitting parameter $\gm$ on numerical accuracy. 

\vskip 8pt
%%%% %%%% %%%% %%%% 
\noindent {\bf Example 5.1.1 }  
In this example, we take $d  = 2$ and verify the analytical results in Theorem \ref{thm1}.  
For $\ap \in (0, 1)$, we choose $s = 1$ in (\ref{fun-u}),  i.e., $u \in C^{0,1}({\mathbb R}^2)$, or equivalently $u \in C^{0,\,\ap + \veps}({\mathbb R}^2)$ with $\veps = 1-\ap$. 
%%%% %%%% %%%% %%%%
Table \ref{Table1} presents the numerical errors $\|e \|_\infty$ and $\|e\|_2$ for various $\ap$. 
%%%% %%%% Table 1: Numerical errors of our method with s = 1 
\begin{table}[htb!]
\begin{center}
\begin{tabular}{|c|l r r r r r r|}
\hline
\diagbox[width=3em]{$\ap$}{$h$} && 1/16 & 1/32 & 1/64 & 1/128 & 1/256  & 1/512 \\
\hline
\multirow{4}{*}{$0.1$} & \multirow{2}{*}{$\|e\|_\infty$}& 5.704E-4
& 3.122E-4 & 1.690E-4 & 9.095E-5 & 4.870E-5 & 2.579E-5  \\
& & c.r. & 0.8694 & 0.8853 & 0.8940 & 0.9012 & 0.9170  \\
\cline{2-8}
&  \multirow{2}{*}{$\|e\|_2$}& 3.687E-4 & 1.424E-4 & 5.450E-5 & 2.075E-5 & 7.861E-6 & 2.946E-6  \\
& & c.r. & 1.3728 & 1.3856 & 1.3932 & 1.4003 & 1.4159 \\
\hline
\multirow{4}{*}{$0.4$} & \multirow{2}{*}{$\|e\|_\infty$}& 4.250E-3
& 2.858E-3 & 1.901E-3 & 1.257E-3 & 8.278E-4 & 5.393E-4 \\
%\cline{2-7}
& & c.r. & 0.5723 & 0.5886 & 0.5963 & 0.6027 & 0.6183  \\
\cline{2-8}
&  \multirow{2}{*}{$\|e\|_2$}&  2.514E-3 & 1.190E-3 & 5.588E-4 & 2.614E-4 & 1.217E-4 & 5.609E-5  \\
%\cline{2-7}
& & c.r. & 1.0793 & 1.0901 & 1.0962 & 1.1023 & 1.1177 \\
\hline
\multirow{4}{*}{$0.7$} & \multirow{2}{*}{$\|e\|_\infty$}& 1.350E-2
& 1.115E-2 & 9.115E-3 & 7.418E-3 & 6.011E-3 & 4.818E-3 \\
%\cline{2-7}
& & c.r. & 0.2757 & 0.2907 & 0.2973 & 0.3033 & 0.3191 \\
\cline{2-8}
&  \multirow{2}{*}{$\|e\|_2$}& 7.590E-3
     &  4.414E-3 & 2.550E-3 & 1.467E-3 & 8.410E-4 & 4.768E-4 \\
%\cline{2-7}
& & c.r. & 0.7820 & 0.7918 & 0.7972 & 0.8030 & 0.8187 \\
\hline
\end{tabular}
\caption{Numerical errors $\| e \|_p= \|(-\Dt)^{\fl{\ap}{2}}u - (-\Dt )_{h, \gm}^{\fl{\ap}{2}}u\|_p$, where $u$ is defined in (\ref{fun-u}) with $d = 2$ and $s = 1$, i.e., $u\in C^{0, 1}({\mathbb R}^2)$. }
\label{Table1}
\end{center}\vspace{-5mm}
\end{table}
It shows that for each $\ap$, the accuracy in $\infty$-norm is ${\mathcal O}(h^{1-\ap})$, confirming our theoretical results in  Theorem \ref{thm1}. 
Moreover,  we find that the accuracy in $2$-norm is ${\mathcal O}(h^{\fl{3}{2}-\ap})$,  and the numerical errors increase with the value of $\ap$. 

%%%% %%%% %%%% %%%% 
On the other hand, Table \ref{Table2} shows the numerical errors for $\ap \in [1, 2)$, where $u$ is defined  in (\ref{fun-u}) with $s = 2$, i.e., $u \in C^{1,1}(\mathbb{R}^2)$, or equivalently $u \in C^{1,\alpha-1+\veps}(\mathbb{R}^2)$ with $\veps = 2-\alpha$.  
%%%% %%%% Table 2: Numerical errors of our method with s = 2 
\begin{table}[htb!]
\begin{center}
\begin{tabular}{|c|l r r r r r r|}
\hline
\diagbox[width=3em]{$\ap$}{$h$} && 1/16 & 1/32 & 1/64 & 1/128 & 1/256  & 1/512 \\
\hline
\multirow{4}{*}{$1$} & \multirow{2}{*}{$\|e\|_\infty$}& 3.017E-3
& 1.434E-3  & 6.584E-4 & 3.063E-4 & 1.452E-4 & 6.923E-5 \\
& & c.r. & 1.0732 & 1.1232 & 1.1041 & 1.0770 & 1.0682 \\
\cline{2-8}
&  \multirow{2}{*}{$\|e\|_2$}& 2.094E-3 & 6.972E-4 & 2.293E-4 & 7.549E-5 & 2.506E-5 & 8.350E-6  \\
& & c.r. & 1.5869 & 1.6044 & 1.6029 & 1.5908 & 1.5856  \\
\hline
\multirow{4}{*}{$1.4$} & \multirow{2}{*}{$\|e\|_\infty$}& 8.871E-3 
& 5.665E-3 & 3.505E-3  & 2.202E-3 & 1.405E-3 & 8.973E-4 \\
%\cline{2-7}
& & c.r. & 0.6470 & 0.6927 & 0.6708 & 0.6481 & 0.6467  \\
\cline{2-8}
&  \multirow{2}{*}{$\|e\|_2$}& 5.568E-3 & 2.270E-3 & 9.677E-4 & 4.243E-4 & 1.898E-4 & 8.532E-5 \\
%\cline{2-7}
& & c.r. & 1.2945 & 1.2301 & 1.1896 & 1.1605 & 1.1535 \\
\hline
\multirow{4}{*}{$1.9$} & \multirow{2}{*}{$\|e\|_\infty$}& 
    1.224E-2 & 7.609E-3 & 7.504E-3 & 6.956E-3 & 6.395E-3 & 5.828E-3 \\
%\cline{2-7}
& & c.r. & 0.6854 & 0.0199 & 0.1094 & 0.1213 & 0.1340 \\
\cline{2-8}
&  \multirow{2}{*}{$\|e\|_2$}& 1.289E-2 & 3.973E-3 & 1.947E-3 & 1.204E-3 & 7.766E-4 & 5.000E-4 \\
%\cline{2-7}
& & c.r. & 1.6975 & 1.0288 & 0.6939 & 0.6323 & 0.6352 \\
\hline
\end{tabular}
\caption{Numerical errors $\| e \|_p= \|(-\Dt)^{\fl{\ap}{2}}u - (-\Dt )_{h, \gm}^{\fl{\ap}{2}}u\|_p$,  where $u$ is defined in (\ref{fun-u}) with $d = 2$ and $s = 2$, i.e., $u\in C^{1, 1}({\mathbb R}^2)$.}
\label{Table2}
\end{center}\vspace{-5mm}
\end{table}
We find that our method has an accuracy of ${\mathcal O}(h^{2-\ap})$ in $\infty$-norm and ${\mathcal O}(h^{\fl{5}{2}-\ap})$ in $2$-norm. 
Our extensive simulations show that the same conclusions can be made for the 3D cases, and here we will omit showing them for the purpose of brevity. 

\vskip 8pt
%%%% %%%% %%%% %%%% 
\noindent{\bf Example 5.1.2  } 
In this example, we verify the second order of accuracy of our methods. 
Figure \ref{Figure2} presents the numerical errors $\|e\|_\infty$ versus different mesh size $h$, where the function $u$ is defined in (\ref{fun-u}) with $s = 2+\ap+\veps$ and $\veps = 0.1$, equivalently,  $u \in C^{2, \ap+\veps}({\mathbb R}^d)$ for $\ap < 1$ and $u \in C^{3, \ap-1+\veps}({\mathbb R}^d)$ for $\ap \ge 1$.
%%%% %%%% %%%% %%%% Figure 2
\begin{figure}[htb!]
\centerline{
(a) \includegraphics[height=4.80cm,width=6.660cm]{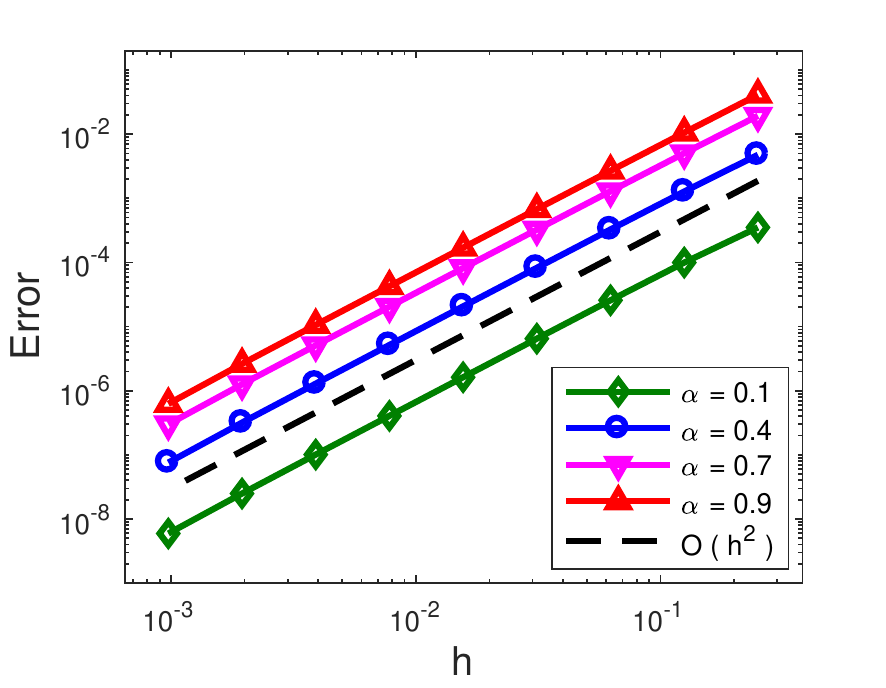}
(b) \includegraphics[height=4.80cm,width=6.660cm]{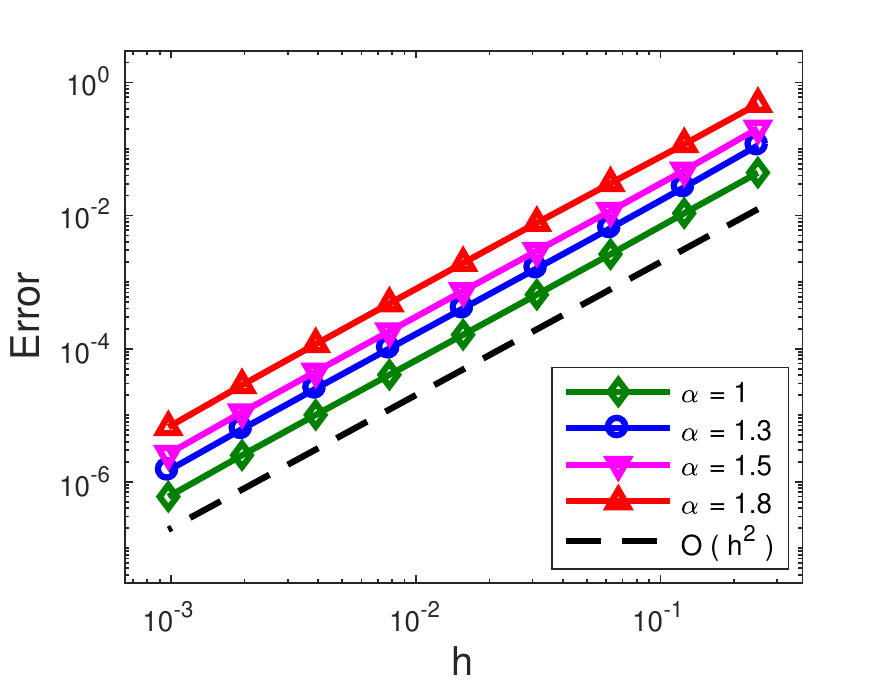}}
\centerline{
(c) \includegraphics[height=4.80cm,width=6.660cm]{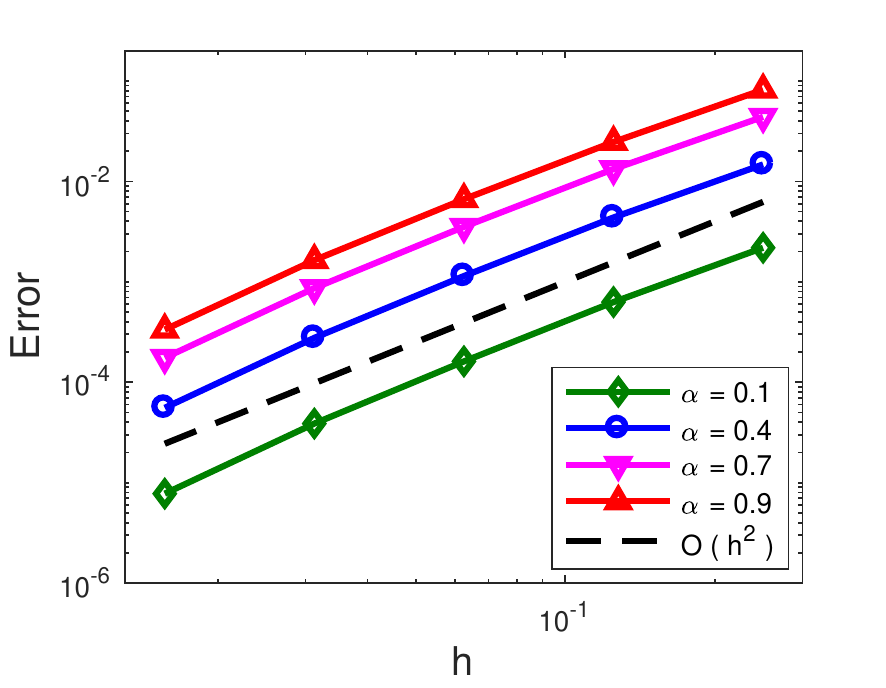}
(d) \includegraphics[height=4.80cm,width=6.660cm]{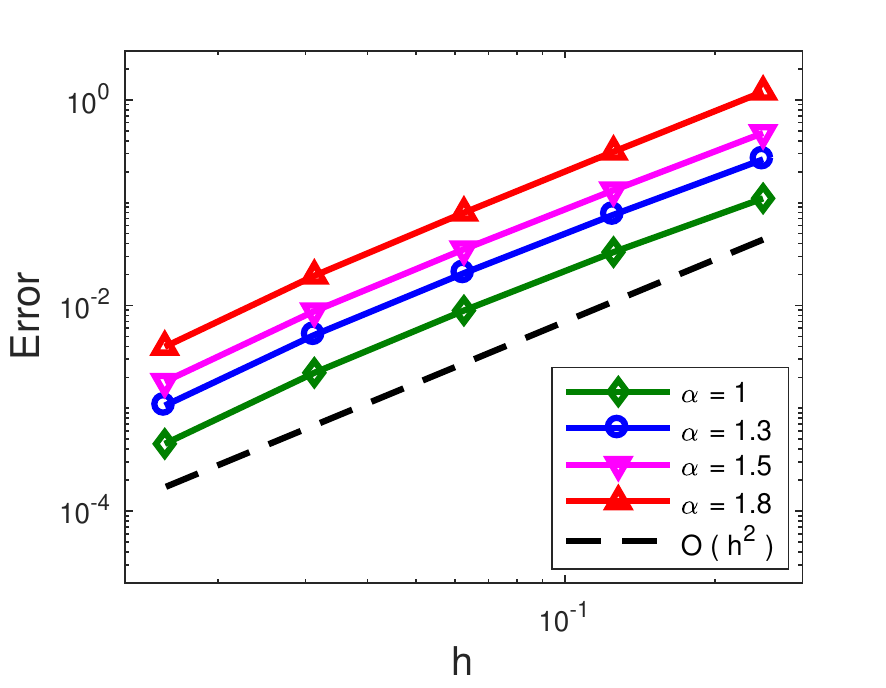}}
\caption{Numerical errors $\|e\|_\infty$ for $u$  in (\ref{fun-u}) with $s = 2+\ap+\veps$ and $\veps = 0.1$, i.e.,   $u \in C^{2, \ap+\veps}({\mathbb R}^d)$ for $\ap\in(0, 1)$ and $u \in C^{3, \ap-1+\veps}({\mathbb R}^d)$ for $\ap\in[1, 2)$, where $d = 2$ (a \& b) and $d = 3$ (c \& d).} \label{Figure2}
\end{figure}
It shows that when the maximum regularity conditions are satisfied, our method has the second order of accuracy  in both $\infty$-norm and $2$-norm, independent of the parameter $\ap$. This observation confirms our analytical results in Theorems \ref{thm2} and \ref{thm4}. 
Moreover, we find that the smaller the parameter $\ap$, the smaller the numerical errors. 

%%%% %%%% %%%% %%%% 
In addition, we compare the numerical errors of our methods for different choices of the splitting parameter $\gm$ in Figure \ref{Figure3}. 
\begin{figure}[htb!]
\centerline{
(a) \includegraphics[height=4.80cm,width=6.660cm]{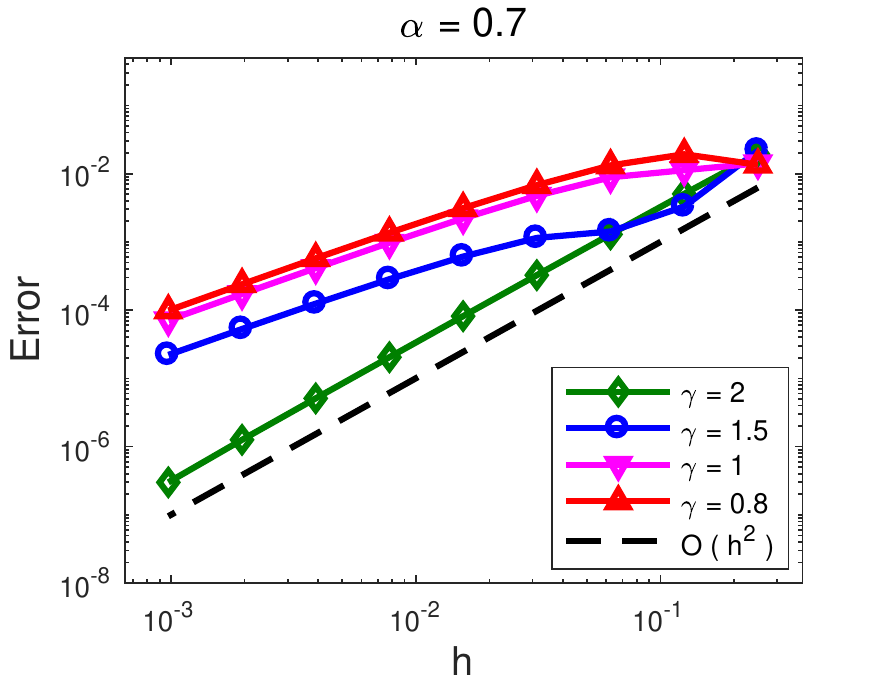}
(b) \includegraphics[height=4.80cm,width=6.660cm]{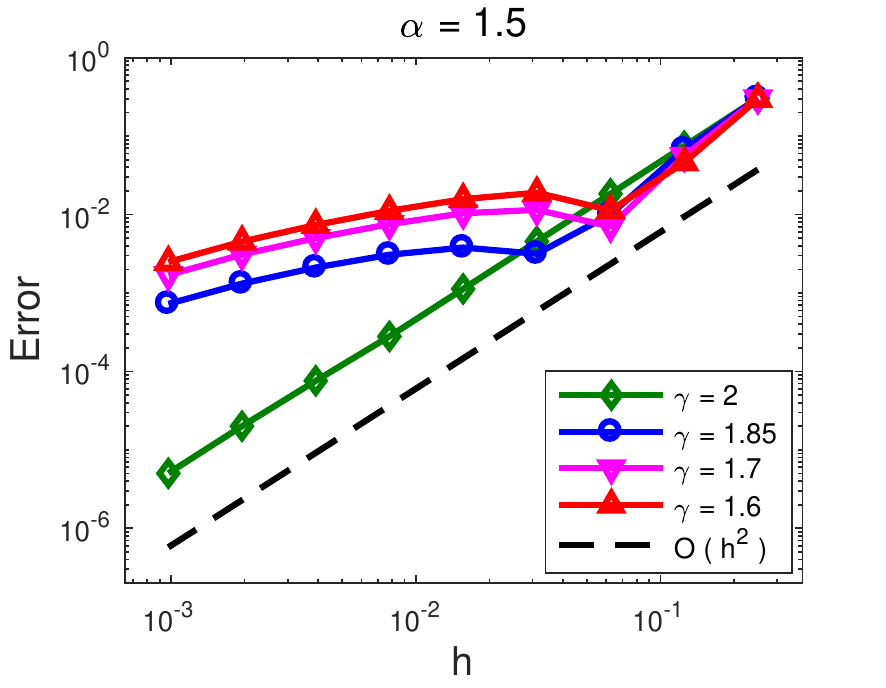}
}
\caption{Numerical errors $\|e\|_\infty$ for $u$  in (\ref{fun-u}) with $s = 2+\ap+\veps$ and $\veps = 0.1$ for various $\gm$.} \label{Figure3}
\end{figure}
It shows that the splitting parameter $\gm$ plays a crucial role in determining the numerical accuracy. 
Among all the choices of $\ap < \gm \le 2$, only the case of $\gm = 2$ leads to the optimal accuracy, i.e., ${\mathcal O}(h^2)$, while the other choices of $\gm \in (\ap, 2)$ yield an $\ap$-dependent rate, e.g., ${\mathcal O}(h^{1.3})$ for $\ap = 0.7$ and ${\mathcal O}(h^{0.5})$ for $\ap = 1.5$.

% --------------------------------------------------------------------------------------------------------------------------
\subsection{Fractional Poisson equations}
\label{section7-2}

%%%% %%%% %%%% %%%%
In this section, we test the performance of our finite difference methods in solving the fractional Poisson problem. 
Without loss of generality, we will consider the 2D fractional Poisson equation:  
\bea\label{Poisson}
(-\Dt)^{\fl{\ap}{2}}u = f(\bx), &\ & \bx \in (-1, 1)^2,\\ 
\label{BCP}
u(\bx) = 0, && \bx \in {\mathbb R}^2\backslash (-1, 1)^2.
\eea
Our extensive studies show that the same conclusions can be obtained in solving the  3D Poisson problems. 
Denote the numerical errors $\|e_u\|_p = \|u - u_h\|_p$,  where $u$ and $u_h$ represent the exact and numerical solutions of (\ref{Poisson})--(\ref{BCP}), respectively.

\vskip 8pt
\noindent{\bf Example 5.2.1.}  We solve the problem (\ref{Poisson})--(\ref{BCP}) with the exact solution 
\bea\label{uu2}
\label{u2}u(\bx) = [(1-x^2)(1-y^2)]^{2}, \qquad \bx\in (-1, 1)^2.
\eea
In practice, the function $f$ in (\ref{Poisson})  is prepared numerically with a fine mesh size $h =  2^{-12}$, i.e., computing $f = (-\Dt)_{h,\gm}^{\fl{\ap}{2}}u$ with $u$ defined in (\ref{uu2}).

In Table \ref{Tab1}, we present both the $\infty$-norm and $2$-norm errors for various $\ap$. 
It shows that even though the solution satisfies $u \in C^{1, 1}({\mathbb R}^2)$, our method can achieve the accuracy ${\mathcal O}(h^2)$, independent of the parameter $\ap$. 
%%%% %%%% Table 1: Numerical errors of our method with s = 2 and \gamma = 2
\begin{table}[htb!]
\begin{center}
\begin{tabular}{|c|l r r r r r r|}
\hline
\diagbox[width=3em]{$\ap$}{$h$} && 1/16 & 1/32 & 1/64 & 1/128 & 1/256  & 1/512 \\
\hline
\multirow{4}{*}{$0.4$} & \multirow{2}{*}{$\|e_u\|_\infty$}& 1.658E-4 &  4.135E-5  &  1.031E-5 &  2.571E-6 &  6.399E-7 & 1.580E-7 \\
& & c.r. & 2.0031 & 2.0040 & 2.0039 & 2.0060 & 2.0182 \\
\cline{2-8}
&  \multirow{2}{*}{$\|e_u\|_2$}& 1.288E-4 & 3.332E-5 & 8.547E-6 &  2.172E-6 &  5.470E-7 & 1.360E-7  \\
& & c.r. & 1.9504 & 1.9627 & 1.9767 & 1.9892 & 2.0080   \\
\hline
\multirow{4}{*}{$1$} & \multirow{2}{*}{$\|e_u\|_\infty$}& 6.547E-4 & 1.584E-4 & 3.872E-5 & 9.520E-6 & 2.347E-6 & 5.754E-7 \\
%\cline{2-7}
& & c.r. & 2.0468 & 2.0331 & 2.0239 & 2.0202 & 2.0282  \\
\cline{2-8}
&  \multirow{2}{*}{$\|e_u\|_2$}& 5.134E-4 & 1.248E-4 & 3.123E-5 & 7.885E-6 & 1.989E-6 & 4.966E-7  \\
%\cline{2-7}
& & c.r. & 2.0400 & 1.9989 & 1.9858 & 1.9870 & 2.0019 \\
\hline
\multirow{4}{*}{$1.4$} & \multirow{2}{*}{$\|e_u\|_\infty$}& 1.184E-3
& 2.790E-4 & 6.650E-5 & 1.600E-5 & 3.872E-6 & 9.349E-7  \\
%\cline{2-7}
& & c.r. & 2.0856 & 2.0687 & 2.0553 & 2.0469 & 2.0502 \\
\cline{2-8}
&  \multirow{2}{*}{$\|e_u\|_2$}&  1.014E-3 & 2.291E-4 & 5.349E-5 & 1.281E-5 &  3.120E-6 & 7.622E-7 \\
%\cline{2-7}
& & c.r. & 2.1460 & 2.0988 &  2.0616 & 2.0380 & 2.0334 \\
\hline
\end{tabular}
\caption{Numerical errors in solving the fractional Poisson problem (\ref{Poisson})--(\ref{BCP}), where $f$ is chosen such that the exact solution is $u(\bx) = (1-x^2)^2(1-y^2)^2$.}
\label{Tab1}
\end{center}\vspace{-5mm}
\end{table}
In other words, to obtain the second order of accuracy,  the regularity of solution $u$ that is required  in solving the fractional Poisson problem is much lower than that required in discretizing the operator $(-\Dt)^{\fl{\ap}{2}}$, which is consistent with the central difference scheme for the classical Poisson problems. 
%%%% %%%% %%%% %%%%
From the above results and our extensive studies, we conclude that our method has the second order of accuracy  in solving the fractional Poisson problem, if the solution satisfy {\it at most}  $u\in C^{1,1}({\mathbb R}^d)$. 

\bigskip
%%%% %%%% %%%% %%%% Example 4
\noindent {\bf Example 5.2.2.}  We solve the problem (\ref{Poisson})--(\ref{BCP}) with the exact solution 
\bea\label{u2}u(\bx) = (1-x^2)(1-y^2), \qquad \bx\in (-1, 1)^2,
\eea
i.e., the solution $u \in C^{0, 1}({\mathbb R}^2)$.   
In this case, the solution has less regularity than that in the previous example. 
The  function $f$ in (\ref{Poisson}) is  similarly computed with the fine mesh size $h = 2^{-12}$. 
Table \ref{Tab2} shows the numerical errors $\|e_u\|_\infty$ and $\|e_u\|_2$ for various $\ap$.
%%%% %%%% Table 1: Numerical errors of our method with s = 2 and \gamma = 2
\begin{table}[htb!]
\begin{center}
\begin{tabular}{|c|c| r r r r r r|}
\hline
\diagbox[width=3em]{$\ap$}{$h$} && 1/16 & 1/32 & 1/64 & 1/128 & 1/256  & 1/512 \\
\hline
\multirow{4}{*}{$0.4$} & \multirow{2}{*}{$\|e\|_\infty$}& 1.551E-3 &  8.011E-4 & 4.071E-4 & 2.051E-4 & 1.026E-4 & 5.075E-5  \\
%\cline{2-7}
& & c.r. & 0.9533 & 0.9766 & 0.9893 & 0.9988 & 1.0160   \\
&  \multirow{2}{*}{$\|e\|_2$}& 1.198E-3 & 4.536E-4
 & 1.674E-4 & 6.079E-5 & 2.181E-5 & 7.699E-6  \\
%\cline{2-7}
& & c.r. & 1.4005 & 1.4380&  1.4616 & 1.4792 & 1.5019 \\
\hline
\multirow{4}{*}{$1$} & \multirow{2}{*}{$\|e\|_\infty$}& 2.424E-3 & 1.265E-3 & 6.456E-4 & 3.259E-4 & 1.632E-4 & 8.070E-5  \\
%\cline{2-7}
& & c.r. & 0.9389 & 0.9701 & 0.9863 & 0.9975 & 1.0162  \\
&  \multirow{2}{*}{$\|e\|_2$}& 2.492E-3 & 1.053E-3 & 4.258E-4 & 1.673E-4 &  6.433E-5 & 2.416E-5 \\
%\cline{2-7}
& & c.r. & 1.2433 &1.3056 &1.3476 &1.3791 &1.4130   \\
\hline
\multirow{4}{*}{$1.4$} & \multirow{2}{*}{$\|e\|_\infty$}& 1.982E-3
& 1.040E-3 & 5.326E-4 & 2.692E-4 & 1.349E-4 & 6.668E-5 \\
%\cline{2-7}
& & c.r. & 0.9301 & 0.9659 &  0.9844 & 0.9968 & 1.0165 \\
&  \multirow{2}{*}{$\|e\|_2$}& 2.583E-3 & 1.209E-3 & 5.397E-4 & 2.336E-4 & 9.886E-5 & 4.087E-5 \\
%\cline{2-7}
& & c.r. & 1.0954 & 1.1632 & 1.2081 & 1.2407 & 1.2744 \\
\hline
\end{tabular}
\caption{Numerical errors in solving the fractional Poisson problem (\ref{Poisson})--(\ref{BCP}), where $f$ is chosen such that the exact solution is $u(\bx) = (1-x^2)(1-y^2)$.}
\label{Tab2}
\end{center}
\end{table}

%%%% %%%% %%%% %%%%
It shows that the accuracy in $\infty$-norm is ${\mathcal O}(h)$, independent of the values of $\ap$, while the accuracy in 2-norm is ${\mathcal O}(h^p)$ for $1 \le p \le \fl{3}{2}$:  the smaller the parameter $\alpha$, the larger the accuracy rate $p$. 
For example, we have  $\|e\|_2 \sim O(h)$  for $\ap = 1.9$, while $\|e\|_2 \sim O(h^{\fl{3}{2}})$ for $\ap = 0.1$.  
Additionally, comparing the results in Table \ref{Tab2} with those in Table \ref{Tab1}, we find that the numerical errors increase when the solution $u$ is less smooth. 

\bigskip
%%%% %%%% %%%% %%%% Example 5:
\noindent{\bf Example 5.2.3 }  We solve the problem (\ref{Poisson})--(\ref{BCP}) with $f = 1$. 
It is well-known that the solution of the 1D fractional Poisson equation with $f = 1$ can be found analytically, which satisfies $u \in C^{0, \fl{\ap}{2}}({\mathbb R})$ \cite{Dyda2012,  Duo2018}. 
By contrast, the exact solution  on the 2D rectangular domain still remains unknown. 
Figure \ref{Fignew} illustrates the numerical solution $u$ for different $\ap$.
%%%% %%%% Figure u for Poisson f = 1.
\begin{figure}[htb!]
\centerline{
\includegraphics[height=3.80cm,width=4.660cm]{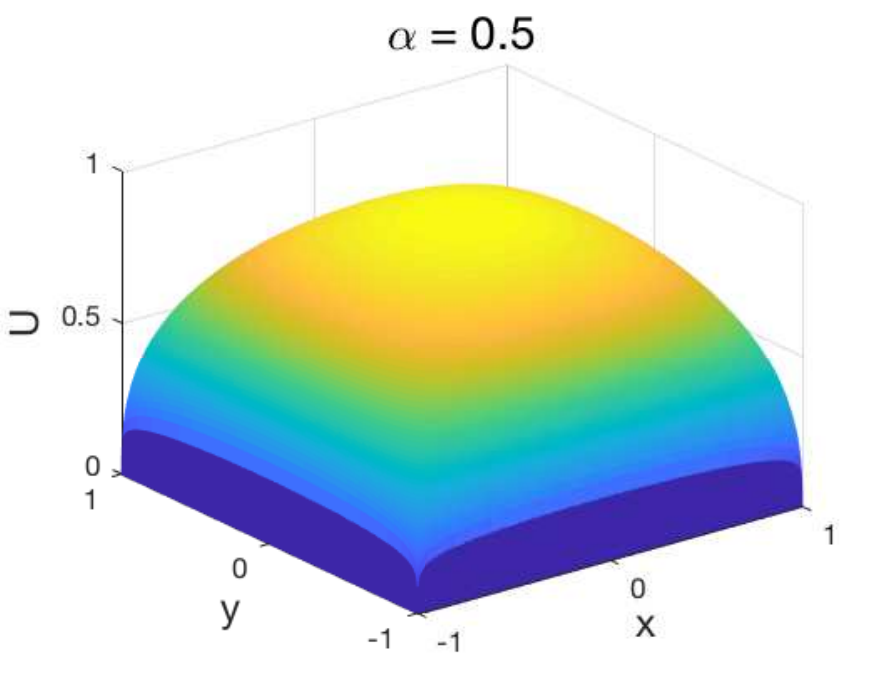}
\includegraphics[height=3.80cm,width=4.660cm]{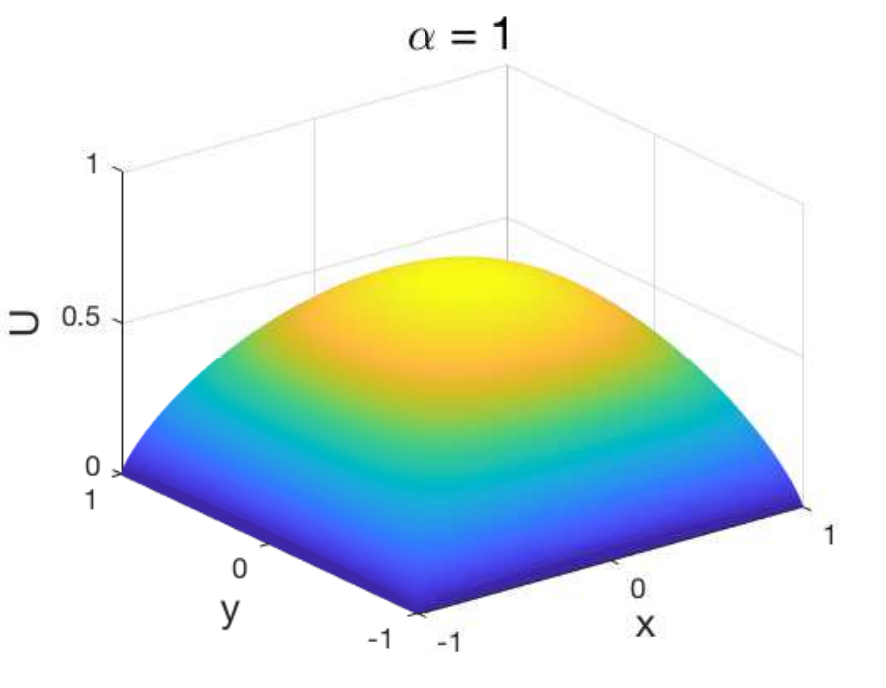}
\includegraphics[height=3.80cm,width=4.660cm]{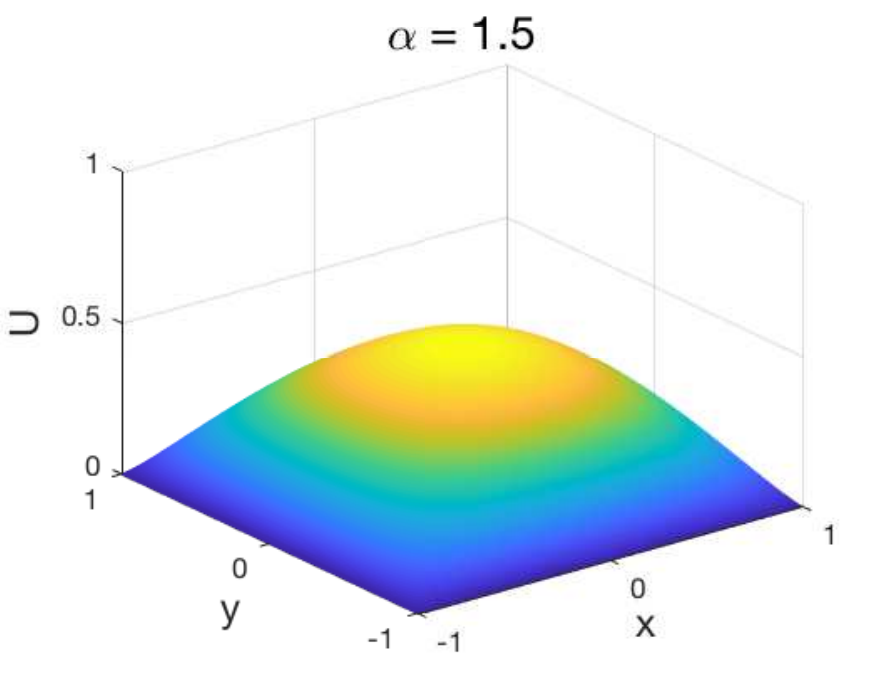}
}
\caption{Numerical solution to the fractional Poisson problem (\ref{Poisson})--(\ref{BCP}) with $f = 1$.} \label{Fignew}
\end{figure}
It shows that the larger the power $\ap$, the smoother the function near the boundary, which is consistent with the observations in 1D cases  \cite{Dyda2012,  Duo2018}. 

%%%% %%%% %%%% %%%%
Table \ref{Tab3} shows the numerical errors $\|e_u\|_\infty$ and $\|e_u\|_2$ for various $\ap$, where we use the numerical solution computed with fine mesh size $h = 2^{-12}$ as our ``exact" solution.
It shows that even though the regularity of the solution in this case is much lower than those in Examples 5.2.1--5.2.2, our method performs effectively in solving the fractional Poisson problem. 
%%%% %%%% Table 1: Numerical errors of our method with s = 1 and \gamma = 2
\begin{table}[htb!]
\begin{center}
\begin{tabular}{|c|l r r r r r r|}
\hline
\diagbox[width=3em]{$\ap$}{$h$} && 1/16 & 1/32 & 1/64 & 1/128 & 1/256  & 1/512 \\
\hline
\multirow{4}{*}{$0.5$} & \multirow{2}{*}{$\|e_u\|_\infty$}& 3.435E-2 &  2.645E-2  & 2.157E-2 & 1.814E-2 &  1.525E-2 & 1.282E-2 \\
%\cline{2-7}
& & c.r. & 0.3772 & 0.2941 & 0.2502 & 0.2501 & 0.2500\\
&  \multirow{2}{*}{$\|e_u\|_2$}& 2.754E-2 & 1.699E-2 & 1.034E-2 & 6.241E-3 & 3.748E-3 &  2.243E-3 \\
%\cline{2-7}
& & c.r. & 0.6970 & 0.7164 & 0.7283 & 0.7357 & 0.7404  \\
\hline
\multirow{4}{*}{$1$} & \multirow{2}{*}{$\|e_u\|_\infty$}&  1.631E-2 & 1.155E-2 & 8.173E-3 & 5.783E-3 & 4.091E-3 &  2.893E-3\\
%\cline{2-7}
& & c.r. & 0.4984 & 0.4985 & 0.4991 & 0.4995 & 0.4997 \\
&  \multirow{2}{*}{$\|e_u\|_2$}& 1.814E-2 & 1.015E-2 & 5.563E-3 & 3.007E-3 & 1.609E-3 & 8.537E-4  \\
%\cline{2-7}
& & c.r. & 0.8388 & 0.8669 & 0.8873 & 0.9025 & 0.9141 \\
\hline
\multirow{4}{*}{$1.5$} & \multirow{2}{*}{$\|e_u\|_\infty$}& 4.727E-3
& 2.812E-3 & 1.675E-3 & 9.967E-4 & 5.930E-4 & 3.527E-4 \\
%\cline{2-7}
& & c.r. & 0.7493 &0.7478 &0.7485 &0.7491 &0.7495 \\
&  \multirow{2}{*}{$\|e_u\|_2$}& 7.105E-3 & 3.737E-3 & 1.940E-3 & 9.966E-4 & 5.080E-4 & 2.575E-4 \\
%\cline{2-7}
& & c.r. & 0.9269 &0.9458 &0.9610 &0.9722 &0.9803 \\
\hline
\end{tabular}
\caption{Numerical errors in solving the fractional  Poisson problem (\ref{Poisson})--(\ref{BCP}) with $f = 1$.}
\label{Tab3}
\end{center}\vspace{-5mm}
\end{table}
We find that the accuracy in $\infty$-norm is ${\mathcal O}(h^{\fl{\ap}{2}})$ for any $\ap \in (0, 2)$, and the maximum errors occur around the boundary of the domain.  While the accuracy in $2$-norm is ${\mathcal O}(h^{\min\{1, \fl{1}{2}+\fl{\ap}{2}\}})$, that is,  ${\mathcal O}(h)$ for $\ap \ge 1$.  

% ----------------------------------------------------------------------------------------------------------------------------
\subsection{Fractional Allen--Cahn equation}
\label{section7-3}

%%%% %%%% %%%% %%%% 
The Allen--Cahn equation has been widely used in modeling phase field problems arising in materials science and fluid dynamics. 
Recently, the fractional analogue of the Allen--Cahn equation was proposed to study phase transition in the presence of anomalous diffusion \cite{Song2016}.  
Here, we apply our method to  study the benchmark problem -- coalescence of two ``kissing" bubbles -- in the phase field models. 
%%%% %%%% %%% %%%%
Consider the fractional Allen--Cahn equation \cite{Song2016}: 
\begin{eqnarray}\label{AC-1}
	&& \p_tu(\bx, t) = -(-\Delta)^{\fl{\alpha}{2}}u - \frac{1}{\dt^\ap}\,u(u^2-1), \qquad \bx \in \Omega, \quad t>0,\qquad\qquad \\ \label{AC-2}
	&& u(\bx, t) = -1, \qquad \bx \in \Omega^{c}, \quad t\geq 0, 
\end{eqnarray}
where the domain $\Og = (0, 1)^d$, and $u$ is the phase field function. 
The constant $\dt > 0$ describes the diffuse interface width.  
The initial condition is chosen  as 
\begin{eqnarray}\label{u0-two_bubble}
	u(\bx, 0) = 
		1-\tanh\bigg(\frac{d(\bx, \bx_1)}{{\dt}} \bigg) - \tanh\bigg(\frac{d(\bx, \bx_2)}{{\dt}} \bigg),  \qquad \bx \in {\mathbb R}^d, 
\end{eqnarray}
with the function $d(\bx, \bx_i) = |\bx -  \bx_i| - {0.12}$. 
Initially, two bubbles, centered at $\bx_1$ and $\bx_2$, respectively, are osculating or ``kissing". 
Note that the boundary condition in (\ref{AC-2}) is nonzero constant. 
Letting $\bar{u} = u + 1$, we can rewrite the problem  (\ref{AC-1}) as an equation of $\bar{u}$ with the extended homogeneous boundary conditions. 
Here, we discretize the problem (\ref{AC-1})--(\ref{u0-two_bubble}) by our finite difference methods in space and the Crank--Nicolson method in time. 
In each iteration step, we use the the conjugate gradient (CG) method coupled with the FFT-based algorithms for matrix-vector multiplication to solve the linear system. 

%%%% Figure 4
\begin{figure}[htb!]
\centerline{
\includegraphics[height=3.4cm,width=3.7cm]{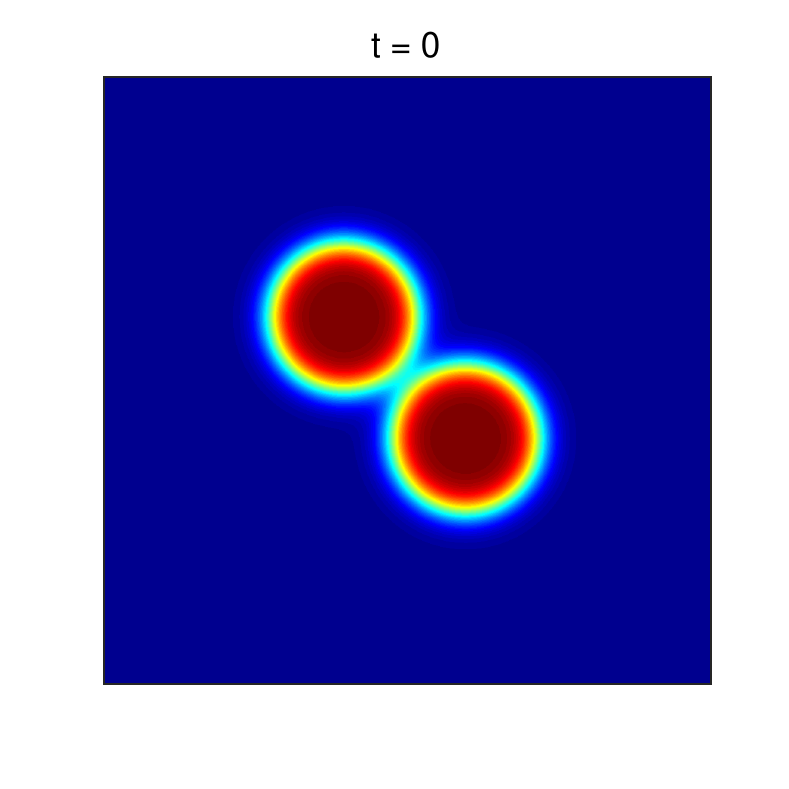}\hspace{-5mm}
\includegraphics[height=3.4cm,width=3.7cm]{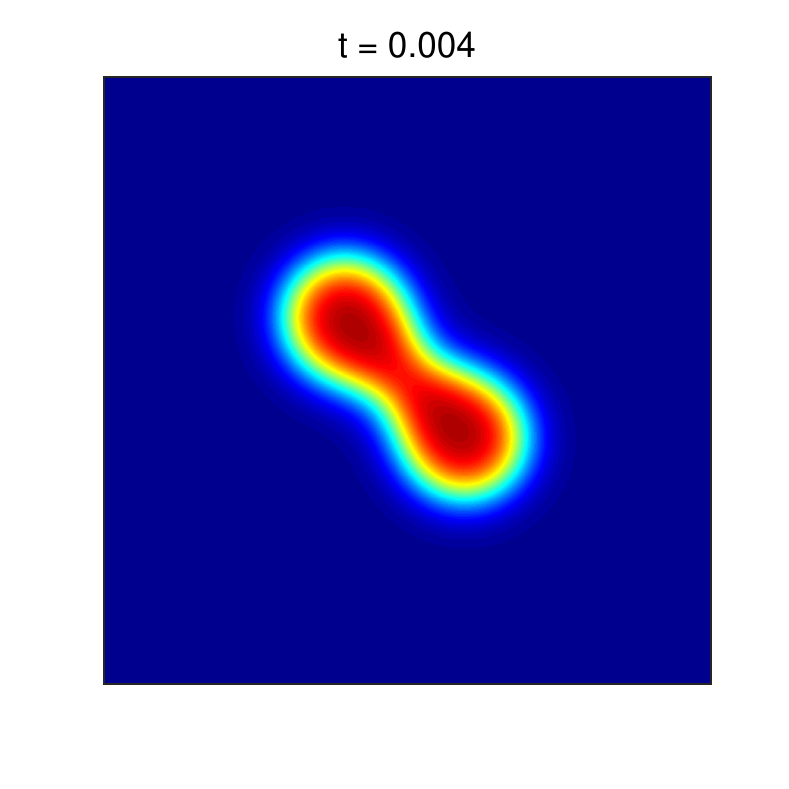}\hspace{-5mm}
\includegraphics[height=3.4cm,width=3.7cm]{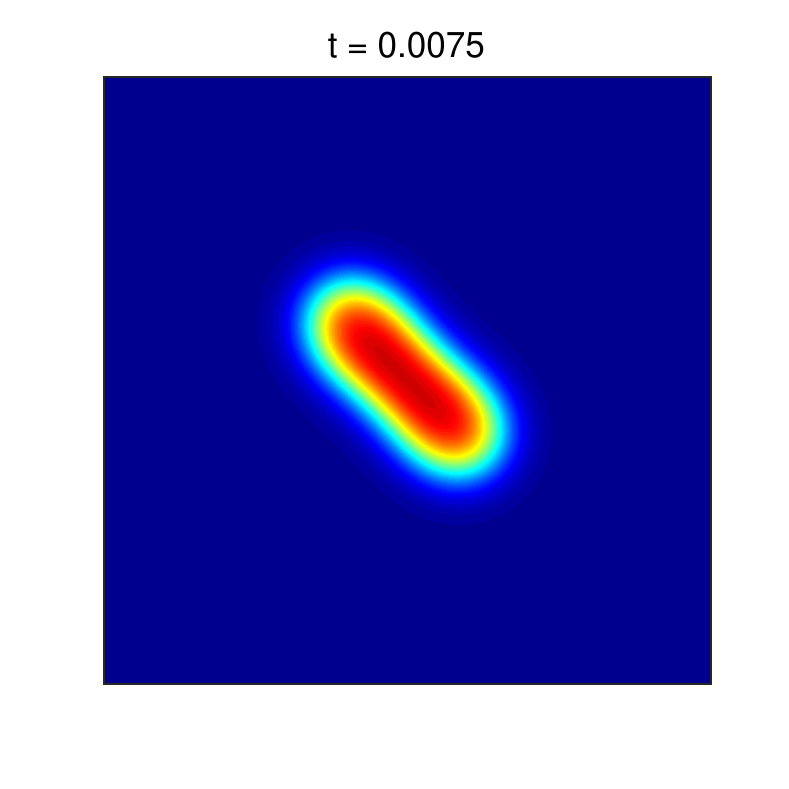}\hspace{-5mm}
\includegraphics[height=3.4cm,width=3.7cm]{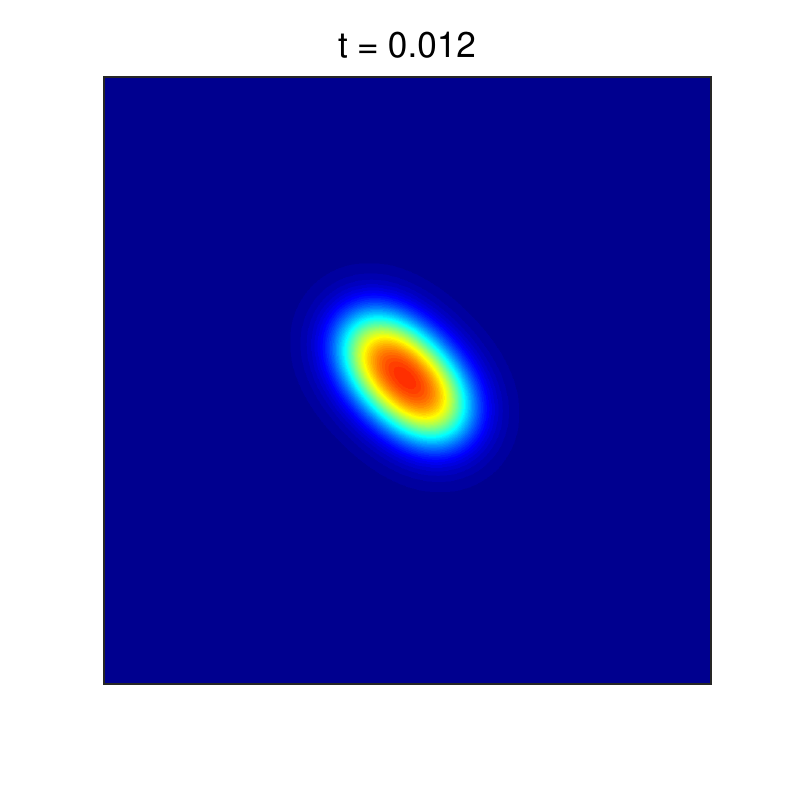}}
\vspace{-5mm}
\centerline{
\includegraphics[height=3.4cm,width=3.7cm]{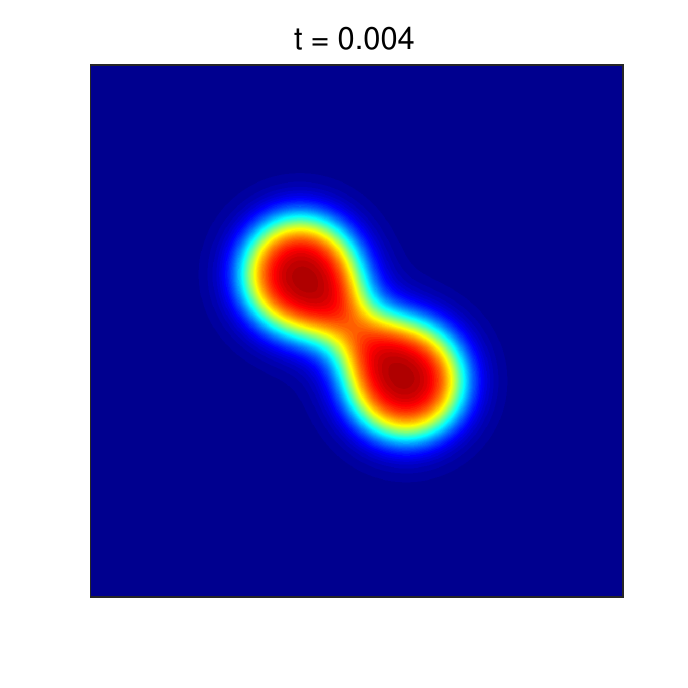}\hspace{-5mm}
\includegraphics[height=3.4cm,width=3.7cm]{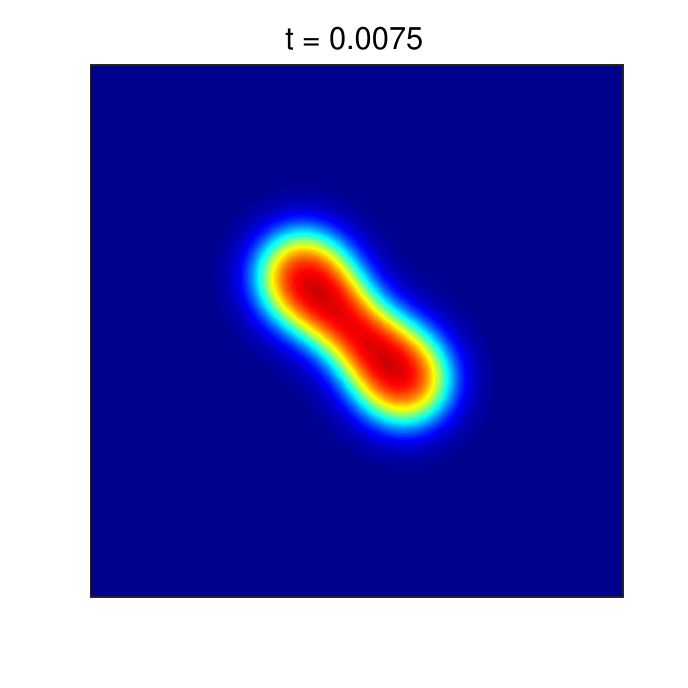}\hspace{-5mm}
\includegraphics[height=3.4cm,width=3.7cm]{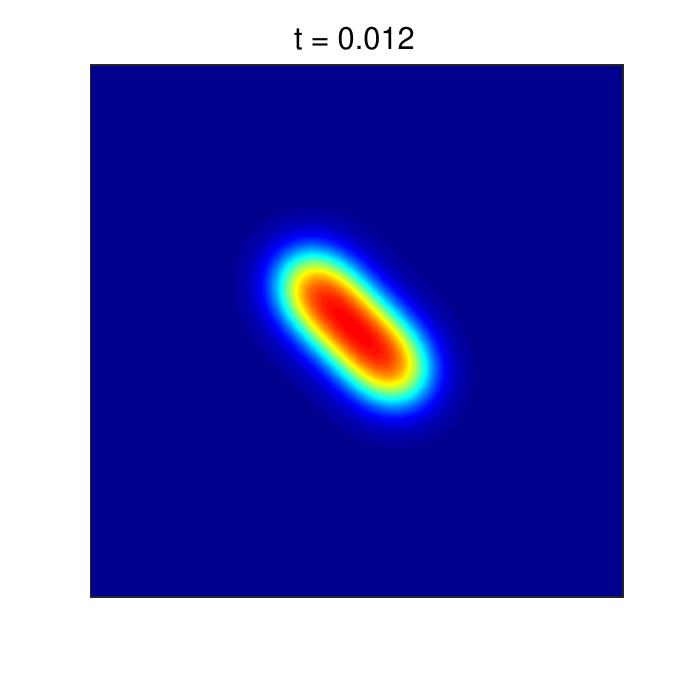}\hspace{-5mm}
\includegraphics[height=3.4cm,width=3.7cm]{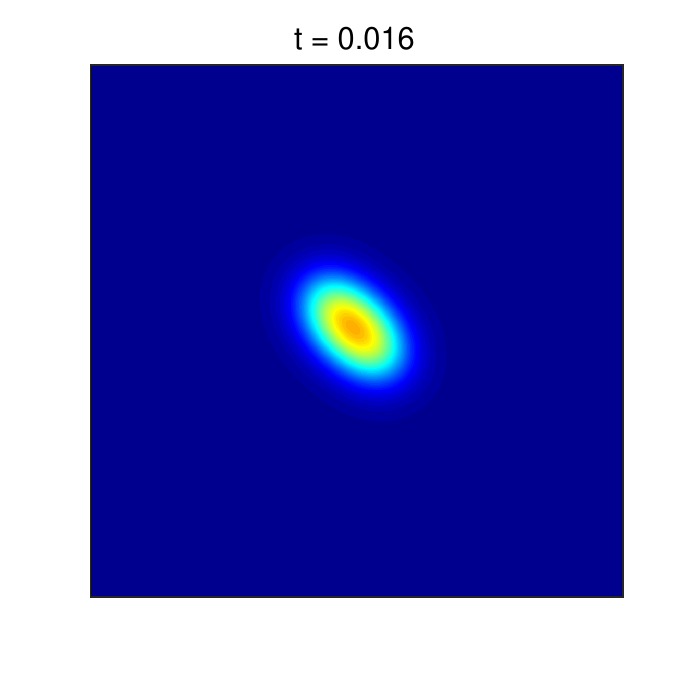}}
\vspace{-5mm}
\centerline{
\includegraphics[height=3.4cm,width=3.7cm]{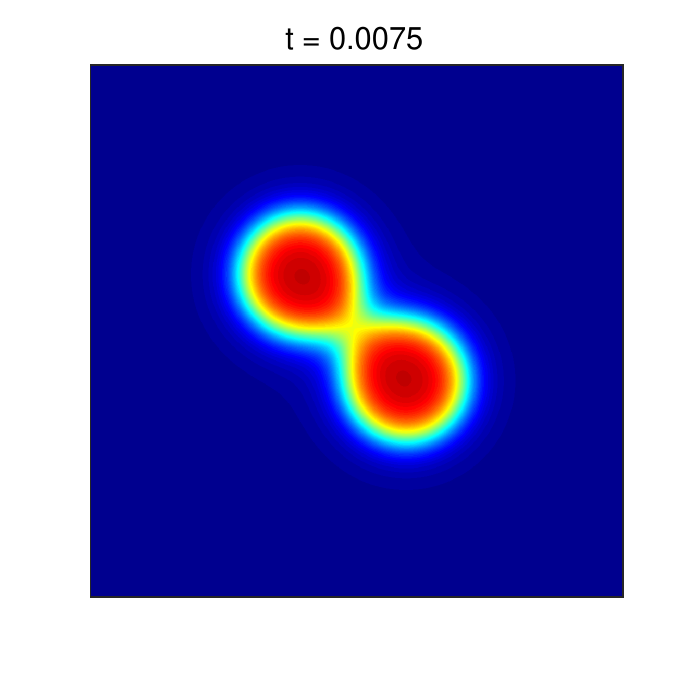}\hspace{-5mm}
\includegraphics[height=3.4cm,width=3.7cm]{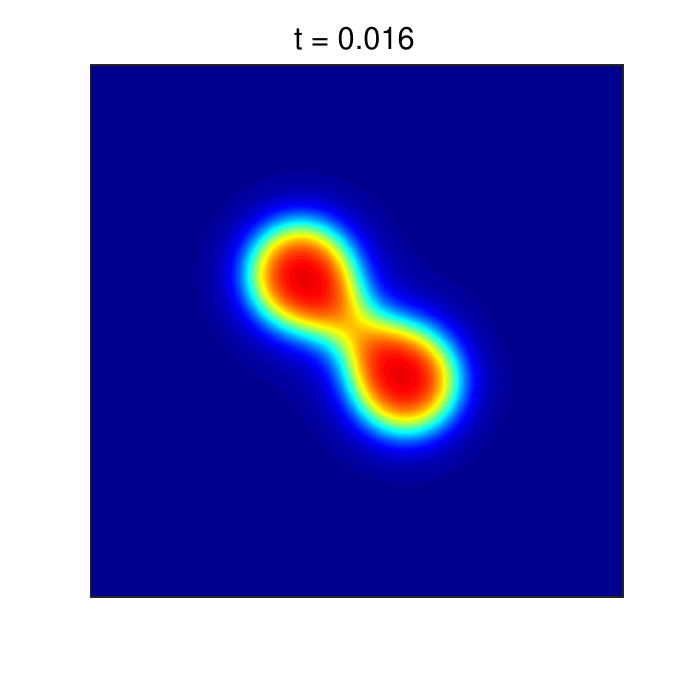}\hspace{-5mm}
\includegraphics[height=3.4cm,width=3.7cm]{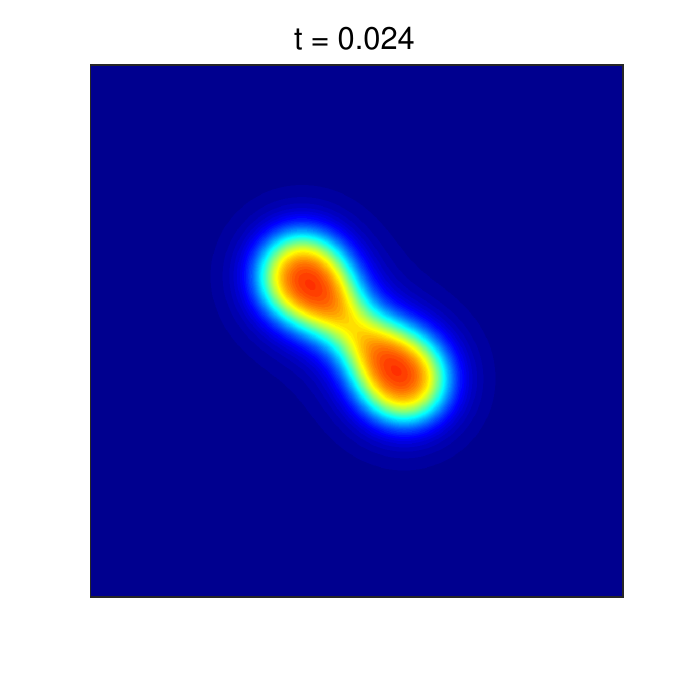}\hspace{-5mm}
\includegraphics[height=3.4cm,width=3.7cm]{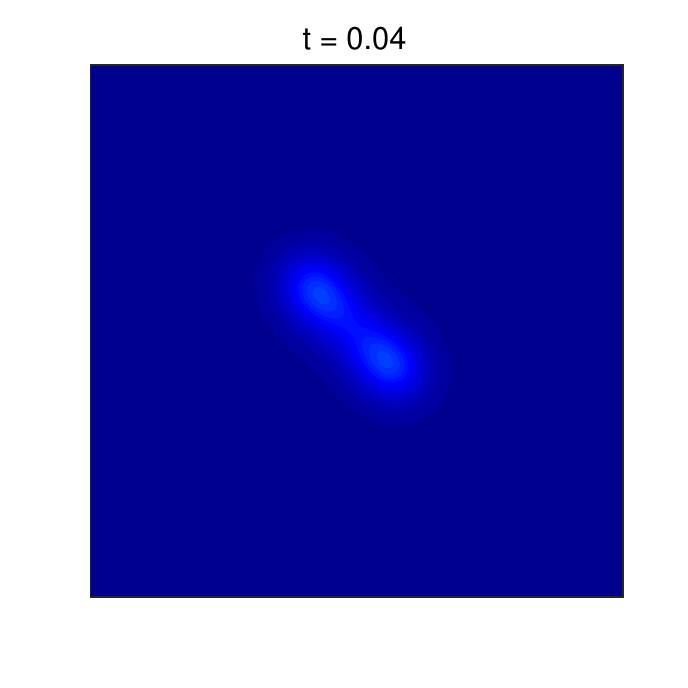}}
\vspace{-5mm}
\centerline{
\includegraphics[height=3.4cm,width=3.7cm]{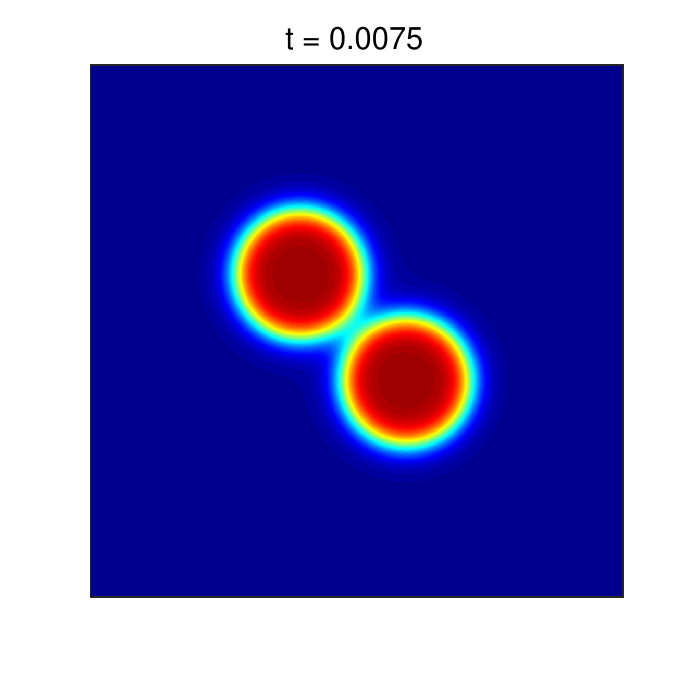}\hspace{-5mm}
\includegraphics[height=3.4cm,width=3.7cm]{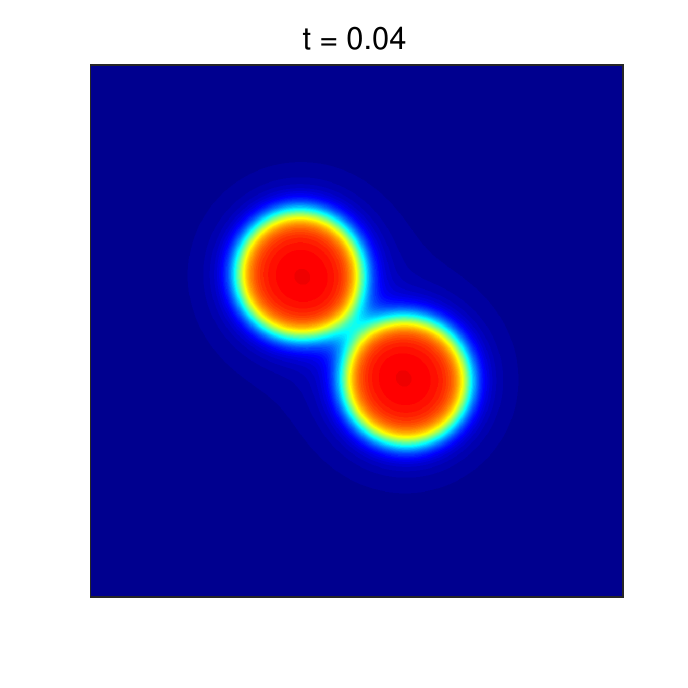}\hspace{-5mm}
\includegraphics[height=3.4cm,width=3.7cm]{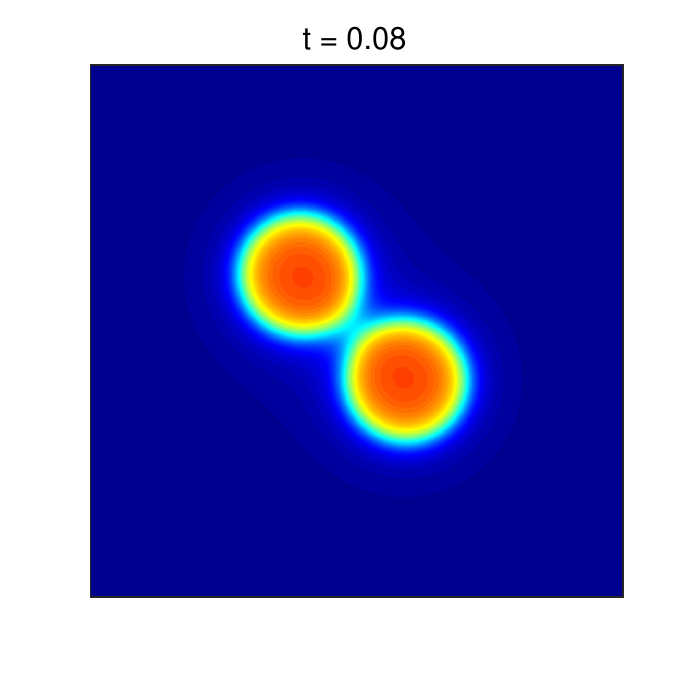}\hspace{-5mm}
\includegraphics[height=3.4cm,width=3.7cm]{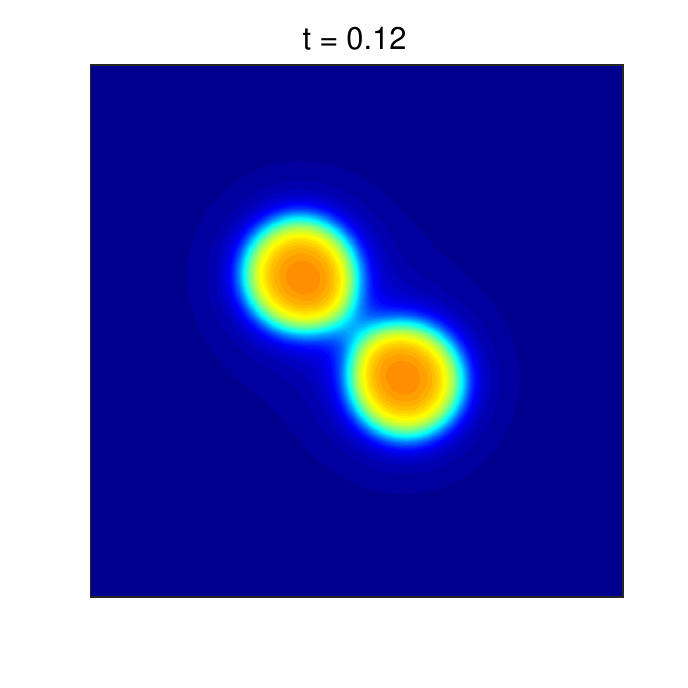}}
\vspace{-5mm}
\caption{Dynamics of the two kissing bubbles in the classical Allen--Cahn equation (top row) and its fractional counterparts with $\ap = 1.9$ (second row), $1.5$ (third row), and $0.7$ (last row).} \label{Fig3}
\end{figure}
%%%% %%%% %%%% %%%
Figure \ref{Fig3} shows the time evolution of the two bubbles in both classical (non-fractional) and fractional  Allen--Cahn equations with $\dt = 0.03$. 
In our simulations, we choose the mesh size $h = 1/1024$ and the time step $\tau = 0.0005$.   
Initially, the two bubbles are centered at $x_1 = (0.4, 0.4)$ and $x_2 = (0.6, 0.6)$, respectively.
In the classical case, the two bubbles first coalesce into one bubble, and then this newly formed bubble shrinks and are eventually absorbed by the fluid (see Fig. \ref{Fig3} top row). 
The dynamics in the fractional cases with $\ap > 1$ are similar to those in the classical cases, but the coalesceness and disappearing of the two bubbles become much slower -- the smaller the power $\ap$, the slower the dynamics. 
%%%% Figure 4
\begin{figure}[htb!]
\centerline{
\includegraphics[height=4.6cm,width=6.7cm]{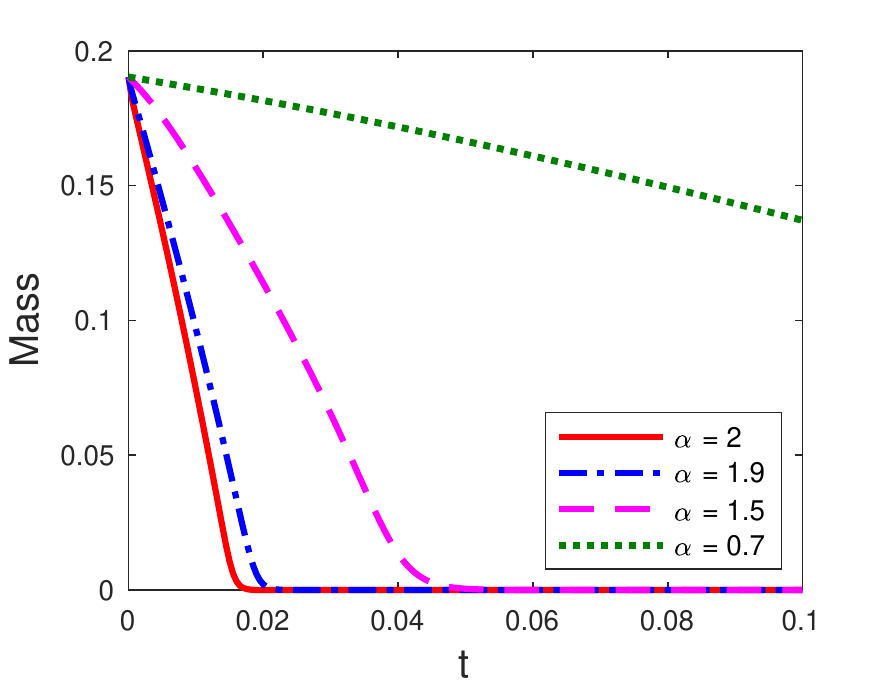}}
\caption{Time evolution of the mass $\int_{{\mathbb R}^2} |u| d{\bf x}$, where $\ap = 2$ represents the classical cases.} \label{Fig3-Mass}
\end{figure}
By contrast, if  $\ap < 1$  the two bubbles do not merge, and they finally vanish at the same time. 
It takes much longer time for two bubbles to vanish for $\ap < 1$ than that took for $\ap > 1$. 
The time evolution of the mass in Figure \ref{Fig3-Mass} further confirms these observations. 

%%%% %%%% %%%% %%%%
To demonstrate the efficiency of our methods, we present in Table \ref{Tab10} the CPU time that is used to solve the  2D fractional Allen--Cahn equations for various numbers of spatial unknowns $M = N_xN_y$, where the time step $\tau = {0.0005}$ is fixed. 
%%%% %%%% Table CPU
\begin{table}[htb!]
\begin{center}
\begin{tabular}{|l|c|c|c|c|}
\hline
CPU time (s) & $M = 1024^2$ & $M = 512^2$  & $M = 256^2$ & $M = 128^2$\\
\hline
$\ap = 1.9$ & 25728 & 4497 & 456 & 86\\
$\ap = 1.5$ & 4896 & 974 & 112 & 27\\
$\ap = 1$    & 648 & 170 & 30 & 9 \\
$\ap = 0.7$ & 430 & 107 & 24 & 7 \\
\hline
\end{tabular}
\caption{CPU time (in seconds) used in solving the 2D fractional Allen--Cahn equations till $t = 0.015$. }\label{Tab10}
\end{center}
\end{table}
It shows that 
%when doubling the number of points $N_x$ and $N_y$,  {\color{blue}the CPU time increases almost four times, since the computational cost is ${\mathcal O}(M\log M)$.}  
for fixed $M$, the CPU time takes longer for larger $\ap$, because the stiffness matrix from larger $\ap$ has bigger conditional number and it will affect the CG iterations. 
Furthermore, due to the symmetric block Toeplitz structure of the stiffness matrix, the storage memory requirement is ${\mathcal O}(M)$. 
Comparing to the 2D cases,  the computations of the 3D fractional Allen--Cahn problems  demand more computational costs and storage memory.
Thanks to the fast algorithms utilizing the {multilevel Toeplitz} structure, our method enables us to simulate the 3D problems efficiently. 
Figure \ref{Fig-AC3} shows the isosurface plots of the solution $u$ and their contour plots $u(x, y, 0)$, for different time $t$.
%%%% %%%% %%%% %%%%
\begin{figure}[htb!]
\centerline{
\includegraphics[height=4.60cm,width=3.90cm]{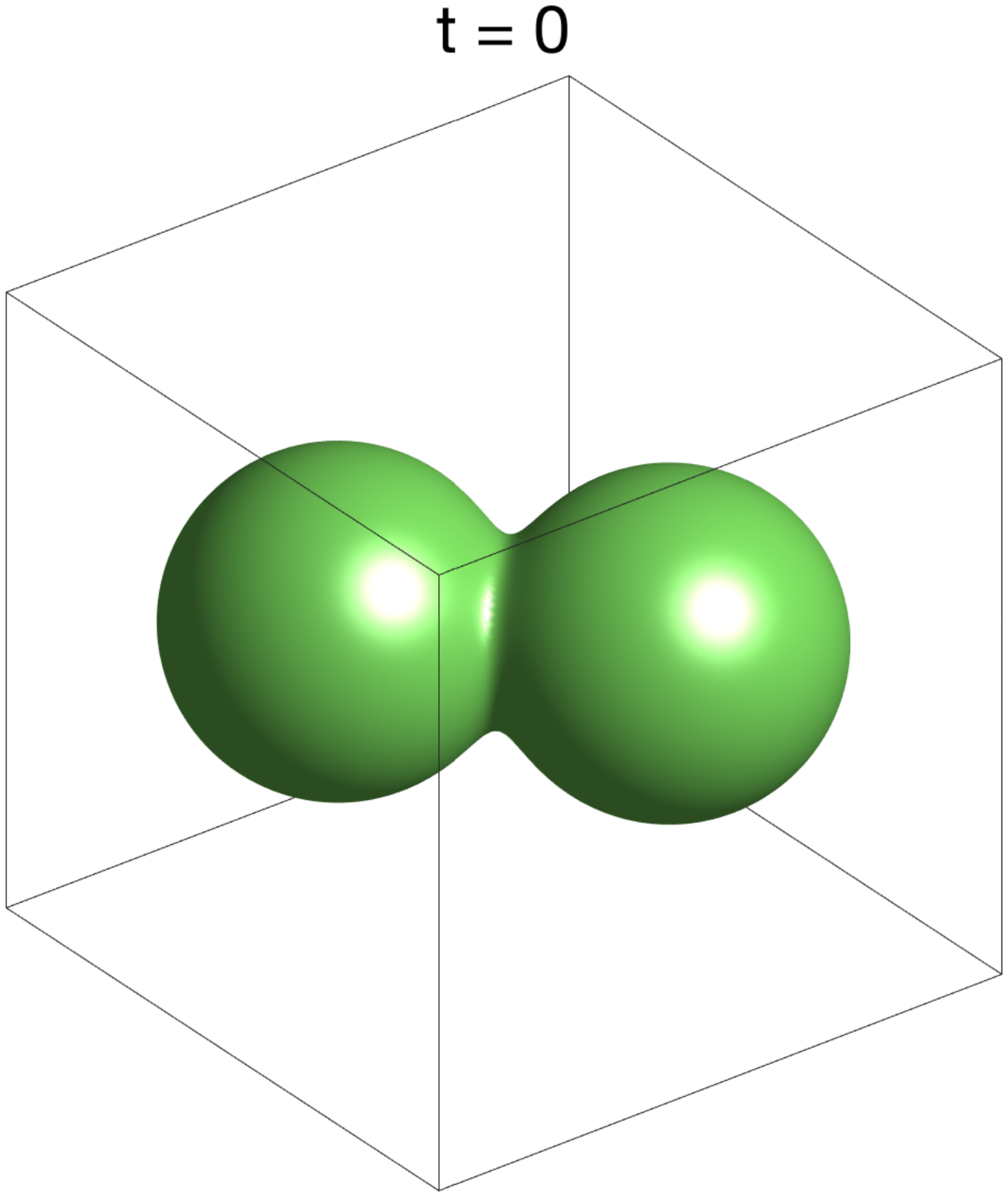}\hspace{-6mm}
\includegraphics[height=4.60cm,width=3.90cm]{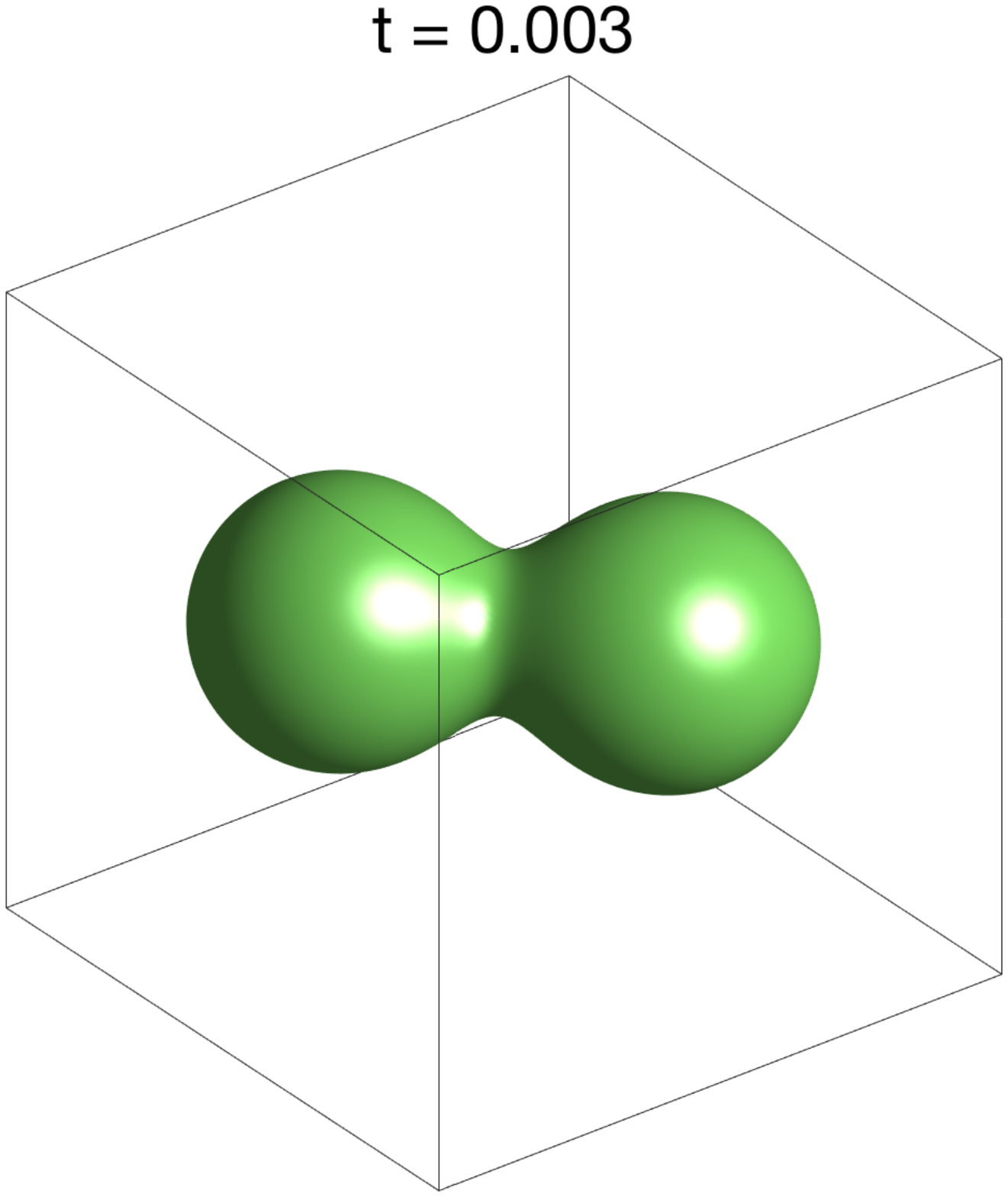}\hspace{-6mm}
\includegraphics[height=4.60cm,width=3.90cm]{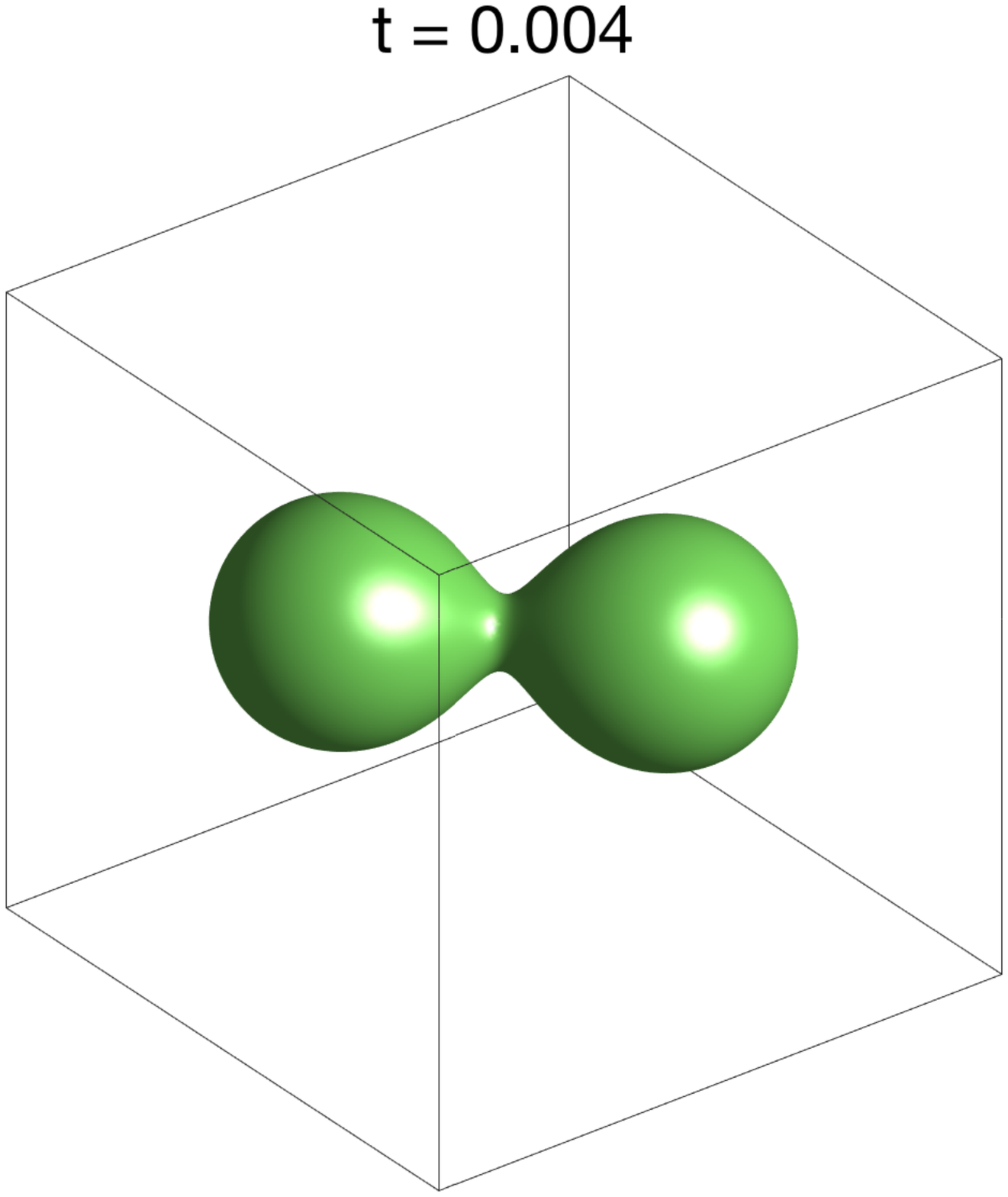}\hspace{-6mm}
\includegraphics[height=4.60cm,width=3.90cm]{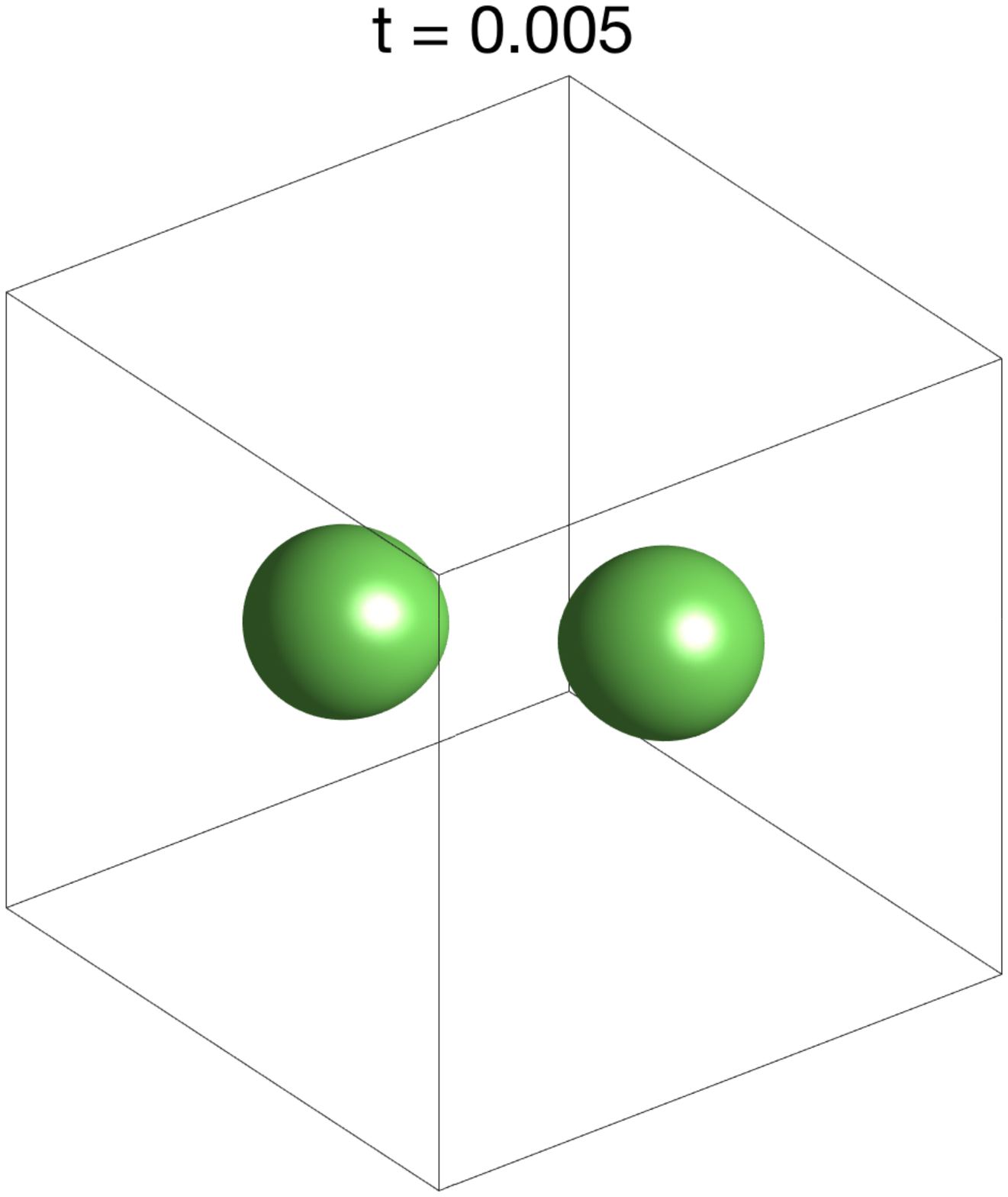}}
\vspace{-6mm}
\centerline{
\includegraphics[height=3.50cm,width=3.80cm]{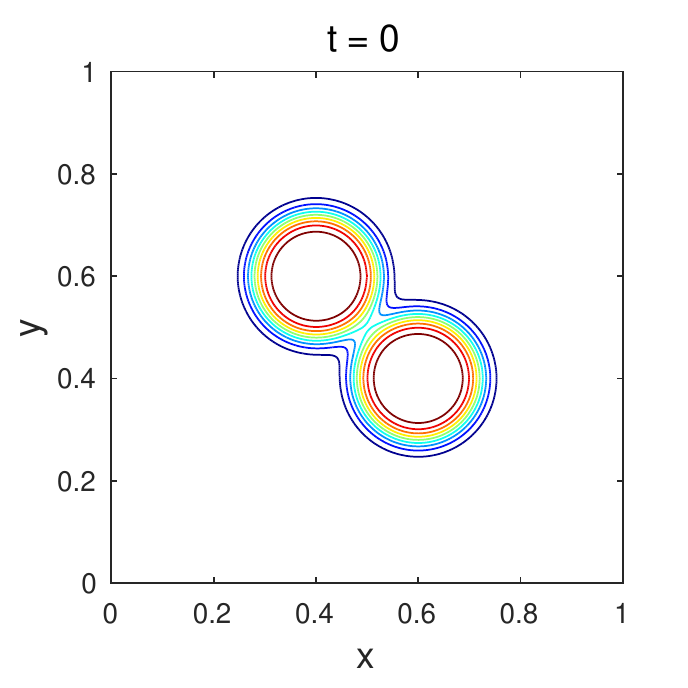}\hspace{-4mm}
\includegraphics[height=3.50cm,width=3.80cm]{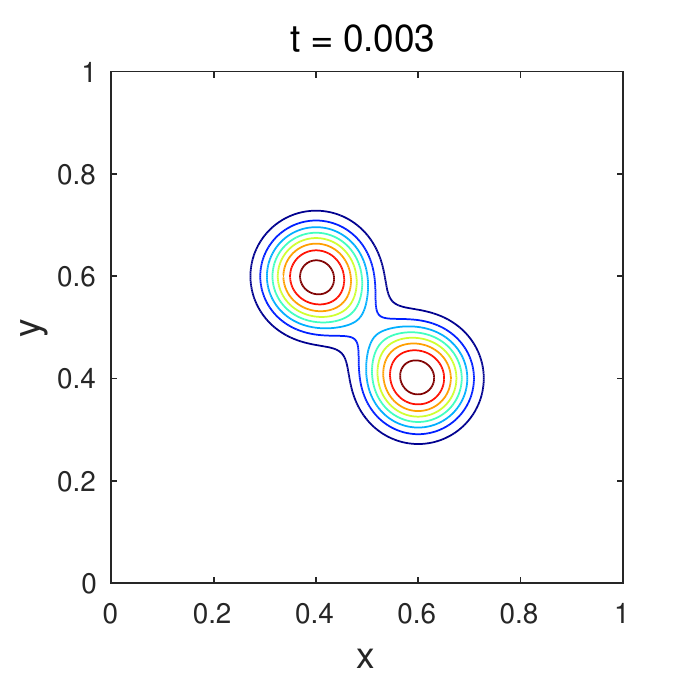}\hspace{-4mm}
\includegraphics[height=3.50cm,width=3.80cm]{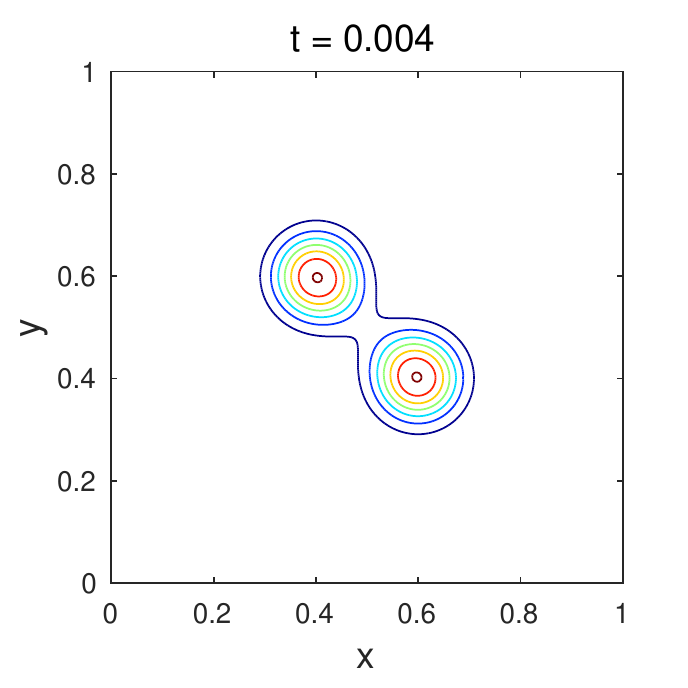}\hspace{-4mm}
\includegraphics[height=3.50cm,width=3.80cm]{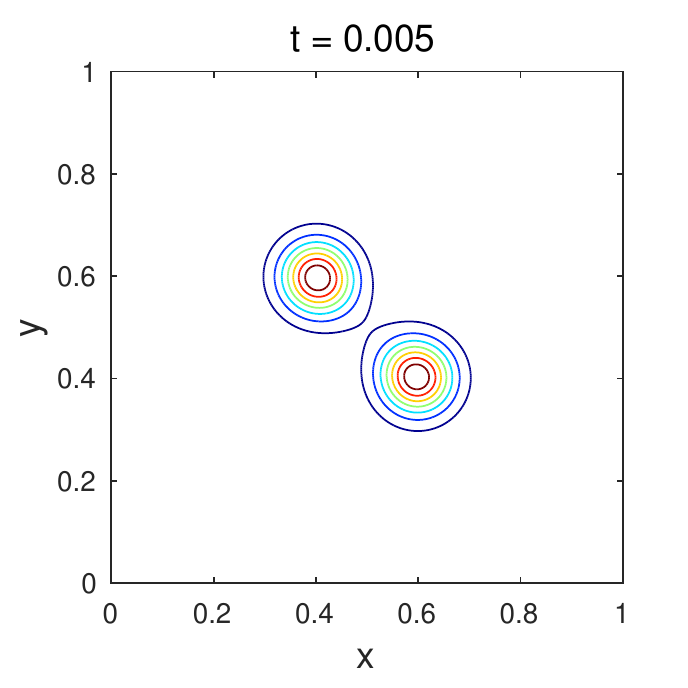}}
\caption{Dynamics of the two kissing bubbles in the 3D fractional Allen--Cahn equation with {$\dt = 0.03$} and $\alpha = 1.9$. Top row:  Isosurface plot for ({$u = 0.5$}); Bottom row: $u(x, y, z = 0)$.}\label{Fig-AC3}
\end{figure}
Here, the initial center of the two bubbles are $\bx_1 = {(0.4, 0.6, 0.5)}$ and $\bx_2 = {(0.6, 0.4, 0.5)}$, and the parameters $\dt = {0.03}$ and $\ap = 1.9$. 
It shows that the two initially kissing bubbles experience a process of merging, shrinking and finally disappearing, a similar process to the 2D cases observed in Fig. \ref{Fig3}.

% =====================================================================
\section{Conclusions}
\label{section8}
\setcounter{equation}{0}

%%%% %%%% %%%%
We proposed accurate and efficient finite difference methods to discretize the two- and three-dimensional integral fractional Laplacian $(-\Dt)^{\fl{\ap}{2}}$ and apply them to solve the fractional Poisson equation and fractional Allen--Cahn equation. 
The key idea of our method is to introduce a splitting parameter $\gm \in (\ap, {2]}$ and then reformulate the fractional Laplacian as the weighted integral of a central difference quotient,  so as to avoid directly approximating the hypersingular kernel. 
Both numerical analysis and simulations show that our methods have the accuracy of ${\mathcal O}(h^{\veps})$ for  $u \in C^{\lfloor\ap\rfloor,\,\ap-\lfloor\ap\rfloor+\veps}({\mathbb R}^d)$, while ${\mathcal O}(h^2)$ for  $u \in C^{2+\lfloor\ap\rfloor,\,\ap-\lfloor\ap\rfloor+\veps}({\mathbb R}^d)$,  with $0< \veps < 1+\lfloor\ap\rfloor-\ap$.  
To the best of our knowledge, the proposed methods are the first finite difference methods for the high-dimensional integral fractional Laplacian (\ref{fL-nD}).
Numerical studies on the fractional Poisson problem and fractional Allen--Cahn equation were presented to test the accuracy and efficiency of our methods. 
It shows that our finite difference methods have the second order of accuracy in solving the fractional Poisson problem, if the solution satisfy  $u \in C^{1,1}({\mathbb R}^d)$.
Numerical studies on the fractional Allen--Cahn equation suggest that our methods can be efficiently implemented via the fast Fourier transform to solve the fractional PDEs, which  significantly reduce the computational costs and storage memory requirements in practice. 

\vskip 10pt
{\bf Acknowledgements.} The authors thank Prof. Hans-Werner van Wyk for the helpful discussion. 
This work was supported by the US National Science Foundation under grant number DMS-1620465.
% =======================================================================
% References;
% =======================================================================

\end{document}